\documentclass[11pt,a4paper]{article}

\usepackage[margin=1in]{geometry}
\usepackage{palatino}
\usepackage{subcaption}

\usepackage{enumitem}

\usepackage[
  colorlinks,
  linkcolor = blue,
  citecolor = blue,
  urlcolor = blue]{hyperref}

\usepackage{amsmath,amssymb,amsthm}
\usepackage{mathtools}
\mathtoolsset{centercolon}

\usepackage[affil-it]{authblk}
\usepackage{tabu}

\newcommand{\N}{\mathbb{N}}
\newcommand{\C}{\mathbb{C}}
\newcommand{\Z}{\mathbb{Z}}

\newcommand{\x}{\otimes}

\newcommand{\be}{\begin{equation}}
\newcommand{\ee}{\end{equation}}

\renewcommand{\vec}[1]{\mathbf{#1}}

\DeclareMathOperator{\tr}{tr}
\DeclareMathOperator{\spn}{span}

\newcommand{\Thm}[1]{\hyperref[thm:#1]{Theorem~\ref*{thm:#1}}}
\newcommand{\Lem}[1]{\hyperref[lem:#1]{Lemma~\ref*{lem:#1}}}
\newcommand{\Cor}[1]{\hyperref[cor:#1]{Corollary~\ref*{cor:#1}}}
\newcommand{\Def}[1]{\hyperref[def:#1]{Definition~\ref*{def:#1}}}
\newcommand{\Obs}[1]{\hyperref[obs:#1]{Observation~\ref*{obs:#1}}}
\newcommand{\Prop}[1]{\hyperref[prop:#1]{Proposition~\ref*{prop:#1}}}
\newcommand{\Rem}[1]{\hyperref[rem:#1]{Remark~\ref*{rem:#1}}}
\newcommand{\Ex}[1]{\hyperref[ex:#1]{Example~\ref*{ex:#1}}}

\newcommand{\Sec}[1]{\hyperref[sec:#1]{Section~\ref*{sec:#1}}}

\newcommand{\Fig}[1]{\hyperref[fig:#1]{Figure~\ref*{fig:#1}}}
\newcommand{\Tab}[1]{\hyperref[tab:#1]{Table~\ref*{tab:#1}}}
\newcommand{\EqRef}[1]{\hyperref[eq:#1]{(\ref*{eq:#1})}}
\newcommand{\Eq}[1]{Equation~\hyperref[eq:#1]{(\ref*{eq:#1})}}

\makeatletter
\newtheorem*{rep@theorem}{\rep@title}
\newcommand{\newreptheorem}[2]{%
\newenvironment{rep#1}[1]{%
 \def\rep@title{#2 \ref{##1}}%
 \begin{rep@theorem}}%
 {\end{rep@theorem}}}
\makeatother

\newreptheorem{theorem}{Theorem}
\newreptheorem{lemma}{Lemma}

\newtheorem{theorem}{Theorem}[section]
\newtheorem*{theorem*}{Theorem}
\newtheorem{lemma}[theorem]{Lemma}
\newtheorem{cor}[theorem]{Corollary}

\newtheorem{prop}[theorem]{Proposition}

\theoremstyle{definition}
\newtheorem{definition}[theorem]{Definition}
\newtheorem{remark}[theorem]{Remark}

\DeclareMathOperator{\aut}{Aut}
\DeclareMathOperator{\qut}{Qut}

\newcommand{\qeds}{\qed\vspace{.2cm}}

\newcommand{\one}{\ensuremath{\vec{1}}}

\DeclareMathOperator{\rel}{rel}

\DeclareMathOperator{\conv}{conv}

\title{Quantum symmetry vs nonlocal symmetry}

\author[1]{David E.~Roberson\footnote{davideroberson@gmail.com}} 
\author[2]{Simon Schmidt\footnote{Simon.Schmidt@glasgow.ac.uk}}

\affil[1]{Department of Applied Mathematics and Computer Science, Technical University of Denmark, DK-2800 Lyngby, Denmark}
\affil[2]{School of Mathematics and Statistics, University of Glasgow, University Place, Glasgow G12 8QQ, United Kingdom}

\begin{document}

\maketitle

\begin{abstract}

We introduce the notion of nonlocal symmetry of a graph $G$, defined as winning quantum correlation for the $G$-automorphism game that cannot be produced classically, i.e., by a local hidden variable theory. Recent connections made between quantum group theory and quantum information show that quantum correlations for this game correspond to tracial states on $C(\qut(G))$ -- the algebra of functions on the quantum automorphism group of $G$. This allows us to also define nonlocal symmetry for any quantum permutation group. We investigate the differences and similarities between this and the already established notion of quantum symmetry, defined as non-commutativity of $C(\qut(G))$. Roughly speaking, quantum symmetry can be viewed as non-classicality of our \emph{model} of reality, whereas nonlocal symmetry is non-classicality of our \emph{observation} of reality respectively.


We show that quantum symmetry is a necessary but not sufficient condition for nonlocal symmetry. In particular, we show that the complete graph on five vertices is the only connected graph on five or fewer vertices with nonlocal symmetry, despite there being a dozen others with quantum symmetry. In particular this shows that the quantum symmetric group on four points, $S_4^+$, does not exhibit nonlocal symmetry, answering a question from the literature. In contrast to quantum symmetry, we show that two \emph{disjoint} classical automorphisms does not guarantee nonlocal symmetry. However, three disjoint automorphisms does suffice. We also give a construction of quantum permutation matrices built from a finite abelian group $\Gamma$ and a permutation $\pi$ on $|\Gamma|$ elements. Our computational evidence suggests that for cyclic groups of increasing size almost all permutations $\pi$ result in quantum permutation matrices that produce nonlocal symmetry. Nonlocal symmetry seems to occur far less often for non-cyclic groups; the extreme case being $\mathbb{Z}^3_2$ for which it never occurs. 
We also investigate under what conditions nonlocal symmetry arises when taking unions or products of graphs/groups. In particular we show that the tensor product of two quantum groups has no nonlocal symmetry if neither factor has it and at least one factor is a classical group. 

\end{abstract}




\section{Introduction}

Recently, researchers have uncovered links between the abstract theory of quantum groups from non-commutative mathematics and entanglement-assisted strategies for nonlocal games, an operationally defined notion from quantum information. Specifically, the works~\cite{qfuncs} and~\cite{qperms} established strong connections between quantum automorphism groups of graphs~\cite{banicahomogeneous} and quantum strategies for the graph isomorphism game~\cite{qiso1}. This work investigates the relationship between the notions of quantumness that have arisen in these two scenarios. Within quantum group theory, a graph is said to have \emph{quantum symmetry} if the algebra defining its quantum automorphism group is non-commutative. On the other hand, an entanglement-assisted strategy for a nonlocal game is said to be \emph{nonlocal} if the correlation (joint conditional probability distribution) produced by this strategy cannot be obtained by a so-called \emph{local} (i.e., classical) strategy. We introduce the term \emph{nonlocal symmetry} to refer to when a graph admits a nonlocal entanglement-assisted strategy for its automorphism game. Both notions are natural in their respective settings: In the physically motivated setting of nonlocal games, it is useful to characterize quantumness from the perspective of a classical observer, who may not trust that the parties are performing any genuinely quantum actions and would only be able to observe 
the resulting correlation. Whereas in the abstract world of quantum groups, one takes a purely mathematical point of view with no reference to ``observers". Thus one has only algebraic considerations. Despite arising in seemingly different contexts, the connections uncovered in~\cite{qperms,qfuncs} allow us to directly compare these two types of quantumness. In fact, it was already noted (without proof) in~\cite{qperms} that these two notions do not coincide. The most succinct way to describe the difference between these notions is that quantum symmetry is a quantumness of our \emph{model} of reality whereas nonlocal symmetry is a quantumness of our \emph{observation} of reality. This is analogous to the relationship between the concepts of entanglement and nonlocality in quantum mechanics.

In this work we make the first serious investigation into the difference between quantum and nonlocal symmetry, and as well study nonlocal correlations of graph automorphism games more generally. Answering a question raised in~\cite{qperms}, we show that the complete graph on four vertices does not admit nonlocal symmetry (Section~\ref{sec:K4}). This implies that no graph on four vertices admits nonlocal symmetry, despite several of them having quantum symmetry. In contrast, we show that the complete graph on five vertices does admit nonlocal symmetry by exhibiting a specific non-classical correlation for its automorphism game (Section~\ref{sec:K5}). Later, in Section~\ref{sec:smallgraphs}, we show that this (and its complement the empty graph) is the only graph on five or fewer vertices with nonlocal symmetry. In Section~\ref{sec:construction} we present two constructions of quantum permutations. For the first, whose ingredients are an abelian group $\Gamma$ and permutation $\pi$ of the dual group $\widehat{\Gamma}$, we show that the produced quantum permutation is non-commutative (i.e., is a genuine quantum symmetry) precisely when $\pi(x) = z\sigma(x)$ for some $z \in \widehat{\Gamma}$ and $\sigma \in \aut(\widehat{\Gamma})$. 
We also give examples of this construction producing non-classical correlations (in particular, our example for the complete graph on five vertices is constructed in this way), but we leave open the question of characterizing when this occurs. We compute all such correlations constructed from abelian groups of order at most 10, showing that most of them do produce non-classical correlations. The second construction we present requires the graph in question to have three ``disjoint" classical automorphisms, but always produces a nonlocal symmetry. A similar construction for quantum symmetry requires only two disjoint automorphisms~\cite{foldedcubes}, but this does not suffice for nonlocal symmetry. Section~\ref{sec:products} investigates quantum and nonlocal symmetries in products, both of graphs and quantum groups. We show that even though the free product of two non-trivial classical groups is always a non-classical quantum group, it does not produce any nonlocal symmetry. 
Next we show that the free wreath product of a classical group with $\mathbb{Z}_2$ does not produce nonlocal symmetry. This and the previous result shows that if graphs $G$ and $H$ have no quantum symmetry, then their disjoint union does not admit nonlocal symmetry unless they are $G$ and $H$ are quantum isomorphic but not isomorphic. We also show that if $\mathbb{G}$ is a quantum group without nonlocal symmetry and $\mathbb{H}$ is a classical group, then the tensor product $\mathbb{G} \times \mathbb{H}$ does not have nonlocal symmetry. Applying this to graphs, we prove that if graphs $G$ and $H$ have no nonlocal symmetry and no quantum symmetry respectively, and their spectrums satisfy a certain property, then their cartesian/categorical product does not admit nonlocal symmetry either. We end with a discussion of our results and open questions.

\section{Background}\label{sec:background}

We make use of basic notions and terminology from graph theory. For us, a graph $G$ is always finite and simple, meaning no multiple edges or loops, and thus consists of a vertex set $V(G)$ and edge set $E(G)$ which is a set of unordered pairs of its vertices. We write $i \sim_G j$ to denote that vertices $i$ and $j$ are adjacent in $G$, and we often simply write $i \sim j$ when $G$ is clear from context. The \emph{adjacency matrix} of a graph $G$, denoted $A_G$, is the symmetric 01-matrix whose $ij$-entry is 1 if and only if $i \sim j$. An \emph{isomorphism} $\varphi$ of graphs $G$ and $H$ is a bijection from $V(G)$ to $V(H)$ such that $\varphi(i) \sim_H \varphi(j)$ if and only if $i \sim_G j$. An automorphism of $G$ is an isomorphism from $G$ to itself, and these form a group under composition known as the \emph{automorphism group} of $G$, denoted $\aut(G)$. The \emph{complement} of a graph $G$, denoted $\overline{G}$, is the graph with the same vertex set as $G$ whose edges are the non-edges of $G$.

\subsection{Quantum groups}\label{sec:qgroups}

In the following, we give a brief introduction to compact matrix quantum groups which were introduced by Woronowicz in \cite{woronowicz}. Note that we denote by $A\otimes B$ the minimal tensor product of the $C^*$-algebras $A$ and $B$.

A \emph{compact matrix quantum group} $\mathbb{G}$ is a pair $(C(\mathbb{G}),U)$, where $C(\mathbb{G})$ is a unital $C^*$-algebra and  $U = (u_{ij}) \in M_n(C(\mathbb{G}))$ is a matrix such that 
\begin{itemize}
    \item the elements $u_{ij}$, $1 \le i,j \le n$, generate $C(\mathbb{G})$, 
    \item the $*$-homomorphism $\Delta: C(\mathbb{G}) \to C(\mathbb{G}) \otimes C(\mathbb{G})$, $u_{ij} \mapsto \sum_{k=1}^n u_{ik} \otimes u_{kj}$ exists,
    \item the matrix $U$ and its transpose $U^{T}$ are invertible. 
\end{itemize}
The matrix $U$ is usually called \emph{fundamental representation} of $\mathbb{G}$. Every compact matrix group $\Gamma$ gives rise to a compact matrix quantum group: The pair consisting of the $C^*$-algebra $C(\Gamma)$ of continuous complex-valued functions on $\Gamma$ and the matrix $U$ with entries $u_{ij}: \Gamma \to \C$, $g \mapsto g_{ij}$ fulfills the conditions above. Furthermore, if $\mathbb{G}=(C(\mathbb{G}),U)$ is a compact matrix quantum group where $C(\mathbb{G})$ is a commutative $C^*$-algebra, then there exists a compact matrix group $\Gamma$ such that $(C(\Gamma),\tilde{U}) \cong (C(\mathbb{G}), U)$. Here $\tilde{U}=(\tilde{u}_{ij})$ is the matrix with entries $\tilde{u}_{ij}: \Gamma \to \C$, $g \mapsto g_{ij}$. Thus, we have a one-to-one correspondence of compact matrix groups and compact matrix quantum groups with commutative $C^*$-algebras. 


The following important example is due to Wang \cite{wang}. We denote by $C^*(G \,|\,R)$ the universal $C^*$-algebra with generators $G$ and relations $R$. For the definition of universal $C^*$-algebras, see for example \cite[Section II.8.3]{Blackadar}.

\begin{definition}
The \emph{quantum symmetric group} $S_n^+= (C(S_n^+),U)$ is the compact matrix quantum group, where
\begin{align*}
C(S_n^+) := C^*(u_{ij}, \, 1 \le i,j \le n \, | \, u_{ij} = u_{ij}^* = u_{ij}^2, \, \sum_{k=1}^n u_{ik} = \sum_{k=1}^n u_{ki} =\vec{1}). 
\end{align*} 
\end{definition}

We have the following connection to the symmetric group $S_n$: If we add the relations $u_{ij}u_{kl}=u_{kl}u_{ij}$ for all $1\le i,j,k,l \le n$ to the relations of $C(S_n^+)$, we obtain the algebra of complex-valued, continuous functions on $S_n$, i.e.,
\begin{align*}
C(S_n) \cong C^*(u_{ij}, \, 1 \le i,j \le n \, | \, u_{ij} = u_{ij}^* = u_{ij}^2, \, \sum_{k=1}^n u_{ik} = \sum_{k=1}^n u_{ki} =\vec{1},u_{ij}u_{kl}=u_{kl}u_{ij}). 
\end{align*}
A matrix $U=(u_{ij})_{1 \le i,j \le n}$, where $u_{ij}$ are elements of some unital $C^*$-algebra is called a \emph{magic unitary} or \emph{quantum permutation} if $u_{ij} = u_{ij}^* = u_{ij}^2$ for all $i,j$ and $\sum_{k=1}^n u_{ik} = \sum_{k=1}^n u_{ki} =\vec{1}$ for all $i$. Thus, $C(S_n^+)$ is the universal $C^*$-algebra generated by the entries of a magic unitary $U=(u_{ij})_{1 \le i,j \le n}$.

Quantum subgroups of $S_n^+$ are called \emph{quantum permutation groups}. The next definition yields many examples of quantum subgroups of $S_n^+$.
 
\begin{definition}[\cite{banicahomogeneous}]
Let $G$ be a finite graph. The \emph{quantum automorphism group} $\qut(G)$ is the compact matrix quantum group $(C(\qut(G)), U)$, where $C(\qut(G))$ is the universal $C^*$-algebra with generators $u_{ij}$, $i,j \in V(G)$ and relations
\begin{align}\allowdisplaybreaks
&u_{ij} = u_{ij}^*= u_{ij}^2, && i,j \in V(G),\label{QA1}\\ 
&\sum_{l=1}^n u_{il} = \vec{1} = \sum_{l=1}^n u_{li}, && i \in V(G),\label{QA2}\\
&U A_G = A_G U \label{QA3},
\end{align}
where \eqref{QA3} is nothing but $\sum_ku_{ik}(A_G)_{kj}=\sum_k(A_G)_{ik}u_{kj}$ for all $i,j \in V(G)$.
\end{definition}

Adding the commutation relations $u_{ij}u_{kl}=u_{kl}u_{ij}$ for all $i,j,k,l  \in V(G)$ now gives us the algebra of complex-valued, continuous functions over the automorphism group $\aut(G)$ of the graph $G$. Sometimes it is also the case that the relations \eqref{QA1}--\eqref{QA3} already imply that the $C^*$-algebra $C(\qut(G))$ is commutative. Then, we say that the graph $G$ has \emph{no quantum symmetry}. Otherwise, if $C(\qut(G))$ is non-commutative, we say the graph $G$ \emph{does have quantum symmetry}.

\subsection{Quantum information and the isomorphism game}\label{sec:qinfo}

In quantum mechanics, two or more systems are \emph{entangled} if their collective state cannot be described purely in terms of the states of the individual constituents. This is the phenomenon that enables Einstein's ``spooky action-at-a-distance". In quantum information, \emph{nonlocal games} provide a rigorous framework in which to study the capabilities and limits of entanglement. In a (2-player) nonlocal game, a referee/verifier sends a question/input to each of two players (Alice and Bob) and they must each respond to with an answer/output. The players then win or lose depending on a function of their questions and answers. This function, as well as the sets of allowable inputs and outputs, are known to the players who are allowed to devise a strategy ahead of time. Crucially though, the players are not allowed to communicate during the game, i.e., after they have received their inputs. The aim of the players is to win with as great a probability as possible\footnote{This success probability depends on the probability distribution according to which the verifier selects pairs of inputs, which is also known to the players. But we will only be interested in whether or not the game can be won perfectly and so we may assume a uniform distribution, or any distribution with full support.}. The power of entanglement is measured by contrasting the performance of classical players with players that are able to make local quantum measurements on a shared entangled state. In the extreme case, quantum players may be able to win with probability one even though no classical strategy achieves this. This phenomenon is usually known as \emph{pseudotelepathy}, as it appears as if the players are telepathic to one unfamiliar with quantum mechanics. The complete details of what constitutes a classical or quantum strategy will be given in the specific context of the isomorphism game below.

For graphs $G$ and $H$ with $|V(G)| = |V(H)|$ the $(G,H)$-isomorphism game~\cite{qiso1} is defined as follows. The referee sends Alice and Bob vertices $g_A$ and $g_B$ of $G$, and they must respond with some vertices $h_A$ and $h_B$ of $H$ respectively. They win if these vertices satisfy $\rel(g_A,g_B) = \rel(h_A,h_B)$, where $\rel$ is a function denoting the \emph{relationship} of two vertices, i.e., whether they are equal, distinct and adjacent, or distinct and non-adjacent. We remark that the game remains the same if one replaces $G$ and $H$ with their complements, since this does not change the truth value of $\rel(g_A,g_B) = \rel(h_A,h_B)$. Also note that if both $G$ and $H$ are complete (or empty) graphs then $\rel(g_A,g_B) = \rel(h_A,h_B)$ reduces to $\delta_{g_Ag_B} = \delta_{h_Ah_B}$. Thus any winning strategy for the $(G,H)$-isomorphism game is also a winning strategy for the $(K_{|V(G)|},K_{|V(H)|})$-isomorphism game, but not vice versa.

The idea behind the game is that the players are attempting to convince the referee that they know an isomorphism between $G$ and $H$. They are only required to play one round of the game, but they aim to optimize their probability of winning, and in fact we are only interested in whether they can win with probability one, in which case we say that they are playing with a \emph{perfect} or \emph{winning} strategy. It is not difficult to see that if $\varphi: V(G) \to V(H)$ is an isomorphism, then responding with $\varphi(g)$ upon receiving $g \in V(G)$ guarantees that the players win regardless of inputs. Conversely, in a \emph{deterministic classical strategy} Alice's output depends only on her input, and similarly for Bob. Thus such a strategy is described by a pair of functions $\varphi_A, \varphi_B : V(G) \to V(H)$ for Alice and Bob respectively. It is straightforward to see that if this is a perfect strategy then we must have that $\varphi_A = \varphi_B$ and these functions are an isomorphism from $G$ to $H$. More generally, a classical strategy may also include a source of randomness shared between and Alice and Bob, but this only allows them to pick a common deterministic classical strategy at random. Thus classical players can win the $(G,H)$-isomorphism game if and only if $G$ and $H$ are isomorphic.

\begin{remark}\label{rem:biggergame}
The definition of the isomorphism game given here differs from the original definition given in~\cite{qiso1}. There, the players may receive and respond with vertices of either $G$ or $H$. Under this definition, the condition $|V(G)| = |V(H)|$ is implied by the existence of any winning quantum (or even non-signalling) strategy and thus need not be explicitly included. However, with this assumption in place the quantum (and thus classical) winning strategies/correlations for the version of the game presented here are simply the quantum winning strategies/correlations of the ``larger" game restricted to inputs from one graph and outputs from the other. Further, the $|V(G)| = |V(H)|$ condition is necessary for the smaller game since otherwise any isomorphism from $G$ to an induced subgraph of $H$ would provide a winning classical strategy. Lastly, we note that these two versions of the game are not equivalent for non-signalling strategies. In fact the version of the $(G,H)$-isomorphism game presented here can be won by a non-signalling strategy whenever $H$ contains a regular induced subgraph that is neither empty nor complete.
\end{remark}

Given any strategy for the $(G,H)$-isomorphism game, there is a corresponding joint conditional probability distribution $p : V(H) \times V(H) \times V(G) \times V(G) \to [0,1]$ known as a \emph{correlation}. The value $p(h,h'|g,g')$ is the probability of Alice and Bob responding with answers $h$ and $h'$ conditioned on them receiving questions $g$ and $g'$ respectively. Note that the correlation determines whether the strategy is perfect: it is equivalent to $p(h,h'|g,g')$ being zero whenever answering $h$ and $h'$ to $g$ and $g'$ causes the players to lose.

For a deterministic classical strategy, the correlation is 01-valued. Indeed $p(h,h'|g,g')$ is 1 if and only if $\varphi_A(g) = h$ and $\varphi_B(g') = h'$ for the functions $\varphi_A$ and $\varphi_B$ described above. Allowing shared randomness, i.e., considering arbitrary classical strategies, results in a correlation that is a convex combination of these deterministic correlations. Thus the set of (perfect) classical correlations is the convex hull of the (perfect) deterministic classical correlations.

A quantum strategy for the $(G,H)$-isomorphism game consists of a unit vector $\psi$ from a Hilbert space $\mathcal{H}$, a quantum measurement $\mathcal{M}_g = \{M_{gh} \in B(\mathcal{H}) : h \in V(H)\}$ for each $g \in V(G)$ for Alice, and a quantum measurement $\mathcal{N}_g = \{N_{gh} \in B(\mathcal{H}) : h \in V(H)\}$ for each $g \in V(G)$ for Bob. Here $B(\mathcal{H})$ is the set of bounded linear operators on $\mathcal{H}$, and a quantum measurement is a collection of positive semidefinite operators that sum to the identity. For a quantum strategy it is also required that $M_{gh}N_{g'h'} = N_{g'h'}M_{gh}$ for all $g,g' \in V(G)$ and $h,h' \in V(H)$. Given such a quantum strategy, the corresponding correlation is
\begin{equation}\label{eq:qcorr}
    p(h,h'|g,g') = \psi^*M_{gh}N_{g'h'}\psi.
\end{equation}
It is well known that the set of quantum correlations for a nonlocal game is a closed convex set.

We remark that the above is often referred as a \emph{quantum commuting} strategy. In contrast, in quantum \emph{tensor} strategies the Hilbert space $\mathcal{H}$ is a tensor product $\mathcal{H}_A \otimes \mathcal{H}_B$ and Alice's measurement operators are from $B(\mathcal{H}_A)$ and analogously for Bob. Then the product of $M_{gh}$ and $N_{g'h'}$ in Equation~\ref{eq:qcorr} is replaced by a tensor product, and we no longer require that $M_{gh}$ and $N_{g'h'}$ commute. Note that a quantum commuting strategy can always be obtained from a quantum tensor strategy by replacing $M_{gh}$ and $N_{g'h'}$ with $M_{gh} \otimes \vec{1}$ and $\vec{1} \otimes N_{g'h'}$ respectively. Thus quantum commuting strategies are potentially more general, and in fact it is known that they are strictly more general for isomorphism games~\cite{qiso1}. One can additionally consider both finite and infinite-dimensional strategies, thus yielding four different types of quantum strategies. However, it is known that for the isomorphism game the infinite-dimensional quantum commuting strategies are the most general and the other three are all equivalent. We remark that in most works ``quantum strategies" refer to quantum tensor strategies, however we will use it to refer to quantum commuting strategies, as these are our main focus. However, all of the constructions we present (such as in Section~\ref{sec:construction}) will use finite dimensional operators and therefore can also be used in the tensor product framework.

The following lemma combines results of~\cite{qiso1} and~\cite{paulsen2016estimating}, characterizing the winning quantum correlations for the isomorphism game.

\begin{lemma}\label{lem:qisocorr}
Let $G$ and $H$ be graphs with $|V(G)| = |V(H)|$. Then $p : V(H)^2 \times V(G)^2 \to [0,1]$ is a winning quantum correlation for the $(G,H)$-isomorphism game if and only if there exists a $C^*$-algebra $\mathcal{A}$ with tracial state $\tau$ and a $V(G) \times V(H)$ quantum permutation matrix $U = (u_{gh}) \in M_n(\mathcal{A})$ satisfying $A_GU = UA_H$ such that
\begin{equation}\label{eq:qisocorr}
    p(h,h'|g,g') = \tau(u_{gh}u_{g'h'}) \text{ for all } g,g' \in V(G), h,h' \in V(H).
\end{equation}
Moreover, if the entries of $U$ commute, then the correlation $p$ is classical.
\end{lemma}

We will see that the converse of the final claim of the above lemma does not hold, i.e., there are classical correlations that arise from tracial states and quantum permutation matrices whose entries do not commute. However, it is well known that a partial converse does hold, since for any classical correlation $p$ that wins the isomorphism game there is \emph{some} quantum permutation matrix with commuting entries that produces $p$. 

Any correlation which is non-classical, i.e., cannot be obtained as a convex combination of deterministic classical correlations, is said to be \emph{nonlocal}, as it cannot be explained by any local realistic theory. We will use nonlocal and non-classical interchangeably in regards to correlations. In the case of the $(G,H)$-isomorphism game for non-isomorphic $G$ and $H$, any winning quantum correlation is nonlocal since there are no winning classical correlations. In this case, the graphs are said to be \emph{quantum isomorphic}\footnote{Isomorphic graphs are also said to be quantum isomorphic, the latter being a strict relaxation of the former.}, and we write $G \cong_{qc} H$. On the other hand if $G$ and $H$ are isomorphic, then there do exist winning classical correlations. In this case, though they cannot rise to the level of pseudotelepathy, there may still exist some nonlocal quantum correlations that win the $(G,H)$-isomorphism game. As they are isomorphic, we may as well assume that $G = H$, in which case we refer to the game as the $G$-automorphism game. We now define the central notion of this work:
\begin{definition}
For a graph $G$, we define a \emph{nonlocal symmetry} of $G$ to be a winning quantum correlation for the $G$-automorphism game that is not classical. If there exists such a correlation for a graph $G$, then we say that $G$ \emph{admits nonlocal symmetry}.
\end{definition}
We can apply Lemma~\ref{lem:qisocorr} to the $G = H$ case to obtain a characterization of the winning quantum correlations for the $G$-automorphism game, but we can also characterize them in terms of the quantum automorphism group of $G$. This is done for the complete graph in Theorem 2.2 of~\cite{bisynchronous}, and the following is a straightforward generalization.

\begin{lemma}\label{lem:qcorrgroup}
Let $G$ be a graph and let $U = (u_{ij})$ be the fundamental representation of $\qut(G)$. Then $p$ is a winning quantum correlation for the $G$-automorphism game if and only if there exists a tracial state $\tau$ on $C(\qut(G))$ such that
\[p(l,k|i,j) = \tau(u_{il}u_{jk}).\]
\end{lemma}

Note that a direct application of Lemma~\ref{lem:qisocorr} would be phrased in terms of the existence of \emph{some} quantum permutation $U'$ and $C^*$-algebra $\mathcal{A}$ with a tracial state $\tau$ satisfying certain properties. In contrast, by invoking the quantum automorphism group of $G$, Lemma~\ref{lem:qcorrgroup} is in terms of the existence of a tracial state on the \emph{universal} $C^*$-algebra $C(\qut(G))$ generated by the entries of the fundamental representation $U$. Though formally different these are equivalent since by the universal property of $C(\qut(G))$ there is necessarily a $*$-homomorphism $\phi : C(\qut(G)) \to \mathcal{A}$ which maps $u_{ij}$ to $u'_{ij}$ and then $\tau \circ \phi$ provides the necessary tracial state on $C(\qut(G))$.

By Lemma~\ref{lem:qcorrgroup} and the final claim of Lemma~\ref{lem:qisocorr}, we see that if $C(\qut(G))$ is commutative (i.e.~$G$ has no quantum symmetry), then $G$ does not admit nonlocal symmetry. Thus if $G$ admits nonlocal symmetry then it necessarily has quantum symmetry. Interestingly, the converse does not hold: the complete graph $K_4$ is known to have quantum symmetry but our Theorem~\ref{sec:K4} shows that it does not admit nonlocal symmetry. Therefore nonlocal symmetry is a strictly stronger property than quantum symmetry.


We will denote the set of all classical, or local, correlations that win the $G$-automorphism game as $L(G)$, and the corresponding quantum set as $Q(G)$. Recall that these are both closed convex sets and that $L(G) \subseteq Q(G)$ for any graph $G$. Moreover, a graph $G$ admits nonlocal symmetry if and only if $Q(G) \ne L(G)$. In the case that $G$ is the complete graph $K_n$, we will simply denote these sets as $L(n)$ and $Q(n)$. We remark that the sets $L(G)$ and $Q(G)$ are invariant under taking the graph complement, i.e., $L(\overline{G}) = L(G)$ and $Q(\overline{G}) = Q(G)$. This is simply because the $(G,H)$-isomorphism game is identical to the $(\overline{G},\overline{H})$-isomorphism game. Also note that if $G$ is a graph on $n$ vertices, then $L(G) \subseteq L(n)$ and $Q(G) \subseteq Q(n)$, since the constraints of the $G$-automorphism game contain the constraints of the $K_n$-automorphism game.

Lemma~\ref{lem:qcorrgroup} motivates us to extend our notion of nonlocal symmetry to quantum permutation groups, i.e, quantum subgroups of the quantum symmetric groups $S_n^+$:
\begin{definition}
For a given quantum permutation group $\mathbb{G} \subseteq S_n^+$ with fundamental representation $U = (u_{ij})$, we say that $\mathbb{G}$ \emph{exhibits nonlocal symmetry} if there exists a tracial state $\tau$ on $C(\mathbb{G})$ such that the correlation $p(l,k|i,j) = \tau(u_{il}u_{jk}) \in Q(n)$ is nonlocal.
\end{definition}

We remark that quantum symmetry, though originating from purely mathematical considerations, does have a quantum mechanical interpretation. A collection of projective measurements, such as the rows/columns of a quantum permutation matrix, is said to be \emph{incompatible} if not all of the measurement operators commute. Thus if $G$ has quantum symmetry, then the rows/columns of the fundamental representation of $\qut(G)$ form a set of incompatible measurements. However, this is still a property of our model of reality rather than our observations of reality like nonlocal symmetry.

\section{Nonlocal symmetry in complete graphs} \label{sec:complete}

Here we will show that $K_4$ does not admit any nonlocal symmetry but $K_5$ does. In terms of quantum groups this says that the quantum symmetric group $S_4^+$ does not exhibit nonlocal symmetry but $S_5^+$ does. This is in contrast to quantum symmetry, which already occurs for the complete graphs on four vertices (but not three vertices). In fact even the graph on four vertices containing a single edge has quantum symmetry. The fundamental representation of the quantum automorphism group of this graph has the form
\[\begin{pmatrix} p & 1 - p & 0 & 0 \\ 1-p & p & 0 & 0 \\  0 & 0 & q & 1 - q \\ 0 & 0 & 1-q & q \end{pmatrix},\]
where $p$ and $q$ are free projections, and thus do not commute. However, our proof that $K_4$ does not admit nonlocal symmetry implies that neither does this graph nor any graph on four vertices. We remark that a direct proof for this graph is much simpler, and follows as a special case of our results in Section~\ref{subsec:freeproduct}. In fact the lack of nonlocal symmetry for $K_4$ is the most challenging proof in this work. We achieve the result by giving a rather nice description of $L(4)$ in terms of linear equalities and inequalities (Theorem~\ref{thm:ineqchar}), and then showing that any quantum correlation from $Q(4)$ must also satisfy these equalities and inequalities (Theorem~\ref{thm:Q4=C4}). To show that $K_5$ does admit nonlocal symmetry, we consider winning correlations for the $K_5$-automorphism game that possess additional invariance properties. We show that any element of $L(5)$ with these invariance properties must satisfy a certain inequality. We then present an explicit quantum correlation that has these invariance properties but does not satisfy the inequality.

\subsection{$K_4$ has no nonlocal symmetry}\label{sec:K4}

We begin by defining a set $B(n)$ of correlations that contain $Q(n)$, inspired by the so-called non-signalling correlations. We will show that $L(4)$ can be obtained from $B(4)$ by imposing an additional family of inequalities. We then show that any quantum correlation from $Q(4)$ must satisfy all of these additional inequalities and thus $Q(4) = L(4)$. 

\begin{definition}\label{def:bicorr}
For $N \in \mathbb{N}$, define the \emph{bijective correlations} on the set $[n]$, denoted $B(n)$, to be the set of correlations $p : [n]^4 \to [0,1]$ satisfying the following:
\begin{enumerate}
\item $p(l,k|i,j) = 0$ if $\delta_{ij} \ne \delta_{lk}$.
\item $p(l,k|i,j) = p(k,l|j,i)$ for all $i,j,l,k \in [n]$.
\item For each $a,b \in [n]$, there exists a \emph{marginal probability}, denoted $p(a|b)$ such that
\begin{align*}
\sum_{k \in V(G)} p(a,k|b,j) &= p(a|b) \text{ for all } j \in V(G), \\
\sum_{j \in V(G)} p(a,k|b,j) &= p(a|b) \text{ for all } k \in V(G), \\
\sum_{l \in V(G)} p(l,a|i,b) &= p(a|b) \text{ for all } i \in V(G), \\
\sum_{i \in V(G)} p(l,a|i,b) &= p(a|b) \text{ for all } l \in V(G).
\end{align*}
The last two conditions above are redundant given Condition (2), but we list them explicitly anyways.
\end{enumerate}
\end{definition}
Note that though we did not explicitly list it, the use of the term ``correlation" in the definition of $B(n)$ means that $\sum_{l,k} p(l,k|i,j) = 1$ for all $i,j \in [n]$ for any $p \in B(n)$. Also note that $B(n)$ is easily seen to be convex. The set $B(n)$ also has the nice property that it is preserved under swapping inputs and outputs:

\begin{prop}
Let $p \in B(n)$. Then $\sum_{i,j \in [n]}p(l,k|i,j) = 1$ for all $l,k \in [n]$ and thus $p' \in B(n)$ where $p'(l,k|i,j) = p(i,j|l,k)$.
\end{prop}
\proof
The definition of $B(n)$ is symmetric with respect to inputs and outputs with the exception of the condition $\sum_{l,k}p(l,k|i,j) = 1$. Thus proving the condition in the proposition statement does indeed imply the conclusion. For any $l,k \in [n]$,
\[\sum_{i,j} p(l,k|i,j) = \sum_i p(l|i) = \sum_{i,k'} p(l,k'|i,j') = \sum_{k'} p(k'|j') = \sum_{l',k'} p(l',k'|i',j') = 1.\]\qeds

The above also shows that the matrix whose $a,b$-entry is $p(a|b)$ is a doubly stochastic matrix. Showing that $Q(n) \subseteq B(n)$ is straightforward:

\begin{lemma}
For any $n \in \mathbb{N}$, we have that $Q(n) \subseteq B(n)$.
\end{lemma}
\proof
Suppose that $p \in Q(n)$. Then $p$ satisfies Condition (1) by definition. By Lemma~\ref{lem:qisocorr} there exists a quantum permutation matrix $U$ and tracial state $\tau$ such that $p(l,k|i,j) = \tau(u_{il}u_{jk})$. Condition (2) now follows immediately from the cyclicity of tracial states. For Condition (3), define $p(a|b) = \tau(u_{ba})$ and note that
\[\sum_{k \in V(G)} p(a,k|b,j) = \sum_{k \in V(G)} \tau(u_{ba}u_{jk}) =  \tau\left(u_{ba}\sum_{k \in V(G)} u_{jk}\right) = \tau(u_{ba}) = p(a|b),\]
and similarly for the other equations in Condition (3).\qeds

We remark that the conditions that the sum $\sum_{k \in V(G)} p(l,k|i,j)$ is independent of $j$ and $\sum_{l \in V(G)} p(l,k|i,j)$ is independent of $i$ are known as the \emph{non-signalling} conditions for correlations. This is the mathematical formalization of the requirement that Alice and Bob are unable to communicate any information to each other, as it says that Alice's probability distribution over her outputs conditioned on her inputs is independent of Bob's input, and vice versa. Note however that $B(n)$ is \emph{not} the set of non-signalling correlations that win the $K_n$-automorphism game, since we impose additional constraints on the elements of $B(n)$. For an explicit example of a non-signalling correlation that wins the $K_3$-automorphism game but is not an element of $B(3)$, see Remark 2.8 of~\cite{bisynchronous}. There they define \emph{bisynchronous} correlations as correlations $p : [m]^2 \times [n]^2 \to [0,1]$ such that $p(l,k|i,j) = 0$ if $\delta_{ij} \ne \delta_{lk}$. Our bijective correlations $B(n)$ are a subset of the bisynchronous correlations in the case $m = n$. 

We will also be concerned with the span of the convex sets $L(n)$, $Q(n)$, and $B(n)$. Note that if $p \in \spn B(n)$, then $p \in B(n)$ if and only if $p(l,k|i,j) \ge 0$ for all $i,j,l,k \in [n]$, and $\sum_{l,k \in V(G)} p(l,k|i,j) = 1$ for all $i,j \in [n]$, but the same is not necessarily (and in fact is not) true for $L(n)$ and $Q(n)$. This means that $\spn B(n)$ consists of precisely the functions $p: [n]^4 \to \mathbb{R}$ that satisfy Conditions (1)--(3) of Definition~\ref{def:bicorr}. In particular, this means that we can define $p(a|b)$ for $a,b \in [n]$ for any $p \in \spn B(n)$ in the obvious way.

Given $p \in \spn B(n)$, we see that $p(l|i) = \sum_{k} p(l,k|i,i) = p(l,l|i,i)$. Thus $p(l,k|i,i) = \delta_{lk}p(l|i)$ and so these values of $p$ are determined by the values of the form $p(l,k|i,j)$ for $i \ne j$ since we can use the latter to determine the marginals $p(l|i)$. Furthermore, by Condition (2) we need only consider the case $i < j$, and by Condition (1) we can further restrict to $l \ne k$. Let $\mathcal{T}(n)$ denote the tuples $(l,k|i,j) \in [n]^4$ such that $i < j$ and $l \ne k$. Then we have the following:

\begin{lemma}\label{lem:Tenough}
Let $p,p' \in \spn B(n)$. Then $p = p'$ if and only if $p(l,k|i,j) = p'(l,k|i,j)$ for all $(l,k|i,j) \in \mathcal{T}(n)$.
\end{lemma}

We will identify the elements of $\mathcal{T}(4)$ with certain pairs of permutations from $S_4$ which we will then take as edges of a graph. Before we do this we will discuss the classical bijective correlations which provide the inspiration.

Given an $n \in \mathbb{N}$, and $\pi \in S_n$, we define the correlation $p_\pi$ as
\begin{equation}\label{eq:detcorrs}
p_\pi(l,k|i,j) = \begin{cases}1 & \text{if } l = \pi(i) \ \& \ k = \pi(j) \\ 0 & \text{o.w.}\end{cases}
\end{equation}
The correlations $p_\pi$ for $\pi \in S_n$ are precisely the deterministic classical correlations that win the $K_n$-automorphism game, and thus we have that
\[L(n) = \conv\{p_\pi : \pi \in S_n\}.\]
Therefore, if $p \in L(n)$ then $p$ can be written as a convex combination $p = \sum_\pi \alpha_\pi p_\pi$. Note that this convex combination need not be unique. For instance, if $n \ge 4$ then the alternating group $A_n$ is doubly transitive and thus taking a uniform mixture of $p_\pi$ for $\pi \in A_n$ will produce the same classical correlation as taking a uniform mixture over the full symmetric group. Despite this non-uniqueness, the following holds for any valid choice of the $\alpha_\pi$:
\begin{equation}\label{eq:classical}
    p(l,k|i,j)= \sum_{\substack{\pi \in S_n \\ \pi(i) = l, \pi(j) = k}} \alpha_\pi.
\end{equation}
Note that this holds also for $p \in\spn L(n)$. In the case of $n = 4$, and assuming that $i \ne j$ and $l \ne k$, the constraints $\pi(i) = l$ and $\pi(j) = k$ are satisfied by precisely two permutations $\pi$. Indeed, if $a,b,c,d \in [4]$ are such that $\{a,b,i,j\} = \{c,d,l,k\} = [4]$, then $\pi(i) = l$ and $\pi(j) = k$ hold if and only if $(\pi(i),\pi(j),\pi(a),\pi(b)) = (l,k,c,d)$ or $(l,k,d,c)$. This allows us to construct a bijection between the elements of $\mathcal{T} := \mathcal{T}(4)$ and certain pairs of permutations from $S_4$, which we consider to be edges of a graph. Given $(l,k|i,j) \in \mathcal{T}$, define $F(l,k|i,j)$ to be equal to the set $\{\pi,\sigma\}$ where $\pi,\sigma \in S_4$ are the two permutations that map $i$ to $l$ and $j$ to $k$. Let $G_4$ be the \emph{transposition graph} of $S_4$. In other words the vertex set of $G_4$ is $S_4$ and $\pi,\sigma \in S_4$ are adjacent if $\pi^{-1}\sigma$ is a transposition.

\begin{lemma}\label{lem:Fbij}
The map $F$ is a bijection between $\mathcal{T}$ and $E(G_4)$.
\end{lemma}
\proof
Let $(l,k|i,j) \in \mathcal{T}$ and let $\pi, \sigma \in S_4$ be such that $F(l,k|i,j) = \{\pi,\sigma\}$. Let $a,b \in [4]$ be such that $a,b,i,j$ are distinct. Then the product $\pi(ab)$ must map $i$ to $l$ and $j$ to $k$ and thus $\sigma = \pi(ab)$, i.e., $\pi^{-1}\sigma = (ab)$. This shows that $F$ maps $\mathcal{T}$ into $E(G_4)$. To see that $F$ is surjective, suppose that $\{\pi,\sigma\} \in E(G_4)$ and let $a,b \in [4]$ be such that $\pi^{-1}\sigma = (ab)$. Let $i<j \in [4]$ be such that that $a,b,i,j$ are distinct. Then $\pi^{-1}\sigma = (ab)$ fixes $i$ and $j$. It follows that $\sigma(i) = \pi(i)$ and $\sigma(j) = \pi(j)$. Let $l,k \in [4]$ be the values that $\pi$ and $\sigma$ map $i$ and $j$ to respectively. We then have that $F(l,k|i,j) = \{\pi,\sigma\}$. Thus we have shown that $F$ is surjective. Injectivity follows from the fact that both sets have cardinality 72.\qeds

The above allows us to identify the correlations of $B(4)$ with functions on the edges of $G_4$. Given $p \in \spn B(4)$, define $\hat{p}: E(G_4) \to \mathbb{R}$ as $\hat{p}(\{\pi,\sigma\}) = p(l,k|i,j)$ such that $F(l,k|i,j) = \{\pi,\sigma\}$. By Lemmas~\ref{lem:Tenough} and~\ref{lem:Fbij} we immediately have the following:

\begin{lemma}\label{lem:hatenough}
Let $p,p' \in \spn B(4)$. Then $p = p'$ if and only if $\hat{p} = \hat{p}'$.
\end{lemma}

If $p \in \spn L(4)$, then $p = \sum_{\pi} \alpha_\pi p_\pi$ for some arbitrary coefficients $\alpha_\pi \in \mathbb{R}$. By Equation~\eqref{eq:classical} and the discussion that followed, we have that $p(l,k|i,j) = \alpha_\pi + \alpha_\sigma$ where $F(l,k|i,j) = \{\pi,\sigma\}$. In other words, $\hat{p}(\{\pi,\sigma\}) = \alpha_\pi + \alpha_\sigma$.

We now show that the set of classical correlations $L(4)$ is essentially the convex hull of the columns of the \emph{incidence matrix} of $G_4$. The incidence matrix of a graph $G$ is a matrix with rows and columns indexed by the edges and vertices of $G$ respectively, such that the $e,v$-entry is 1 if $v$ is an endpoint of $e$ and is 0 otherwise. 

\begin{lemma}\label{lem:linbij}
Let $M$ be the incidence matrix of $G_4$. The map $p \mapsto \hat{p}$ on $\spn B(4)$ is a linear injection. Moreover, its restriction to $L(4)$ (resp.~$\spn L(4)$) is a bijection to the convex hull (resp.~span) of the columns of $M$.
\end{lemma}
\proof
The fact that the map $p \mapsto \hat{p}$ is linear is immediate since the map is simply a restriction to a subset of the values/coordinates of $p$. The fact that it is injective follows from Lemma~\ref{lem:hatenough}.

To prove the bijection claims, it suffices to show that the columns of $M$ are precisely the elements $\hat{p}_\pi$ for $\pi \in S_4$. The claims then follow by convexity/linearity respectively.

Let $\pi \in S_4$ and consider $p_\pi \in L(4)$. We will show that $\hat{p}_\pi$ is the $\pi$ column of $M$. If $\hat{p}_\pi(e) = 1$ then $e = \{\sigma_1,\sigma_2\} = F(l,k|i,j)$ for some $(l,k|i,j) \in \mathcal{T}$ such that $p_\pi(l,k|i,j) = 1$. This implies that $\pi(i) = l$ and $\pi(j) = k$, and thus $\pi \in \{\sigma_1,\sigma_2\}$ by the definition of $F$. In other words, if $\hat{p}_\pi(e) =1$, then $e$ is an edge incident to $\pi$ in $G_4$. Conversely, if $e$ is an edge incident to $\pi$ in $G_4$, then $e = \{\pi,\sigma\}$ and there exists $(l,k|i,j) \in \mathcal{T}$ such that $F(l,k|i,j) = \{\pi,\sigma\}$ by Lemma~\ref{lem:Fbij}. Moreover, we must have that $\pi(i) = l = \sigma(i)$ and $\pi(j) = k = \sigma(j)$. In particular this means that $p_\pi(l,k|i,j) = 1$ and thus $\hat{p}_\pi(e) = 1$. Thus $\hat{p}_\pi(e) = 1$ if and only if $e$ is an edge incident to $\pi$ in $G_4$, and $\hat{p}_\pi$ is zero elsewhere since $p$ is 01-valued. In other words, $\hat{p}_\pi$ is precisely the column of $M$ corresponding to $\pi \in V(G_4)$.

By the above, we have that $\{\hat{p} : p \in L(4)\}$ is the convex hull of the columns of $M$ and $\{\hat{p} : p \in \spn L(4)\}$ is the span of the columns of $M$. The bijectivity of the map $p \mapsto \hat{p}$ on these sets follows from its injectivity on $\spn B(4)$ which contains these sets.\qeds

The above allows us to determine the dimension of $\spn L(4)$ since it must be equal to the rank of the incidence matrix of $G_4$. The rank of the incidence matrix of a graph on $n$ vertices is $n-b$ where $b$ is the number of connected components which are bipartite. The graph $G_4$ is connected since the transpositions generate the whole group $S_4$. Moreover, it is bipartite since if $\pi$ and $\sigma$ are adjacent then $\pi = \sigma(ab)$ for some $a,b \in [4]$, and thus exactly one of them is an odd permutation and the other even. Therefore $\spn L(4)$ has dimension $24 - 1 = 23$. Knowing this allows us to prove the following:

\begin{lemma}\label{lem:spacesequal}
$\spn L(4) = \spn B(4)$.
\end{lemma}
\proof
The containment $\spn L(4) \subseteq \spn B(4)$ is trivial. To show the equality, it suffices to show that the dimension of $\spn B(4)$ is at most the dimension of $\spn L(4)$, which is 23 by the discussion above.

To show that $\spn B(4)$ has dimension at most 23, it suffices to show that there is a subset $S \subseteq [4]^4$ of size 23 such that for $p,p' \in \spn B(4)$ we have that $p = p'$ whenever $p(s) = p'(s)$ for all $s \in S$. We take
\[S = \{(l,k|1,4) : l \ne k\} \cup \{(l,k|1,3) : k \equiv l+1 \text{ or } l+2 \ \mathrm{mod} \ 4\} \cup \{(1,3|2,3), (4,3|2,3), (2,4|2,4)\}.\]
We remark that $S \subseteq \mathcal{T}$. In order to prove the claim, we must show that if we know the values of $p(s)$ for $s \in S$, then we can determine all values $p(l,k|i,j)$ for $i,j,l,k \in [4]$. This is done by using the linear equalities in Definition~\ref{def:bicorr} that define $\spn B(4)$. There appears to be no non-tedious way of doing this. Thus instead of presenting every detail, we provide Table~\ref{tab:spantree}. This is a $16 \times 16$ table split into $4 \times 4$ blocks. The $l,k$-cell in the $i,j$ block represents $(l,k|i,j)$. We have placed X's in the cells corresponding to the elements of $S$, and the elements of $S^T = \{(k,l|j,i) : (l,k|i,j) \in S)$ since these must have the same value by Condition (2) of Definition~\ref{def:bicorr}. We have also placed 0's in the cells corresponding to $(l,k|i,j)$ where $p \in \spn B(4)$ must be zero by Condition (1). Starting from this table, one can repeatedly use Condition (3) to fill in more and more empty cells that correspond to values of $p \in \spn B(4)$ that are determined by the values of $p$ on the already filled in cells. For example, since the cells corresponding to $(1,1|1,4), (1,2|1,4), (1,3|1,4), (1,4|1,4)$ and $(1,1|1,3), (1,2|1,3), (1,3|1,3)$ are filled in, we can fill in the cell corresponding to $(1,4|1,3)$ since 
\[\sum_k p(1,k|1,3) = p(1|1) = \sum_k p(1,k|1,4)\]
and the only term corresponding to a unfilled cell in the equation above is $p(1,4|1,3)$. Similarly, since the above shows that $p(1|1)$ must be determined and the cells corresponding to $(1,2|1,1), (1,2|1,3), (1,2|1,4)$ are filled, we may fill in the cell corresponding to $(1,2|1,2)$. Continuing on, it is straightforward to fill in all of the cells (one can also use Condition (2)), but we leave this as an exercise for the reader. We do provide an explanation for how the set $S$ was found in Remark~\ref{rem:findS} below.\qeds

\begin{table}[ht]
    \centering
    \resizebox{\columnwidth}{!}{
    \begin{tabu}{|[2pt]c|c|c|c|[2pt]c|c|c|c|[2pt]c|c|c|c|[2pt]c|c|c|c|[2pt]}
        \tabucline[2pt]{-}
          & 0 & 0 & 0  &  0 &   &   &    &  0 & X & X &    &  0 & X & X & X \\
        \hline
        0 &   & 0 & 0  &    & 0 &   &    &    & 0 & X & X  &  X & 0 & X & X \\
        \hline
        0 & 0 &   & 0  &    &   & 0 &    &  X &   & 0 & X  &  X & X & 0 & X \\
        \hline
        0 & 0 & 0 &    &    &   &   & 0  &  X & X &   & 0  &  X & X & X & 0 \\
        \tabucline[2pt]{-}
        0 &   &   &    &    & 0 & 0 & 0  &  0 &   & X &    &  0 &   &   &   \\
        \hline
          & 0 &   &    &  0 &   & 0 & 0  &    & 0 &   &    &    & 0 &   & X \\
        \hline
          &   & 0 &    &  0 & 0 &   & 0  &    &   & 0 &    &    &   & 0 &   \\
        \hline
          &   &   & 0  &  0 & 0 & 0 &    &    &   & X & 0  &    &   &   & 0 \\
        \tabucline[2pt]{-}
        0 &   & X & X  &  0 &   &   &    &    & 0 & 0 & 0  &  0 &   &   &   \\
        \hline
        X & 0 &   & X  &    & 0 &   &    &  0 &   & 0 & 0  &    & 0 &   &   \\
        \hline
        X & X & 0 &    &  X &   & 0 & X  &  0 & 0 &   & 0  &    &   & 0 &   \\
        \hline
          & X & X & 0  &    &   &   & 0  &  0 & 0 & 0 &    &    &   &   & 0 \\
        \tabucline[2pt]{-}
        0 & X & X & X  &  0 &   &   &    &  0 &   &   &    &    & 0 & 0 & 0 \\
        \hline
        X & 0 & X & X  &    & 0 &   &    &    & 0 &   &    &  0 &   & 0 & 0 \\
        \hline
        X & X & 0 & X  &    &   & 0 &    &    &   & 0 &    &  0 & 0 &   & 0 \\
        \hline
        X & X & X & 0  &    & X &   & 0  &    &   &   & 0  &  0 & 0 & 0 &   \\
        \tabucline[2pt]{-}
    \end{tabu}
    }
    \caption{The elements of $S$ and $S^T$ from the proof of Lemma~\ref{lem:spacesequal}.}
    \label{tab:spantree}
\end{table}

\begin{remark}\label{rem:findS}
The set $S$ used in the proof of Lemma~\ref{lem:spacesequal} above was not simply guessed, nor found by computer search. The idea is that, if the claim of the lemma holds, then any set of coordinates that determines any element of $\spn L(4)$ must also work for $\spn B(4)$. For the former, by Lemma~\ref{lem:linbij} such a set corresponds to a subset of 23 rows of the incidence matrix $M$ of $G_4$ such that the submatrix of $M$ consisting of those rows still has rank 23. Recall that the rows of $M$ correspond to edges of $G_4$ and thus the submatrix of $M$ consisting of the rows indexed by a subset $T \subseteq E(G_4)$ is the incidence matrix of the subgraph $H$ of $G_4$ with vertex set $V(G_4)$ and edge set $T$. Since $G_4$ is bipartite, any such subgraph $H$ is also bipartite, and thus for its incidence matrix to have rank 23 it must be connected. As we want $|T| = 23$ we must choose $H$ with 23 edges, and since $H$ has 24 vertices this means that $H$ must be a spanning tree of $G_4$. Conversely, any spanning tree $H$ of $G_4$ will suffice, and the subset $S \subseteq \mathcal{T}$ such that $F(S) = E(H)$ will provide us with a set of 23 coordinates that determine any element of $\spn L(4)$. The same set will determine any element of $\spn B(4)$ if and only if the claim holds, and then it is just a matter of checking the former in the manner described in the proof.
\end{remark}



The above two lemmas will be an important part of proving that $L(4)$ consists of the correlations in $B(4)$ that satisfy a certain family of inequalities. These inequalities were found using a computer, but we give a proof that they characterize $L(4)$.

\begin{theorem}\label{thm:ineqchar}
Let $p \in B(4)$. Then $p \in L(4)$ if and only if for every $\pi \in S_4$ and 4-cycle $(abcd) \in S_4$,
\[p(\pi(c),\pi(d)|c,d) - p(\pi(b),\pi(d)|a,d) + p(\pi(b),\pi(c)|a,b) \ge 0.\]
\end{theorem}
\proof
Suppose that $p \in B(4)$. By Lemma~\ref{lem:spacesequal}, we can write $p$ as a \emph{linear} combination $p = \sum_{\sigma \in S_4} \alpha_\sigma p_\sigma$. Let
\[\pi_1 = \pi, \quad \pi_2 = \pi(ab), \quad \pi_3 = \pi_2(bc) = \pi(abc), \quad \pi_4 = \pi_3(cd)= \pi(abcd).\]
It is straightforward to check that
\begin{align*}
p(\pi(c),\pi(d)|c,d) &= \hat{p}(\{\pi_1,\pi_2\}) = \alpha_{\pi_1} + \alpha_{\pi_2} \\
p(\pi(b),\pi(d)|a,d) &= \hat{p}(\{\pi_2,\pi_3\}) = \alpha_{\pi_2} + \alpha_{\pi_3} \\
p(\pi(b),\pi(c)|a,b) &= \hat{p}(\{\pi_3,\pi_4\}) = \alpha_{\pi_3} + \alpha_{\pi_4}.
\end{align*}
It follows that
\begin{align*}
p(\pi(c),\pi(d)|c,d) - p(\pi(b),\pi(d)|a,d) + &p(\pi(b),\pi(c)|a,b)\\
&= (\alpha_{\pi_1} + \alpha_{\pi_2}) - (\alpha_{\pi_2} + \alpha_{\pi_3}) + (\alpha_{\pi_3} + \alpha_{\pi_4}) \\
&= \alpha_{\pi_1} + \alpha_{\pi_4} \\
&= \alpha_\pi + \alpha_{\pi(abcd)}.
\end{align*}
Now if $p \in L(4)$, then we could have chosen the $\alpha_\sigma$ to be nonnegative for all $\sigma \in S_4$. Then we would have $\alpha_{\pi_1} + \alpha_{\pi_4} \ge 0$ and thus the inequalities in the theorem statement will hold.

Conversely, suppose that $p$ satisfies all of the inequalities described in the theorem statement, i.e., that $\alpha_\pi + \alpha_{\pi(abcd)} \ge 0$ for all $\pi \in S_4$ and all 4-cycles $(abcd) \in S_4$. We will show that the values of the $\alpha_\sigma$ can be adjusted so that they are all nonnegative and thus $p \in L(4)$. By Lemmas~\ref{lem:linbij} and~\ref{lem:spacesequal}, we have that $\hat{p}$ is in the column space of the incidence matrix $M$ of $G_4$. In fact, since $p = \sum_\sigma \alpha_\sigma p_\sigma$, we have that
\[\hat{p} = \sum_{\sigma} \alpha_\sigma \hat{p}_\sigma = M\alpha,\]
where $\alpha = (\alpha_\sigma)_{\sigma \in S_4}$. We now aim to show that if $\pi$ and $\sigma$ are even and odd permutations respectively, then $\alpha_\pi + \alpha_\sigma \ge 0$. Since $\pi$ and $\sigma$ are even and odd respectively, we have that $\pi^{-1}\sigma$ is odd. There are only two types of odd permutations in $S_4$: the 4-cycles and the transpositions. If $\pi^{-1}\sigma$ is equal to a 4-cycle $(abcd)$, then we have that $\sigma = \pi(abcd)$ and so $\alpha_\pi + \alpha_\sigma = \alpha_\pi + \alpha_{\pi(abcd)} \ge 0$ by our assumption. If $\pi^{-1}\sigma$ is a transposition, then $\pi$ and $\sigma$ are adjacent via an edge $e$ in $G_4$. It follows that
\[\hat{p}(e) = \left(M\alpha\right)_e = \alpha_\pi + \alpha_\sigma.\]
Since $\hat{p}(e) = p(l,k|i,j) \ge 0$ for some $i,j,l,k \in [4]$, we have that $\alpha_\pi + \alpha_\sigma \ge 0$ as desired. Thus $\alpha_\pi + \alpha_\sigma \ge 0$ whenever $\pi$ and $\sigma$ are even and odd respectively.

Now let $\gamma = (\gamma_\pi)_{\pi \in S_4}$ be the vector defined as $\gamma_\pi = (-1)^{\mathrm{sgn}(\pi)}$, i.e., $\gamma$ is $-1$ on odd permutations and $+1$ on even permutations. It is easy to see that $M \gamma = 0$. Thus $M(\alpha + c\gamma) = \hat{p}$ for any constant $c$. Pick $c$ such that $\min\{(\alpha + c\gamma)_\pi : \pi \in S_4 \text{ is even}\} = 0$, and let $\pi^*$ be an even permutation such that $(\alpha +c\gamma)_{\pi^*} = 0$. Define $\beta = \alpha +c \gamma$. Then by choice of $c$, we have that $\beta_\pi \ge 0$ for all even permutations $\pi$. On the other hand, if $\sigma$ is an odd permutation, then by the above we have that $\beta_{\sigma} =  \beta_{\pi^*} + \beta_{\sigma} \ge 0$. Thus we have that $\beta$ is a nonnegative vector satisfying $M\beta = \hat{p}$ and therefore $p = \sum_{\pi} \beta_\pi p_\pi$.

Lastly, we must show that $\sum_{\pi \in S_4} \beta_\pi = 1$. We have that
\[\sum_{e \in E(G_4)} \left(M\beta\right)_e = \sum_{e \in E(G_4)} \hat{p}(e).\]
The lefthand side is equal to
\[\sum_{e \in E(G_4)} \sum_{\pi \in S_4} M_{e, \pi} \beta_{\pi} = \sum_{\pi \in S_4} \beta_\pi \left(\sum_{e \in E(G_4)} M_{e,\pi}\right) = 6 \sum_{\pi \in S_4} \beta_\pi,\]
since $G_4$ is 6-regular. By Lemma~\ref{lem:Fbij}, the righthand side is equal to
\[\sum_{i < j \in [4]} \sum_{\ell,k \in [4]} p(l,k|i,j) = \sum_{i < j \in [4]} 1 = 6.\]
Thus we have that $\sum_{\pi \in S_4} \beta_\pi = 1$ and therefore $p$ is a convex combination of the deterministic classical correlations $p_\pi$. Therefore $p \in L(4)$ as desired.\qeds

Now that we have shown that $L(4)$ consists of the correlations from $B(4)$ that satisfy the inequalities described in Theorem~\ref{thm:ineqchar}, in order to show that $Q(4) = L(4)$ (i.e., that $K_4$ does not admit nonlocal symmetry) we only need to show that every element of $Q(4)$ satisfies these inequalities. This we proceed to do.

In the proof of the following lemma, as well as throughout the rest of this subsection, we will make frequent use of the relations $\sum_{j} u_{ij} = \vec{1} = \sum_{i} u_{ij}$ and $u_{ij}u_{lk} = 0$ if $\delta_{il} \ne \delta_{jk}$ in the proof throughout the rest of this subsection.

\begin{lemma}\label{lem:triples}
Let $\mathcal{A}$ be a $C^*$-algebra, and $U = (u_{ij}) \in M_4(\mathcal{A})$ be a quantum permutation matrix. If $\pi \in S_4$ and $\{a,b,c,d\} = [4]$, then
\[u_{a\pi(a)}u_{b\pi(b)}u_{c\pi(c)} = u_{a\pi(a)}u_{d\pi(d)}u_{c\pi(c)}.\]
Moreover, this implies that if $\tau$ is a tracial state on $\mathcal{A}$ then
\[\tau\left(u_{a\pi(a)}u_{b\pi(b)}u_{c\pi(c)}\right) = \tau\left(u_{x\pi(x)}u_{y\pi(y)}u_{z\pi(z)}\right)\]
for any distinct $x,y,z \in [4]$.
\end{lemma}
\proof
We have that
\begin{align*}
u_{a\pi(a)}u_{b\pi(b)}u_{c\pi(c)} &= u_{a\pi(a)}\left(\vec{1} - u_{a\pi(b)} - u_{c\pi(b)} - u_{d\pi(b)}\right)u_{c\pi(c)} \\
&= u_{a\pi(a)}\left(\vec{1} - u_{d\pi(b)}\right)u_{c\pi(c)} \\
&= u_{a\pi(a)}\left(u_{d\pi(a)} + u_{d\pi(c)} + u_{d\pi(d)}\right)u_{c\pi(c)} \\
&= u_{a\pi(a)}u_{d\pi(d)}u_{c\pi(c)}.
\end{align*}
Now, by the above and the cyclicity of the trace, we have that
\begin{align*}
\tau\left(u_{a\pi(a)}u_{b\pi(b)}u_{c\pi(c)}\right) &= \tau\left(u_{a\pi(a)}u_{d\pi(d)}u_{c\pi(c)}\right) \\
&= \tau\left(u_{c\pi(c)}u_{a\pi(a)}u_{d\pi(d)}\right) \\
&= \tau\left(u_{c\pi(c)}u_{b\pi(b)}u_{d\pi(d)}\right) \\
&= \tau\left(u_{b\pi(b)}u_{d\pi(d)}u_{c\pi(c)}\right) \\
&= \tau\left(u_{b\pi(b)}u_{a\pi(a)}u_{c\pi(c)}\right).
\end{align*}
This, plus cyclicity of the trace, shows that the order of $u_{a\pi(a)},u_{b\pi(b)},u_{c\pi(c)}$ in the trace does not matter. Moreover, by the previous part, the elements $a,b,c$ can be replaced by any distinct elements $x,y,z \in [4]$.\qeds

Recall from Theorem~\ref{thm:ineqchar} that we need to prove that any quantum correlation $p \in Q(4)$ satisfies $p(\pi(c),\pi(d)|c,d) - p(\pi(b),\pi(d)|a,d) + p(\pi(b),\pi(c)|a,b) \ge 0$ for any $\pi \in S_4$ and choice of $\{a,b,c,d\} = [4]$. Towards this we prove the following:

\begin{lemma}\label{lem:probs2triples}
Let $\mathcal{A}$ be a $C^*$-algebra with tracial state $\tau$, and $U = (u_{ij}) \in M_4(\mathcal{A})$ be a quantum permutation matrix. If $\pi \in S_4$ and $\{a,b,c,d\} = [4]$, then
\[\tau\left(u_{c\pi(c)}u_{d\pi(d)} - u_{a\pi(b)}u_{d\pi(d)} + u_{a\pi(b)}u_{b\pi(c)}\right) = \tau\left(u_{a\pi(a)}u_{b\pi(b)}u_{c\pi(c)} + u_{a\sigma(a)}u_{b\sigma(b)}u_{c\sigma(c)}\right),\]
where $\sigma = \pi(abcd)$.
\end{lemma}
\proof
We have that
\begin{align*}
\tau&\left(u_{c\pi(c)}u_{d\pi(d)} - u_{a\pi(b)}u_{d\pi(d)} + u_{a\pi(b)}u_{b\pi(c)}\right)\\
&= \tau\left(u_{c\pi(c)}u_{d\pi(d)} - u_{d\pi(d)}u_{a\pi(b)} + u_{b\pi(c)}u_{a\pi(b)}\right) \\
&= \tau\left(u_{c\pi(c)}u_{d\pi(d)} - \left(\sum_{x \in [4]}u_{x\pi(c)}\right)u_{d\pi(d)}u_{a\pi(b)} + u_{b\pi(c)}u_{a\pi(b)}\right) \\
&= \tau\left(u_{c\pi(c)}u_{d\pi(d)} - \left(u_{c\pi(c)} + u_{b\pi(c)}\right)u_{d\pi(d)}u_{a\pi(b)} + u_{b\pi(c)}u_{a\pi(b)}\right) \\
&= \tau\left(u_{c\pi(c)}u_{d\pi(d)}\left(\vec{1}-u_{a\pi(b)}\right) + u_{b\pi(c)}\left(\vec{1} - u_{d\pi(d)}\right)u_{a\pi(b)}\right) \\
&= \tau\left(u_{c\pi(c)}u_{d\pi(d)}\left(u_{a\pi(a)} + u_{a\pi(c)} + u_{a\pi(d)}\right)\right) \\
&\phantom{=} + \tau\left(u_{b\pi(c)}\left(u_{a\pi(d)} + u_{b\pi(d)} + u_{c\pi(d)}\right)u_{a\pi(b)}\right) \\
&= \tau\left(u_{c\pi(c)}u_{d\pi(d)}u_{a\pi(a)} + u_{b\pi(c)}u_{c\pi(d)}u_{a\pi(b)}\right) \\
&= \tau\left(u_{a\pi(a)}u_{b\pi(b)}u_{c\pi(c)} + u_{a\sigma(a)}u_{b\sigma(b)}u_{c\sigma(c)}\right),
\end{align*}
where the last equality uses Lemma~\ref{lem:triples} and the fact that $\pi(c) = \sigma(b)$, $\pi(d) = \sigma(c)$, and $\pi(b) = \sigma(a)$.\qeds


We are now able to show that the quantum bijective correlations on four points are always classical.

\begin{theorem}\label{thm:Q4=C4}
$Q(4) = L(4)$.
\end{theorem}
\proof
Let $p \in Q(4)$. Then there exists a $C^*$-algebra $\mathcal{A}$ with tracial state $\tau$ and quantum permutation matrix $U = (u_{ij}) \in M_4(\mathcal{A})$ such that $p(l,k|i,j) = \tau(u_{il}u_{jk})$ for all $i,j,l,k \in [4]$. Thus for any permutation $\pi \in S_4$ and $\{a,b,c,d\} = [4]$, we have that 
\begin{align*}
p(\pi(c),\pi(d)|c,d) - p(\pi(b),\pi(d)|a,d) + &p(\pi(b),\pi(c)|a,b)\\
&= \tau\left(u_{c\pi(c)}u_{d\pi(d)} - u_{a\pi(b)}u_{d\pi(d)} + u_{a\pi(b)}u_{b\pi(c)}\right).
\end{align*}
Thus by Theorem~\ref{thm:ineqchar} and Lemma~\ref{lem:probs2triples}, to prove $p \in L(4)$ it suffices to show that
\[\tau\left(u_{a\pi(a)}u_{b\pi(b)}u_{c\pi(c)} + u_{a\sigma(a)}u_{b\sigma(b)}u_{c\sigma(c)}\right) \ge 0,\]
where $\sigma = \pi(abcd)$. This we proceed to do. For convenience we will denote the quantity $\tau\left(u_{a\pi(a)}u_{b\pi(b)}u_{c\pi(c)} + u_{a\sigma(a)}u_{b\sigma(b)}u_{c\sigma(c)}\right)$ by $\gamma$. Let $(w,x,y,z) = (\pi(a),\pi(b), \pi(c), \pi(d))$. Since $\sigma=\pi(abcd)$, we have that
\begin{align*}
\gamma &= \tau(u_{a\pi(a)}u_{b\pi(b)}u_{c\pi(c)}+u_{a\pi(b)}u_{b\pi(c)}u_{c\pi(d)}) \\
&= \tau(u_{aw}u_{bx}u_{cy}+u_{ax}u_{by}u_{cz}).
\end{align*}
By Lemma \ref{lem:triples} and cyclicity of $\tau$, we get that
\begin{align}
\gamma &= \tau(u_{aw}u_{bx}u_{cy}+u_{ax}u_{by}u_{cz})\nonumber\\
&=\tau(u_{cy}u_{dz}u_{bx}+u_{dw}u_{by}u_{cz})\nonumber\\
&=\tau(u_{cy}u_{dz}u_{bx}(u_{dw}+u_{dz})+u_{dw}u_{by}u_{cz}(u_{by}+u_{bx})).\label{Ea}
\end{align}
We also have
\begin{align}
u_{bx}u_{dw}u_{cy}u_{dz} &= u_{bx}u_{dw}(\vec{1} - u_{ay}-u_{by})u_{dz}\nonumber\\
&=-u_{bx}u_{dw}(u_{ay}+u_{by})u_{dz}\nonumber\\
&= -u_{bx}u_{dw}\left((\vec{1} - u_{aw} - u_{ax} - u_{az}) + (\vec{1} - u_{bw} - u_{bx} - u_{bz})\right)u_{dz}\nonumber\\
&= u_{bx}u_{dw}(u_{ax}+u_{bx})u_{dz}\label{Eb}
\end{align}
and 
\begin{align}
\tau(u_{bx}u_{dw}u_{ax}u_{dz}) 
&= \tau(u_{bx}u_{dw}u_{ax}(\vec{1} - u_{cz}))\nonumber\\
&= \tau(- u_{bx}u_{dw}u_{ax}u_{cz}).\label{Ec}
\end{align}
Using the cyclicity of the tracial state, $u_{dw}u_{by}u_{cz}=u_{dw}u_{ax}u_{cz}$ (by Lemma \ref{lem:triples}) and Equations \eqref{Eb}, \eqref{Ec}, we get
\begin{align}
\tau(&u_{cy}u_{dz}u_{bx}u_{dw}+u_{dw}u_{by}u_{cz}u_{bx})\nonumber\\
&= \tau(u_{bx}u_{dw}u_{bx}u_{dz}+u_{bx}u_{dw}u_{ax}u_{dz} +u_{dw}u_{by}u_{cz}u_{bx})\nonumber\\
&= \tau(u_{bx}u_{dw}u_{bx}u_{dz}-u_{bx}u_{dw}u_{ax}u_{cz} +u_{dw}u_{by}u_{cz}u_{bx})\nonumber\\
&= \tau(u_{bx}u_{dw}u_{bx}u_{dz}-u_{bx}u_{dw}u_{ax}u_{cz} +u_{dw}u_{ax}u_{cz}u_{bx})\nonumber\\
&=\tau(u_{bx}u_{dw}u_{bx}u_{dz}).\label{Ed}
\end{align}
Now, the Equations \eqref{Ea} and \eqref{Ed} yield
\begin{align*}
\gamma &=\tau(u_{cy}u_{dz}u_{bx}(u_{dw}+u_{dz})+u_{dw}u_{by}u_{cz}(u_{by}+u_{bx}))\\
&=\tau(u_{cy}u_{dz}u_{bx}u_{dz}+u_{dw}u_{by}u_{cz}u_{by}\\
&\qquad+u_{cy}u_{dz}u_{bx}u_{dw}+u_{dw}u_{by}u_{cz}u_{bx})\\
&=\tau(u_{cy}u_{dz}u_{bx}u_{dz}+u_{dw}u_{by}u_{cz}u_{by}+u_{bx}u_{dw}u_{bx}u_{dz})\\
&=\tau(u_{cy}u_{dz}u_{bx}u_{dz}u_{cy}+u_{dw}u_{by}u_{cz}u_{by}u_{dw}\\
&\qquad+u_{dz}u_{bx}u_{dw}u_{bx}u_{dz}).
\end{align*}
Since all the elements appearing in the sum are positive, for example 
\begin{align*}
u_{cy}u_{dz}u_{bx}u_{dz}u_{cy}
&= u_{cy}u_{dz}u_{bx}u_{bx}u_{dz}u_{cy}\\
&= (u_{bx}u_{dz}u_{cy})^*u_{bx}u_{dz}u_{cy},
\end{align*}
 we deduce 
\begin{align*}
\gamma &= \tau(u_{a\pi(a)}u_{b\pi(b)}u_{c\pi(c)}+u_{a\sigma(a)}u_{b\sigma(b)}u_{c\sigma(c)})
\geq 0.
\end{align*}\qeds

The above shows that any quantum correlation $p \in Q(4)$ can in fact be written as a convex combination of classical deterministic correlations. Equivalently, there exists an entrywise nonnegative vector $\alpha = (\alpha_\pi)_{\pi \in S_4}$ such that $\sum_{\pi \in S_4} \alpha_\pi = 1$ and $M\alpha = \hat{p} := \Phi(p)$, where $M$ is the incidence matrix of $G_4$. The vector $\alpha$ can be found efficiently using linear programming, but if we know the quantum permutation matrix from which $p$ arises, we can compute $\alpha$ directly:
\begin{lemma}
Let $p \in Q(4)$ be a quantum correlation with corresponding quantum permutation matrix $U = (u_{ij}) \in M_4(\mathcal{A})$ for a $C^*$-algebra $\mathcal{A}$ with tracial state $\tau$. Then the vector $\alpha = (\alpha_\pi)$ defined as $\alpha_\pi = \tau(u_{a\pi(a)}u_{b\pi(b)}u_{c\pi(c)})$ for any distinct $a,b,c \in [4]$ satisfies $M\alpha = \hat{p}$.
\end{lemma}
\proof
Let $e = \{\pi,\sigma\}$ be an edge of $G_4$. Then WLOG we have that $\sigma = \pi(cd)$ for some distinct $c,d \in [4]$. Let $a,b \in [4]$ be such that $\{a,b,c,d\} = [4]$. This implies that $\pi$ and $\sigma$ agree on $a,b$ and that $\sigma(c) = \pi(d)$ and $\sigma(d) = \pi(c)$. Then
\begin{align*}
\left(M\alpha\right)_{e} &= \tau(u_{a\pi(a)}u_{b\pi(b)}u_{c\pi(c)}) + \tau(u_{a\sigma(a)}u_{b\sigma(b)}u_{c\sigma(c)}) \\
&= \tau\left(u_{a\pi(a)}u_{b\pi(b)}u_{c\pi(c)} + u_{a\pi(a)}u_{b\pi(b)}u_{c\pi(d)}\right) \\
&= \tau\left(u_{a\pi(a)}u_{b\pi(b)}(u_{c\pi(c)} + u_{c\pi(d)})\right) \\
&= \tau\left(u_{a\pi(a)}u_{b\pi(b)}(u_{c\pi(a)} + u_{c\pi(b)} + u_{c\pi(c)} + u_{c\pi(d)})\right) \\
&= \tau(u_{a\pi(a)}u_{b\pi(b)}) \\
&= p(\pi(a),\pi(b)|a,b) = \hat{p}(e).
\end{align*}\qeds

The above fits with the fact that classically (i.e., when the quantum permutation matrix $U$ has commutative entries), the value of $\tau(u_{a\pi(a)}u_{b\pi(b)}u_{c\pi(c)}) = \tau(u_{a\pi(a)}u_{b\pi(b)}u_{c\pi(c)}u_{d\pi(d)})$ really can be thought of as the probability with which you use the permutation $\pi$. However, it is possible that $\tau(u_{a\pi(a)}u_{b\pi(b)}u_{c\pi(c)}) < 0$ for an arbitrary quantum permutation, and thus to obtain a nonnegative vector $\beta$ satisfying $M\beta = \hat{p}$ from the $\alpha$ in the lemma above, we must sometimes add a multiple of the vector $\gamma$ as in the proof of Theorem~\ref{thm:ineqchar}. Is it an interesting question what vectors $\alpha = (\alpha_\pi)$ can be constructed via the method in the above Lemma? Similarly, what vectors $\alpha$ satisfy the property that $M\alpha \in B(4)$?

\subsection{$K_5$ does have nonlocal symmetry} \label{sec:K5}

In this subsection we will show that, in contrast to $K_4$, the graph $K_5$ does admit nonlocal symmetry. To do this we will give a specific $5 \times 5$ quantum permutation matrix that produces a nonclassical correlation. The quantum permutation matrix we use arises from a \emph{quantum Latin square}~\cite{qlatin}. A quantum Latin square is an $n \times n$ array of unit vectors from $\mathbb{C}^n$ such that each row and column forms an orthonormal basis. A quantum permutation matrix can then be constructed by taking as its $ij$-entry the projection onto the $ij$-entry of the quantum Latin square. We say that the quantum Latin square is \emph{classical} if all of the entries of this quantum permutation matrix commute, i.e., the set of vectors in any row/column differ from the set of vectors in any other row/column only by some scalar multiple (possibly a different scalar for each vector). The quantum Latin square we make use of is presented in Figure~\ref{fig:qlatin} (note that we have omitted a factor of $1/\sqrt{5}$ needed to make the vectors have unit norm).

\begin{figure}[!ht]
    \centering
    \[
    \begin{array}{|c|c|c|c|c|}
        \hline
         (1, 1, 1, 1, 1) & (1, \omega^4, \omega^3, \omega^2, \omega) & (1, \omega^3, \omega, \omega^4, \omega^2) & (1, \omega^2, \omega^4, \omega, \omega^3) & (1, \omega, \omega^2, \omega^3, \omega^4) \\
         \hline
         (1, \omega, \omega^2, \omega^4, \omega^3) & (1, 1, 1, \omega, \omega^4) & (1, \omega^4, \omega^3, \omega^3, 1) & (1, \omega^3, \omega, 1, \omega) & (1, \omega^2, \omega^4, \omega^2, \omega^2) \\
         \hline
         (1, \omega^2, \omega^4, \omega^3, \omega) & (1, \omega, \omega^2, 1, \omega^2) & (1, 1, 1, \omega^2, \omega^3) & (1, \omega^4, \omega^3, \omega^4, \omega^4) & (1, \omega^3, \omega, \omega, 1) \\
         \hline
         (1, \omega^3, \omega, \omega^2, \omega^4) & (1, \omega^2, \omega^4, \omega^4, 1) & (1, \omega, \omega^2, \omega, \omega) & (1, 1, 1, \omega^3, \omega^2) & (1, \omega^4, \omega^3, 1, \omega^3) \\
         \hline
         (1, \omega^4, \omega^3, \omega, \omega^2) & (1, \omega^3, \omega, \omega^3, \omega^3) & (1, \omega^2, \omega^4, 1, \omega^4) & (1, \omega, \omega^2, \omega^2, 1) & (1, 1, 1, \omega^4, \omega) \\
         \hline
    \end{array}
    \]
    \caption{Quantum Latin square of order five where $\omega = e^{2\pi i/5}$.}
    \label{fig:qlatin}
\end{figure}

To prove that a given correlation is nonlocal, it suffices to supply a separating hyperplane, i.e., to give an inequality satisfied by all classical correlations but not by the given correlation. We take this approach, but modified slightly. We will see that the correlation produced by our quantum permutation matrix exhibits a certain type of symmetry, and we will then give an inequality satisfied by any classical correlation that also has this symmetry. Finally, we will show that our quantum correlation does not satisfy this inequality, proving it is nonlocal.

In this subsection and going forward, for a set (or often a group) $X$, we will use the notation $B(X)$ to denote the set of bijective correlations $B(|X|)$ but where the underlying set is $X$ instead of $[|X|]$. We similarly use $L(X)$ and $Q(X)$ to denote the classical and quantum correlations in $B(X)$ respectively. The following definition introduces the relevant type of symmetry mentioned above. 

\begin{definition}
Suppose that $\Gamma$ is a finite group and that $p \in B(\Gamma)$. Then we say that $p$ is \emph{$\Gamma$-invariant} if $p(w,x|y,z)$ only depends on the values of $w^{-1}x$ and $y^{-1}z$.
\end{definition}

Equivalently, a correlation $p$ is $\Gamma$-invariant if $p(w,x|y,z) = p(aw,ax|by,bz)$ for all $a, b, w, x,$ $y, z \in \Gamma$. The correlation produced by the quantum Latin square shown in Figure~\ref{fig:qlatin} is $\mathbb{Z}_5$-invariant. Instead of proving this directly, in Section~\ref{sec:groupconstruct} we will present the construction of quantum Latin squares we used to produce the example in Figure~\ref{fig:qlatin}, and then we will show that this construction always results in group invariant correlations.

For any group $\Gamma$, there is a natural correlation that is $\Gamma$-invariant. For any $a \in \Gamma$ let $\sigma_a \in Sym(\Gamma)$\footnote{As with sets, we use $Sym(\Gamma)$ to denote the group of permutations of the elements of $\Gamma$. Importantly, this is not the group of \emph{automorphisms} of $\Gamma$.} denote the permutation given by $\sigma_a(x) = ax$. Recall that for any $\pi \in Sym(\Gamma)$ the correlation $p_\pi$ is the classical deterministic correlation such that $p_\pi(w,x|y,z)$ is 1 if and only if $\pi(y) = w$ and $\pi(z) = x$. Now define the correlation $p_\Gamma$ as
\begin{equation}\label{eq:pGamma1}
    p_\Gamma = \frac{1}{|\Gamma|} \sum_{a \in \Gamma}p_{\sigma_{a}}.
\end{equation}
Note that the above definition implies that $p_\Gamma$ is always classical. It is easy to see that
\[p_\Gamma(w,x|y,z) = \begin{cases}\frac{1}{|\Gamma|} & \text{if } w^{-1}x = y^{-1}z \\ 0 & \text{o.w.}\end{cases}\]
and thus $p_\Gamma$ is $\Gamma$-invariant. In the upcoming Lemma~\ref{lem:invariant} we make use of the \emph{composition} of correlations. Given correlations $p,p'\in B(X)$, their composition $p \circ p'$ is given by $p\circ p'(l,k|i,j) = \sum_{s,t} p(l,k|s,t)p'(s,t|i,j)$. Operationally, the composition $p \circ p'$ corresponds to first using $p'$ to obtain outputs in the isomorphism game, and then using these outputs as inputs for $p$ to obtain your final outputs that are sent back to the referee. We remark that $p_\Gamma$ is idempotent respect to this composition: $p_\Gamma \circ p_\Gamma = p_\Gamma$.

\begin{lemma}\label{lem:invariant}
Let $p \in B(\Gamma)$. Then the following are equivalent:
\begin{enumerate}
    \item $p$ is $\Gamma$-invariant;
    \item $p = p_\Gamma \circ p \circ p_\Gamma$;
    \item $p = p_{\sigma_b} \circ p \circ p_{\sigma_a}$ for all $a,b \in \Gamma$.
\end{enumerate}
\end{lemma}
\proof
We will show that $(1) \Rightarrow (2) \Rightarrow (3) \Rightarrow (1)$. First, suppose that $p$ is $\Gamma$-invariant. Then there exist numbers $D_{\alpha,\beta}$ for $\alpha, \beta \in \Gamma$ such that $p(w,x|y,z) = D_{\alpha, \beta}$ whenever $w^{-1}x = \alpha$ and $y^{-1}z = \beta$. Fixing $w,x,y,z \in \Gamma$ and letting $\alpha = w^{-1}x$ and $\beta = y^{-1}z$, we have that
\begin{align*}
    p_\Gamma \circ p \circ p_\Gamma (w,x|y,z) &= \sum_{q,r,s,t \in \Gamma} p_\Gamma(w,x|q,r)p(q,r|s,t)p_\Gamma(s,t|y,z) \\
    &= \frac{1}{|\Gamma|^2} \sum_{\substack{q,r,s,t \in \Gamma \\ q^{-1}r = \alpha, \ s^{-1}t = \beta}} p(q,r|s,t) \\
    &= \frac{1}{|\Gamma|^2} \sum_{\substack{q,r,s,t \in \Gamma \\ q^{-1}r = \alpha, \ s^{-1}t = \beta}} D_{\alpha,\beta} \\
    &= D_{\alpha,\beta} \\
    &= p(w,x|y,z).
\end{align*}
Thus we have proven the first implication.

Now suppose that $(2)$ holds. Since $p_{\sigma_a} \circ p_{\sigma_b} = p_{\sigma_{ab}}$, from Equation~\eqref{eq:pGamma1} we see that $p_\Gamma = p_\Gamma \circ p_{\sigma_a} = p_{\sigma_a}\circ p_\Gamma$ for all $a \in \Gamma$. Therefore,
\[p_{\sigma_b} \circ p \circ p_{\sigma_a} = p_{\sigma_b} \circ p_\Gamma \circ p \circ p_\Gamma \circ p_{\sigma_a} = p_\Gamma \circ p \circ p_\Gamma = p,\]
for all $a,b \in \Gamma$. Thus we have shown that $(2) \Rightarrow (3)$.

Lastly, suppose that $(3)$ holds. Fix $w,x,y,z \in \Gamma$ and let $\alpha = w^{-1}x$ and $\beta = y^{-1}z$. Then
\begin{align*}
    p(w,x|y,z) &= p_{\sigma_b} \circ p \circ p_{\sigma_a}(w,x|y,z) \\
    &= \sum_{q,r,s,t \in\Gamma} p_{\sigma_b}(w,x|q,r)p(q,r|s,t)p_{\sigma_a}(s,t|y,z) \\
    &= p(b^{-1}w,b^{-1}x|ay,az).
\end{align*}
Since this holds for all $a,b \in \Gamma$, we have that $p$ is $\Gamma$-invariant.\qeds

The above lemma reveals a justification for investigating $\Gamma$-invariant correlations. If $G$ and $H$ are two Cayley graphs for the group $\Gamma$ and $p$ is a correlation corresponding to a quantum isomorphism from $G$ to $H$, then $p_\Gamma \circ p \circ p_\Gamma$ is a $\Gamma$-invariant quantum correlation giving a quantum isomorphism from $G$ to $H$. Thus we can restrict our search for quantum isomorphisms between Cayley graphs for a group $\Gamma$ to $\Gamma$-invariant quantum correlations. Notably, the smallest known pair of non-isomorphic graphs that are quantum isomorphic are Cayley graphs for $S_4$~\cite{qiso1}. In fact one of them is the transposition graph $G_4$ used in Section~\ref{sec:K4}.

Using Lemma~\ref{lem:invariant} we prove the following:

\begin{lemma}\label{lem:unicoeffs}
Let $p \in L(\Gamma)$ be $\Gamma$-invariant. Then there exist nonnegative coefficients $\alpha_\pi$ for $\pi \in Sym(\Gamma)$ that satisfy
\begin{equation}\label{eq:unicoeffs}
    \alpha_\pi = \alpha_{\sigma_b \pi \sigma_a} \text{ for all } a,b \in \Gamma,
\end{equation}
and such that $\sum_\pi \alpha_\pi = 1$ and
\[p = \sum_{\pi \in Sym(\Gamma)} \alpha_\pi p_\pi.\]
\end{lemma}
\proof
Let $p \in L(\Gamma)$ be $\Gamma$-invariant. Since $p$ is classical, there exist nonnegative coefficients $\beta_\pi$ for $\pi \in Sym(\Gamma)$ such that $\sum_\pi \beta_\pi = 1$ and
\[p = \sum_{\pi \in Sym(\Gamma)} \beta_\pi p_\pi.\]
By Lemma~\ref{lem:invariant}, we have that $p = p_\Gamma \circ p \circ p_\Gamma$. Therefore,
\begin{align*}
    p &= p_\Gamma \circ p \circ p_\Gamma \\
    &= \frac{1}{|\Gamma|^2}\sum_{c,d \in \Gamma} p_{\sigma_d} \circ p \circ p_{\sigma_c} \\
    &= \frac{1}{|\Gamma|^2}\sum_{c,d \in \Gamma}\sum_{\pi \in Sym(\Gamma)} \beta_\pi p_{\sigma_d} \circ p_\pi \circ p_{\sigma_c} \\
    &= \frac{1}{|\Gamma|^2} \sum_{c,d \in \Gamma} \sum_{\pi \in Sym(\Gamma)} \beta_\pi p_{\sigma_d \pi \sigma_c} \\
    &= \frac{1}{|\Gamma|^2} \sum_{c,d \in \Gamma} \sum_{\pi \in Sym(\Gamma)} \beta_{\sigma_{d^{-1}}\pi\sigma_{c^{-1}}} p_{\pi} \\
    &= \frac{1}{|\Gamma|^2}\sum_{\pi \in Sym(\Gamma)} \left(\sum_{c,d \in \Gamma} \beta_{\sigma_{d^{-1}}\pi\sigma_{c^{-1}}}\right) p_{\pi}.
\end{align*}
Letting $\alpha_\pi = \frac{1}{|\Gamma|^2} \sum_{c,d} \beta_{\sigma_{d^{-1}}\pi\sigma_{c^{-1}}}$ we obtain the lemma statement.\qeds

We are now able to give an inequality satisfied by any classical $\mathbb{Z}_5$-invariant correlation. We remark that the inequality was found using a computer, but we give an analytical proof in the lemma below. 


\begin{lemma}\label{lem:C5ineq}
If $p \in L(\mathbb{Z}_5)$ is $\mathbb{Z}_5$-invariant, then it satisfies the following:
\begin{equation}\label{eq:C5ineq}
    2 p(0,2|0,2) + p(0,3|0,1) - p(0,4|0,4) \ge 0.
\end{equation}
\end{lemma}
\proof
As $p$ is a classical $\mathbb{Z}_5$-invariant correlation, we may write $p = \sum_\pi \alpha_\pi p_\pi$ where the $\alpha_\pi$ are nonnegative coefficients satisfying~\eqref{eq:unicoeffs}. For the purposes of the proof it improves clarity to write the coefficient $\alpha_\pi$ as $\alpha(\pi(0)\pi(1)\pi(2)\pi(3)\pi(4))$. Note that using this terminology, the condition of~\eqref{eq:unicoeffs} says that the coefficient $\alpha(abcde)$ is invariant under cyclically shifting the order of $a,b,c,d,e$ (this corresponds to right multiplication by $\sigma_a$) and under adding the same element of $\mathbb{Z}_5$ to each of $a,b,c,d,e$ (this corresponds to left multiplication by $\sigma_b$). In other words
\[\alpha(abcde) = \alpha(bcdea) = \alpha((a+1)(b+1)(c+1)(d+1)(e+1)).\]
To prove the inequality~\eqref{eq:C5ineq}, we first expand each term according to Equation~\eqref{eq:classical}:
\begin{align*}
    p(0,2|0,2) &= \alpha(01234) + \alpha(01243) + \alpha(03214) + \alpha(03241) + \alpha(04213) + \alpha(04231), \\
    p(0,3|0,1) &= \alpha(03124) + \alpha(03142) + \alpha(03214) + \alpha(03241) + \alpha(03412) + \alpha(03421), \\
    p(0,4|0,4) &= \alpha(01234) + \alpha(01324) + \alpha(02134) + \alpha(02314) + \alpha(03124) + \alpha(03214). \\
\end{align*}
Let $\gamma(p)$ denote the lefthand side of~\eqref{eq:C5ineq}. Using the above expressions we obtain
\begin{align}\label{eq:ineqexp}
     \gamma(p) &= \alpha(01234) + 2\alpha(01243) + 2\alpha(03214) + 2\alpha(03241) + 2\alpha(04213) + 2\alpha(04231)  \nonumber \\
     & \phantom{=} + \alpha(03142) + \alpha(03241) + \alpha(03412) + \alpha(03421) - \alpha(01324) - \alpha(02134) - \alpha(02314) \\
     &\ge 2\alpha(01243) + \alpha(03412) - \alpha(01324) - \alpha(02134) - \alpha(02314), \nonumber
\end{align}
where the inequality follows from the fact that we only removed nonnegative terms. We now aim to show that all of the terms in the final expression above can be cancelled through the use of condition~\eqref{eq:unicoeffs}.

For $\alpha(01324)$ we cyclically shift one position to the right and then add 1 everywhere to obtain
\[\alpha(01324) = \alpha(40132) = \alpha(01243).\]
For $\alpha(02134)$ we cyclically shift two positions to the right and then add 2 everywhere to obtain
\[\alpha(02134) = \alpha(34021) = \alpha(01243).\]
For $\alpha(02314)$ we cyclically shift two positions to the right and then subtract 1 everywhere to obtain
\[\alpha(02314) = \alpha(14023) = \alpha(03412).\]
Thus we have that 
\begin{align*}
    \gamma(p) &\ge 2\alpha(01243) + \alpha(03412) - \alpha(01324) - \alpha(02134) - \alpha(02314) \\
    &= 2\alpha(01243) + \alpha(03412) - \alpha(01243) - \alpha(01243) - \alpha(03412) \\
    &= 0.
\end{align*}\qeds

We are now able to show that $K_5$ admits nonlocal symmetry.

\begin{theorem}\label{thm:Q5neqC5}
There exists a quantum correlation $p \in Q(5)$ that is not contained in $L(5)$. In other words, the complete graph $K_5$ admits nonlocal symmetry.
\end{theorem}
\proof
Consider the quantum Latin square in Figure~\ref{fig:qlatin} to have its rows and columns indexed by $\mathbb{Z}_5$. Let $\psi_{ab}$ denote the (normalized) vector in the $ab$-entry of this quantum Latin square. Then $U = (u_{ab})$ such that $u_{ab} = \psi_{ab}\psi_{ab}^\dagger$ is a quantum permutation matrix which produces a quantum correlation
\[q(a,b|c,d) = \tr(u_{ca}u_{db}) = \frac{1}{5} |\langle \psi_{ca},\psi_{db}\rangle |^2.\]
As mentioned, we will not show that $q$ is $\mathbb{Z}_5$-invariant directly. Rather in Section~\ref{sec:groupconstruct} we will present the general construction used to produce the quantum Latin square of Figure~\ref{fig:qlatin}, and then show that the resulting correlation is always group invariant in Lemma~\ref{lem:Ginvariant}. Thus it is only left for us to show that $q$ does not satisfy the inequality in~\eqref{eq:C5ineq}.

Direct computation shows that
\begin{align*}
    q(0,2|0,2) &= \frac{1}{5} |\langle \psi_{00},\psi_{22}\rangle |^2 = \frac{1}{125} |3 + \omega^2 + \omega^3|^2= \frac{1}{25}\left(2 + \omega^2 + \omega^3\right) \\
    q(0,3|0,1) &= \frac{1}{5} |\langle \psi_{00},\psi_{13}\rangle |^2 = \frac{1}{125} |2 + 2\omega + \omega^3|^2 = \frac{1}{25} \\
    q(0,4|0,4) &= \frac{1}{5} |\langle \psi_{00},\psi_{44}\rangle |^2 = \frac{1}{125}|2 - \omega^2 - \omega^3|^2 = \frac{1}{25}\left(1 - \omega^2 - \omega^3\right).
\end{align*}
Thus
\begin{align*}
    2q(0,2|0,2) + q(0,3|0,1) - q(0,4|0,4) &= \frac{1}{25}\left(4 + 3\omega^2 + 3 \omega^3\right) \\
    &= \frac{1}{50}\left(5 - 3 \sqrt{5}\right) \\
    &< 0.
\end{align*}
Therefore, by Lemma~\ref{lem:C5ineq}, $q \in Q(5)$ is not a classical correlation and so $K_5$ does admit nonlocal symmetry.\qeds

\begin{cor}
The complete graph $K_n$ admits nonlocal symmetry if and only if $n \ge 5$.
\end{cor}

In terms of quantum groups, the above corollary says that $S_n^+$ exhibits nonlocal symmetry if and only if $n \ge 5$.

\section{Constructing nonlocal symmetry}\label{sec:construction}

In this section we present two constructions of quantum permutation matrices. The first construction, presented in Section~\ref{sec:groupconstruct}, produces a $\Gamma$-invariant correlation for any finite abelian group $\Gamma$ and permutation $\pi$ of the characters of $\Gamma$. This construction can produce nonlocal correlations, but this property is not guaranteed and it depends on the group $\Gamma$ and permutation $\pi$. In particular, for $\Gamma = \mathbb{Z}_5$ there are permutations $\pi$ such that the resulting correlation is nonlocal, and the correlation used in Theorem~\ref{thm:Q5neqC5} is constructed in this way. In Section~\ref{sec:disjointautos} we use classical automorphisms of a given graph $G$ to produce a quantum correlation that wins the $G$-automorphism game. This construction requires the existence of three pairwise ``disjoint" automorphisms of $G$, but nonlocality of the resulting correlation is guaranteed.

\subsection{Group invariant quantum symmetries}\label{sec:groupconstruct}

A \emph{flat unitary} matrix is a unitary matrix $U$ whose entries all have the same modulus (necessarily equal to $1/\sqrt{n}$ when $U$ is $n \times n$). Flat unitaries are also sometimes called \emph{complex Hadamard matrices}. If $U$ and $V$ are two flat unitaries, and we define $\psi_{ij}$ to be equal to $\sqrt{n}$ times the entrywise product of the $i^\text{th}$ column of $U$ with the $j^\text{th}$ column of $\overline{V}$ (the entrywise conjugate of $V$)\footnote{There is no real need to use $\overline{V}$ instead of $V$ here, but it makes the $U = V$ case a bit nicer since then $\psi_{ii}$ is the constant vector for all $i \in [n]$.}, then these vectors form a quantum Latin square, i.e., the subset of the $\psi_{ij}$ obtained by fixing one index and varying over the other is an orthonormal basis. To see this note that if $D_i$ is the diagonal matrix whose diagonal is equal to the $i^\text{th}$ column of $U$, then $\sqrt{n}D_i$ is unitary and thus so is $\sqrt{n}D_i\overline{V}$. The columns of the latter are simply the vectors $\psi_{ij}$ for $j \in [n]$ and thus they form an orthonormal basis by unitarity. A similar argument shows that for each $j \in [n]$ the vectors $\psi_{ij}$ for $i \in [n]$ form an orthonormal basis. Since the $\psi_{ij}$ form a quantum Latin square, the block matrix whose $ij$-block is the projection onto $\psi_{ij}$ is a quantum permutation matrix. We remark that this is similar to a construction of quantum Latin squares from~\cite{qlatin}, though there only one flat unitary is used.



We will make use of a specific instance of the above construction based on character tables of finite abelian groups. Recall that a character of a finite abelian group $\Gamma$ is a homomorphism $\chi : \Gamma \to S^1$, where $S^1$ is the group of complex numbers with modulus 1 under multiplication. The character table of $\Gamma$ has its columns indexed by $\Gamma$ and its rows by the characters of $\Gamma$, so that the $(\chi,a)$-entry is $\chi(a)$. Note that a finite abelian group has precisely as many characters as elements, and thus the character table is square. We will denote the character table of $\Gamma$ (which we think of as a matrix) by $\mathcal{C}_\Gamma$ or simply $\mathcal{C}$ when the group is clear from context. Since characters are homomorphisms, they satisfy $\chi(ab) = \chi(a)\chi(b)$ for all $a,b \in \Gamma$, and since they are homomorphisms to $S^1$, they satisfy $\overline{\chi(a)} = \chi(a^{-1})$. Letting $\circ$ denote entrywise product, it follows that
\begin{equation}\label{eq:columns}
    (\mathcal{C}e_a) \circ (\mathcal{C}e_b) = (\mathcal{C}e_{ab}), \text{ and } \  \overline{\mathcal{C}e_a} = \mathcal{C}e_{a^{-1}}.
\end{equation}
Note that it is well known that for any two characters $\chi,\chi'$ of a finite abelian group $\Gamma$ we have $\sum_{a \in \Gamma} \chi(a)\overline{\chi'(a)} = |\Gamma|\delta_{\chi,\chi'}$ and therefore the matrix $\frac{1}{\sqrt{|\Gamma|}}\mathcal{C}_\Gamma$ is a flat unitary. The characters themselves also form a group under pointwise multiplication, and this is known as the \emph{dual group}, denoted $\widehat{\Gamma}$. The dual group $\widehat{\Gamma}$ is isomorphic to $\Gamma$.

\begin{definition}\label{def:groupconstruct}
For a finite abelian group $\Gamma$ and permutation $\pi \in Sym(\widehat{\Gamma})$, let $\Psi^\Gamma_\pi = (\psi_{a,b})$ denote the $\Gamma \times \Gamma$ quantum Latin square where $\psi_{a,b} = \frac{1}{\sqrt{|\Gamma|}}(P^\pi\mathcal{C}e_a) \circ (\overline{\mathcal{C}e_b})$, where $P^\pi$ is the permutation matrix sending $e_\chi$ to $e_{\pi(\chi)}$. Furthermore, let $q^\Gamma_\pi$ be the quantum bijective correlation given by
\[q^\Gamma_\pi(a,b|c,d) = \frac{1}{|\Gamma|}|\langle\psi_{c,a},\psi_{d,b}\rangle|^2.\]
\end{definition}

We remark that when the group $\Gamma$ is clear from context we may just use $\Psi_\pi$ and $q_\pi$. Note that $\Psi^\Gamma_\pi$ is indeed a quantum Latin square by the discussion at the beginning of this subsection and the fact that both $\frac{1}{\sqrt{|\Gamma|}}\mathcal{C}$ and $\frac{1}{\sqrt{|\Gamma|}}P^\pi\mathcal{C}$ are flat unitaries. It then follows that $q^\Gamma_\pi$ is the quantum bijective correlation corresponding to the quantum permutation matrix $(\psi_{a,b}\psi_{a,b}^\dagger)_{a,b \in \Gamma}$. The quantum Latin square in Figure~\ref{fig:qlatin} is constructed by taking $\Gamma = \mathbb{Z}_5$, whose character table is given by the Fourier matrix $(\omega^{ij})_{i,j = 0,1,\ldots,4}$ where $\omega = e^{2\pi i/5}$. The permutation $\pi \in Sym(\{0,1,2,3,4\})$ used is the one swapping $3$ and $4$. We now show that $q^\Gamma_\pi$ is $\Gamma$-invariant.

\begin{lemma}\label{lem:Ginvariant}
Let $\Gamma$ be a finite abelian group, and $\pi \in Sym(\widehat{\Gamma})$. Then $q^\Gamma_\pi$ is $\Gamma$-invariant.
\end{lemma}
\proof
Let $\psi_{a,b}$ be the $a,b$ entry of $\Psi_\pi^\Gamma$. In order to show that $q_\pi^\Gamma$ is $\Gamma$-invariant, it suffices to show the stronger statement that $\langle \psi_{c,a},\psi_{d,b}\rangle$ only depends on $a^{-1}b$ and $c^{-1}d$.

Recall that $\psi_{c,a} = \frac{1}{\sqrt{|\Gamma|}}(P^\pi\mathcal{C}e_c) \circ (\overline{\mathcal{C}e_a})$. Let $\text{sum}(\psi)$ denote the sum of the entries of the vector $\psi$. Using Equation~\eqref{eq:columns}, which also holds for $P^\pi\mathcal{C}$, we have that
\begin{align}
\begin{split}\label{eq:Uvalues}
    |\Gamma|\langle \psi_{c,a},\psi_{d,b}\rangle &= |\Gamma|\psi_{c,a}^\dagger \psi_{d,b} \\
    &= \text{sum}\left(\overline{\psi_{c,a}} \circ \psi_{d,b}\right) \\
    &= \text{sum}\left(\overline{P^\pi\mathcal{C}e_c} \circ \mathcal{C}e_a \circ P^\pi\mathcal{C}e_d \circ \overline{\mathcal{C} e_b}\right) \\
    &= \text{sum}\left(\mathcal{C}e_a \circ  \mathcal{C} e_{b^{-1}} \circ P^\pi\mathcal{C}e_{c^{-1}} \circ P^\pi\mathcal{C}e_d\right) \\
    &= \text{sum}\left(\mathcal{C}e_{ab^{-1}} \circ P^\pi\mathcal{C}e_{c^{-1}d}\right) \\
    &= \text{sum}\left(\overline{\mathcal{C}e_{a^{-1}b}} \circ P^\pi\mathcal{C}e_{c^{-1}d} \right) \\
    &= (\mathcal{C}e_{a^{-1}b})^\dagger(P^\pi\mathcal{C}e_{c^{-1}d}) \\
    &= e_{a^{-1}b}^\dagger\left(\mathcal{C}^\dagger P^\pi \mathcal{C}\right)e_{c^{-1}d}
\end{split}
\end{align}
which clearly depends only on $a^{-1}b$ and $c^{-1}d$.\qeds


The above proves that the correlation used in Theorem~\ref{thm:Q5neqC5} to prove $K_5$ has nonlocal symmetry is indeed $\mathbb{Z}_5$-invariant, thus providing the missing piece from that proof.

By definition, a $\Gamma$-invariant bijective correlation $p$ has many redundant values. Thus we introduce a more efficient description of such correlations. Given a $\Gamma$-invariant correlation $p \in B(\Gamma)$, its \emph{characteristic matrix}, denoted $D^p$, is the $\Gamma \times \Gamma$ matrix whose $a,b$ entry is $|\Gamma|p(w,x|y,z)$ where $w^{-1}x = a$ and $y^{-1}z = b$. The factor of $|\Gamma|$ makes the matrix doubly stochastic. Indeed, its entries are nonnegative and its column sums are equal to
\[\sum_{a \in \Gamma} D^p_{a,b} = |\Gamma|\sum_{a \in \Gamma} p(w,wa|y,yb) = |\Gamma|\sum_{x \in \Gamma} p(w,x|y,yb) = |\Gamma|p(w|y) = |\Gamma|p(w,w|y,y) = 1,\]
where the last equality uses the fact that $p(w,w|y,y)$ must be independent of $w$ and thus equal to $1/|\Gamma|$. The proof that the row sums are 1 is similar. 

We remark that the characteristic matrix of the correlation $p_\Gamma$ defined in Equation~\eqref{eq:pGamma1} is the identity matrix. In the case that a $\Gamma$-invariant correlation $p$ comes from a quantum Latin square $\Psi = (\psi_{a,b})_{a,b \in \Gamma}$, we note that $D^p_{a,b} = |\langle\psi_{y,w},\psi_{z,x}\rangle|^2$ for $w^{-1}x = a$ and $y^{-1}z = b$. The following lemma gives a simple formula for the characteristic matrix of the correlation $q_\pi^\Gamma$, which we denote simply by $D^\pi$. It follows immediately from Equation~\eqref{eq:Uvalues}.

\begin{lemma}\label{lem:charmat}
Let $\Gamma$ be a finite abelian group and $\pi \in Sym(\widehat{\Gamma})$. Then
\[D^{\pi} = \frac{1}{|\Gamma|^2}\left(\mathcal{C}^\dagger P^\pi \mathcal{C}\right) \circ \left(\overline{\mathcal{C}^\dagger P^\pi \mathcal{C}}\right).\]\qed
\end{lemma}


The characteristic matrix is not simply a condensed encoding of a $\Gamma$-invariant correlation, it also behaves nicely with respect to composition:

\begin{lemma}\label{lem:Dcomp}
Let $p,p' \in B(\Gamma)$ be two $\Gamma$-invariant correlations. Then $p\circ p'$ is $\Gamma$-invariant and $D^{p\circ p'} = D^p D^{p'}$.
\end{lemma}
\proof
For $w^{-1}x = a$ and $y^{-1}z = b$, we have that
\begin{align*}
    p\circ p'(w,x|y,z) &= \sum_{s,t \in \Gamma} p(w,x|s,t)p'(s,t|y,z) \\
    &= \sum_{c \in \Gamma} \sum_{s,t:s^{-1}t = c} p(w,x|s,t)p'(s,t|y,z) \\
    &= \sum_{c \in \Gamma} \sum_{s,t:s^{-1}t = c} \frac{1}{|\Gamma|^2}D^p_{a,c}D^{p'}_{c,b} \\
    &= \sum_{c \in \Gamma} \frac{1}{|\Gamma|}D^p_{a,c}D^{p'}_{c,b} \\
    &= \frac{1}{|\Gamma|} \left(D^p D^{p'}\right)_{a,b}.
\end{align*}\qeds

The characteristic matrix of the quantum bijective correlation arising from the quantum Latin square in Figure~\ref{fig:qlatin} has characteristic matrix equal to the following, where $\varphi$ denotes the golden ratio $\frac{1+\sqrt{5}}{2}$:
\[
\frac{1}{5}\begin{pmatrix}
5 & 0 & 0 & 0 & 0 \\
0 & 1 + \varphi & 1 & 1 & 2 - \varphi \\
0 & 1 & 2 - \varphi & 1 + \varphi & 1 \\
0 & 1 & 1 + \varphi & 2 - \varphi & 1 \\
0 & 2 - \varphi & 1 & 1 & 1 + \varphi
\end{pmatrix}.
\]

The main result of this subsection is a characterization of precisely when the construction from Definition~\ref{def:groupconstruct} produces commutative quantum permutation matrices. To obtain this characterization, we make use of the following basic result in group theory:

\begin{lemma}\label{lem:charautos}
Let $\Gamma$ be an abelian finite group and $\widehat{\Gamma}$ its group of characters. For any $\sigma \in \aut(\Gamma)$ and $\chi \in \widehat{\Gamma}$, the function $\chi \circ \sigma$ is a character of $\Gamma$. Moreover, the map $\hat{\sigma}$ defined as $\hat{\sigma}(\chi) = \chi\circ \sigma$ is an automorphism of $\widehat{\Gamma}$ and the map $\sigma \mapsto \hat{\sigma}$ is an isomorphism of $\aut(\Gamma)$ and $\aut(\widehat{\Gamma})$.
\end{lemma}

For any $\sigma \in \aut(\Gamma)$ we will use $\hat{\sigma} \in \aut(\widehat{\Gamma})$ to denote the map defined in the above lemma.

\begin{theorem}\label{thm:permcorrs}
Let $\Gamma$ be a finite abelian group with character table $\mathcal{C}$ and $\pi \in Sym(\widehat{\Gamma})$. Also let $\psi_{a,b}$ denote the $a,b$ entry of the quantum Latin square $\Psi_\pi^\Gamma$. Then the following are equivalent:
\begin{enumerate}
    \item the entries of the quantum permutation matrix $\mathcal{P} = (\psi_{a,b}\psi_{a,b}^\dagger)_{a,b \in \Gamma}$ commute;
    \item there exist $\sigma \in \mathrm{Aut}(\Gamma)$ and $z \in \widehat{\Gamma}$ such that $\pi(x) = z\hat{\sigma}(x)$;
    \item the characteristic matrix $D^{\pi}$ is a permutation matrix.
\end{enumerate}
In the case that these hold, $D^{\pi} = P^{\sigma^{-1}} = (P^{\sigma})^T$.
\end{theorem}
\proof

Since two rank one projections commute if and only if they are orthogonal or equal, the projections $\psi_{a,b}\psi_{a,b}^\dagger$ commute if and only if every entry of $D^{\pi}$ is 0 or 1. Since $D^{\pi}$ is doubly stochastic, we have that (1) and (3) are equivalent.

Now we will show that $(3) \Rightarrow (2)$. By Lemma~\ref{lem:charmat} we have that $\frac{1}{|\Gamma|^2}\left(\mathcal{C}^\dagger P^\pi \mathcal{C}\right) \circ \left(\overline{\mathcal{C}^\dagger P^\pi \mathcal{C}}\right)$ is a permutation matrix. Thus every row/column of $\mathcal{C}^\dagger P^\pi \mathcal{C}$ has precisely one nonzero entry. As the columns of $\mathcal{C}$ form an orthogonal basis, this means that there exists a permutation $\sigma \in Sym(\Gamma)$ and function $z : \Gamma \to S^1$ such that
\begin{equation}\label{eq:relatecols}\mathcal{C}e_a = z(a)P^\pi\mathcal{C}e_{\sigma(a)} \text{ for all } a \in \Gamma.\end{equation}
We will show that $\sigma$ is an automorphism of $\Gamma$ and that $z$ is a character. Since $P^\pi u \circ P^\pi v = P^\pi(u \circ v)$ for any vectors $u$ and $v$, by Equations~\eqref{eq:columns} and~\eqref{eq:relatecols} we have that
\[z(ab) P^\pi\mathcal{C}e_{\sigma(ab)} = \mathcal{C}e_{ab} = \mathcal{C}e_a \circ \mathcal{C}e_b = z(a)P^\pi\mathcal{C}e_{\sigma(a)} \circ z(b)P^\pi\mathcal{C}e_{\sigma(b)} = z(a)z(b)P^\pi\mathcal{C}e_{\sigma(a)\sigma(b)}.\]
Since the columns of $\mathcal{C}$, and thus $P^\pi\mathcal{C}$, form an orthogonal basis, this implies that $\sigma(ab) = \sigma(a)\sigma(b)$ and it follows that $z(ab) = z(a)z(b)$ for all $a,b \in \Gamma$. In other words, $\sigma$ is an automorphism and $z$ is a character of $\Gamma$.

By Lemma~\ref{lem:charautos} $\hat{\sigma}$ is an automorphism of $\widehat{\Gamma}$ such that $x(\sigma(a)) = \hat{\sigma}(x)(a)$ for all $x \in \widehat{\Gamma}$ and $a \in \Gamma$. We now show that $\pi(x) = z\hat{\sigma}(x)$. For any $x \in \widehat{\Gamma}$ and $a \in \Gamma$, we have that
\begin{align*}
x(a) &= e^\dagger_x\mathcal{C}e_a \\
&= z(a)e_x^\dagger P^\pi \mathcal{C}e_{\sigma(a)} \\
&= z(a)e_{\pi^{-1}(x)}^\dagger \mathcal{C}e_{\sigma(a)} \\
&= z(a)\pi^{-1}(x)(\sigma(a)) \\
&= z(a)\hat{\sigma}(\pi^{-1}(x))(a) \\
&= (z\hat{\sigma}(\pi^{-1}(x)))(a).
\end{align*}
Thus we have that $x = z\hat{\sigma}(\pi^{-1}(x))$ for all $x \in \widehat{\Gamma}$. This implies that $\pi(x) = z\hat{\sigma}(x)$ as desired.

We will now show that if (2) holds then $D^{\pi} = P^{\sigma^{-1}}$. This additionally shows that $(2) \Rightarrow (3)$ and thus completes the proof. Assuming that $\pi(x) = z\hat{\sigma}(x)$ for all $x \in \widehat{\Gamma}$, we have that
\begin{align}
\begin{split}\label{eq:permcorrs}
    \left(\mathcal{C}^\dagger P^\pi \mathcal{C}\right)_{a,b} &= \sum_{x,y \in \widehat{\Gamma}}\overline{x(a)}(P^\pi)_{x,y}y(b) \\
    &= \sum_y \overline{\pi(y)(a)}y(b) \\
    &= \sum_y \overline{(z\hat{\sigma}(y))(a)}y(b) \\
    &= \sum_y \overline{z(a)y(\sigma(a))}y(b) \\
    &= \overline{z(a)}\sum_y \overline{y(\sigma(a))}y(b) \\
    &= \overline{z(a)} \left(\mathcal{C}^\dagger\mathcal{C}\right)_{\sigma(a),b} \\
    &= \overline{z(a)}|\Gamma|\delta_{\sigma(a),b} \\
    &= \overline{z(a)}|\Gamma|\left(P^{\sigma^{-1}}\right)_{a,b}.
\end{split}
\end{align}
Therefore, by Lemma~\ref{lem:charmat}, $D^{\pi} = \frac{1}{|\Gamma|^2}\left(\mathcal{C}^\dagger P^\pi\mathcal{C}\right) \circ \left(\overline{\mathcal{C}^\dagger P^\pi\mathcal{C}}\right) = P^{\sigma^{-1}}$.\qeds

\paragraph{Nonlocality of the correlations $q_\pi^\Gamma$.} The above shows that for a finite abelian group $\Gamma$, there are precisely $|\Gamma|\cdot |\aut(\Gamma)|$ elements $\pi \in Sym(\Gamma)$ such that the quantum Latin square $\Psi_\pi^\Gamma$ is classical, i.e., the projections onto its entries pairwise commute. In all other cases $\Psi_\pi^\Gamma$ is non-classical, but it may still be the case that the correlation $q_\pi^\Gamma$ that it produces is classical. This happens for instance when $\Gamma = \mathbb{Z}_4$, and $\pi$ is not of the form $\pi(x) = z\hat{\sigma}(x)$ for $\sigma \in \aut(\mathbb{Z}_4)$. Indeed, the resulting correlation $q_\pi^{\mathbb{Z}_4}$ must necessarily be classical by Theorem~\ref{thm:Q4=C4}. However, our computations indicate that the correlations $q_\pi^{\Gamma}$ seem to be nonlocal with high probability when $\Gamma$ is a cyclic group of sufficient size. In order to reduce the number of permutations $\pi \in Sym(\widehat{\Gamma})$ we need to check, we prove Lemma~\ref{lem:reducecorrs} below. Here, for $z \in \widehat{\Gamma}$ we use $\pi_z$ to denote the permutation of $\widehat{\Gamma}$ given by $\pi_z(x) = zx$, and we use $\eta$ be the permutation given by $\eta(x) = x^{-1}$

\begin{lemma}\label{lem:reducecorrs}
Let $\Gamma$ be a finite abelian group with character table $\mathcal{C}$ and let $\pi \in Sym(\widehat{\Gamma})$. Then we have the following:
\begin{enumerate}
    \item $D^{\pi\circ \pi_z} = D^{\pi} = D^{{\pi_z \circ \pi}}$ for all $z \in \widehat{\Gamma}$;
    \item $D^{{\pi \circ \hat{\sigma}}} = D^{\pi}D^{{\hat{\sigma}}} = D^{\pi}P^{\sigma^{-1}}$ and $D^{{\hat{\sigma} \circ \pi}} = D^{{\hat{\sigma}}}D^{\pi} = P^{\sigma^{-1}}D^{\pi}$ for all $\sigma \in \aut(\Gamma)$;
    \item $D^{{\eta \circ \pi \circ \eta}} = D^{\pi}$.
\end{enumerate}
\end{lemma}
\proof
We use Equation~\eqref{eq:permcorrs} to prove the statements. For (1), let $z \in \widehat{\Gamma}$ and define $Z$ to be the diagonal matrix whose $a,a$ entry is $\overline{z(a)}$. Then Equation~\eqref{eq:permcorrs} shows that $\mathcal{C}^\dagger P^{\pi_z}\mathcal{C} = |\Gamma|Z$. Therefore,
\[\mathcal{C}^\dagger P^{\pi \circ \pi_z}\mathcal{C} = \mathcal{C}^\dagger P^{\pi}P^{\pi_z}\mathcal{C} = \frac{1}{|\Gamma|} \mathcal{C}^\dagger P^{\pi}\mathcal{C}\mathcal{C}^\dagger P^{\pi_z}\mathcal{C} = \mathcal{C}^\dagger P^{\pi}\mathcal{C} Z.\]
Since $Z$ is a diagonal matrix with complex units on its diagonal, the above implies that each entry of $\mathcal{C}^\dagger P^{\pi \circ \pi_z}\mathcal{C}$ is the same as the corresponding entry of $\mathcal{C}^\dagger P^{\pi}\mathcal{C}$ up to multiplication of a complex unit. Since $\alpha\beta \cdot \overline{\alpha\beta} = \beta\overline{\beta}$ for any complex unit $\alpha$, it follows that
\[D^{{\pi\circ \pi_z}} = \mathcal{C}^\dagger P^{\pi \circ \pi_z}\mathcal{C} \circ \overline{\mathcal{C}^\dagger P^{\pi \circ \pi_z}\mathcal{C}} = \mathcal{C}^\dagger P^{\pi}\mathcal{C} \circ \overline{\mathcal{C}^\dagger P^{\pi}\mathcal{C}} = D^{\pi}.\]
The proof that $D^{{\pi_z \circ \pi}} = D^{\pi}$ is similar.

For (2), we again use Equation~\eqref{eq:permcorrs} which shows that $\mathcal{C}^\dagger P^{\hat{\sigma}} \mathcal{C} = |\Gamma|P^{\sigma^{-1}}$ for $\sigma \in \aut(\Gamma)$. Similarly to the above, this allows us to conclude that $\mathcal{C}^\dagger P^{\pi \circ \hat{\sigma}}\mathcal{C} = \mathcal{C}^\dagger P^{\pi}\mathcal{C} P^{\sigma^{-1}}$ and thus that
\[D^{{\pi \circ \hat{\sigma}}} = \mathcal{C}^\dagger P^{\pi \circ \hat{\sigma}}\mathcal{C} \circ \overline{\mathcal{C}^\dagger P^{\pi \circ \hat{\sigma}}\mathcal{C}} = \mathcal{C}^\dagger P^{\pi}\mathcal{C}P^{\sigma^{-1}} \circ \overline{\mathcal{C}^\dagger P^{\pi}\mathcal{C}}P^{\sigma^{-1}} = \left(\mathcal{C}^\dagger P^{\pi}\mathcal{C} \circ \overline{\mathcal{C}^\dagger P^{\pi}\mathcal{C}}\right)P^{\sigma^{-1}} = D^{\pi}P^{\sigma^{-1}}.\]
The proof that $D^{{\hat{\sigma} \circ \pi}} = P^{\sigma^{-1}}D^{\pi}$ is similar.

Finally, for (3), we use the fact that $P^{\eta}\mathcal{C} = \overline{\mathcal{C}}$, which follows from the fact that the inverse of a character $\chi$ of $\Gamma$ is also the pointwise conjugate of $\chi$. Note that since $\eta$ is an involution we have that $(P^\eta)^\dagger = P^\eta$ and thus $\mathcal{C}^\dagger P^\eta = \overline{\mathcal{C}^\dagger}$. Therefore we have that

\[\mathcal{C}^\dagger P^{\eta \circ \pi \circ \eta}\mathcal{C} = \mathcal{C}^\dagger P^\eta P^\pi P^\eta\mathcal{C} = \overline{\mathcal{C}^\dagger}P^\pi\overline{\mathcal{C}} = \overline{\mathcal{C}^\dagger P^\pi\mathcal{C}},\]
since $P^\pi$ is real. Thus we conclude that

\[D^{{\eta \circ \pi \circ \eta}} = \mathcal{C}^\dagger P^{{\eta \circ \pi \circ \eta}}\mathcal{C} \circ \overline{\mathcal{C}^\dagger P^{{\eta \circ \pi \circ \eta}}\mathcal{C}} = \overline{\mathcal{C}^\dagger P^\pi\mathcal{C}} \circ \mathcal{C}^\dagger P^\pi\mathcal{C} = D^{\pi}.\]\qeds

As a corollary, we obtain the following:

\begin{cor}\label{cor:nonlocalauts}
Let $\Gamma$ be a finite abelian group and let $\pi \in Sym(\widehat{\Gamma})$ and $\sigma \in \aut(\Gamma)$. Then the following are equivalent:
\begin{enumerate}
    \item $q_\pi$ is nonlocal;
    \item $q_{\pi \circ \hat{\sigma}}$ is nonlocal;
    \item $q_{\hat{\sigma} \circ \pi}$ is nonlocal.
\end{enumerate}
\end{cor}
\proof
We prove the contrapositive of $(1) \Leftrightarrow (2)$. Note that $q_{\hat{\sigma}}$ is a local correlation by Theorem~\ref{thm:permcorrs}. By Lemma~\ref{lem:reducecorrs} we have that $D^{{\pi \circ \hat{\sigma}}} = D^{\pi}D^{q_{\hat{\sigma}}}$. It then follows from Lemma~\ref{lem:Dcomp} that $q_{\pi \circ \hat{\sigma}} = q_\pi \circ q_{\hat{\sigma}}$. Since the composition of two classical correlations is classical, we have that $\neg (1) \Rightarrow \neg (2)$. The other direction follows from the same argument but with $\pi$ replaced with $\pi \circ \hat{\sigma}$ and $\sigma$ replaced with $\sigma^{-1}$. The proof for $(1) \Leftrightarrow (3)$ is identical.\qeds

Corollary~\ref{cor:nonlocalauts} and parts (1) and (3) of Lemma~\ref{lem:reducecorrs} allow us to reduce the number of permutations $\pi \in Sym(\widehat{\Gamma})$ that we need to check when looking for nonlocal correlations. For example, for $\Gamma = \mathbb{Z}_{10}$, instead of checking $10! = 3620000$ permutations, we only need to check 2375. We remark that for all of the groups we have checked, almost every case where $q_\pi = q_{\pi'}$ follows from some combination of the equivalences in items (1) and (3) of Lemma~\ref{lem:reducecorrs}. The exception is for eight pairs of equal correlations in the case $\Gamma = \mathbb{Z}_{10}$ and 256 pairs of equal correlations for $\Gamma = \mathbb{Z}_3^2$. For $\mathbb{Z}_{10}$ this is in a sense a single case, since all eight pairs are of the form $(\hat{\sigma}_1 \circ \pi \circ \hat{\sigma}_2, \hat{\sigma}_1 \circ \pi' \circ \hat{\sigma}_2)$ for some $\sigma_1, \sigma_2 \in \aut(\mathbb{Z}_{10})$ and a fixed pair $\pi,\pi' \in Sym(\widehat{\mathbb{Z}_{10}})$. Similarly, for $\mathbb{Z}_3^2$, the 256 pairs of permutations producing the same correlations is really just two distinct pairs.

\begin{table}[ht]
\begin{center}
\begin{tabu}{|c|[2pt]c|c|c|c|}
\hline
$\Gamma$ & \# distinct $q^\Gamma_\pi$  & \# classical $\Psi^\Gamma_\pi$ & \# local $q^\Gamma_\pi$ & \# nonlocal $q_\pi^\Gamma$ \\
\tabucline[2pt]{-}
$\mathbb{Z}_4$ & 3 & 2 & 3 & 0 \\
\hline
$\mathbb{Z}_5$ & 8 & 4 & 4 & 4 \\
\hline
$\mathbb{Z}_6$ & 20 & 2 & 5 & 15 \\
\hline
$\mathbb{Z}_7$ & 78 & 6 & 12 & 66 \\
\hline
$\mathbb{Z}_8$ & 380 & 4 & 10 & 370 \\
\hline
$\mathbb{Z}_9$ & 2438 & 6 & 14 & 2424 \\
\hline
$\mathbb{Z}_{10}$ & 18736 & 4 & 22 & 18714 \\
\hline
$\mathbb{Z}_2^2$ & 6 & 6 & 6 & 0 \\
\hline
$\mathbb{Z}_2^3$ & 924 & 168 & 924 & 0 \\
\hline
$\mathbb{Z}_2 \times \mathbb{Z}_4$ & 460 & 8 & 284 & 176 \\
\hline
$\mathbb{Z}_3^2$ & 2240 & 48 & 944 & 1296 \\
\hline

\end{tabu}
\end{center}
\vspace{0.1cm}
\caption{Summary of computations on the correlations $q_\pi^\Gamma$ for abelian groups $\Gamma$.}\label{tab:qpicomps}
\end{table}

Note that by Theorem~\ref{thm:permcorrs} the number of \emph{distinct} correlations $q_\pi$ corresponding to quantum Latin squares that are classical is $|\aut(\Gamma)|$. However, as we noted above, there may be other correlations $q_\pi$ that are local despite the corresponding quantum Latin square $\Psi_\pi$ being nonclassical. We have checked all such correlations for $\Gamma$ being an abelian group of order at most 10. For order at most 3, it is immediate that the quantum Latin squares $\Psi_\pi$ are classical since all quantum permutation matrices of order at most 3 are classical. For order between 4 and 10, Table~\ref{tab:qpicomps} summarizes our results. In the first column of this table is the given abelian group $\Gamma$. The second column contains the number of \emph{distinct} correlations $q_\pi^\Gamma$ for that group. 
The third column gives the number of permutations $\pi$ (giving distinct $q_\pi$) for which the quantum Latin square $\Psi_\pi^\Gamma$ is classical (this is just $|\aut(\Gamma)|$). The fourth column is the number of distinct local correlations $q_\pi^\Gamma$ for the given group $\Gamma$, and the fifth column gives the number of nonlocal correlations. 
It is clear from the table that, at least for cyclic groups, the proportion of permutations $\pi$ such that $q_\pi^\Gamma$ is nonlocal seems to quickly approach one as the size of the group grows. Interestingly, the non-cyclic groups seem to behave quite differently. Thus this property seems to be highly dependent on the group, not just its order.

\subsection{Three disjoint automorphisms criterion}\label{sec:disjointautos}

In \cite{foldedcubes}, the second author showed that a graph has quantum symmetry if its automorphism group contains a pair of non-trivial, disjoint automorphisms. This criterion is not enough for a graph to have nonlocal symmetry. For example, the complete graph $K_4$ has two disjoint automorphisms, but we know that it has no nonlocal symmetry by Theorem~\ref{thm:Q4=C4}. But we will see that graphs with \emph{three} disjoint automorphisms do have nonlocal symmetry.
The next definition can be found in \cite{foldedcubes}, we extend it to several permutations. 

\begin{definition}
Let $V$ be a finite set and let $\sigma_1, \dots, \sigma_n$ be permutations on $V$. We say that permutations $\sigma_1, \dots, \sigma_n$ are \emph{disjoint} if $\sigma_i(a) \neq a$ 
implies $\sigma_j(a)=a$ for all $j\neq i$.
\end{definition}

The following theorem shows that a graph has nonlocal symmetry if its automorphism group contains three non-trivial, disjoint automorphisms. 

\begin{theorem}\label{thm:threedisjaut}
Let $G$ be a finite graph. If its automorphism group $\aut(G)$ contains three non-trivial, disjoint automorphisms $\sigma_1, \sigma_2$ and $\sigma_3$, then $G$ has nonlocal symmetry. 
\end{theorem}

\proof
Let $\sigma_1, \sigma_2$ and $\sigma_3$ be non-trivial, disjoint automorphisms of $G$. Let $q_i$, $1\le i \le 3$, be the projection onto $v_i$ defined as
\[v_1 = \begin{pmatrix}1 \\ 0\end{pmatrix}, \quad v_2 = \begin{pmatrix}-1/2 \\ \sqrt{3}/2\end{pmatrix}, \quad v_3 = \begin{pmatrix}-1/2 \\ -\sqrt{3}/2\end{pmatrix}.\]
Note that $\|v_i\| = 1$ for $i = 1,2,3$ and $v_i^Tv_j = -1/2$ for $i \ne j$. Letting $\tr$ denote the normalized trace on $M_2(\mathbb{C})$, i.e.~$\tr(\vec{1}) = 1$, it follows that $\tr(q_i) = 1/2$ for $i = 1,2,3$ and $\tr(q_iq_j) = 1/8$ for $i \ne j$. We first show that the matrix 
\begin{align*}
U:=\left(\sum_{k=1}^3 P^{\sigma_k} \otimes q_k\right)+ \vec{1}\otimes (\vec{1}-q_1-q_2-q_3)\in \mathrm{M}_n(\C)\otimes \mathrm{M}_2(\C)\cong \mathrm{M}_n(\mathrm{M}_2(\C)))
\end{align*}
is a magic unitary that commutes with $(A_G\otimes \vec{1})$. We have
\begin{align}
u_{ij}&=\left(\sum_{k=1}^3 \delta_{i\sigma_k(j)} q_k\right)+ \delta_{ij} (\vec{1}-q_1-q_2-q_3)\nonumber\\
&=\begin{cases} q_k, &\text{if $\sigma_k(j)=i$ for some $k$ and $j \neq i$,}\\ \vec{1}-q_k, &\text{if $\sigma_k(j)\neq j$ for some $k$ and $i=j$,}\\ \delta_{ij}\vec{1},&\text{otherwise.}\label{eq:entries} \end{cases}
\end{align}
(Note that since $\sigma_1$, $\sigma_2$ and $\sigma_3$ are disjoint, $\sigma_k(j) \neq j$ implies $\sigma_l(j)=j$ for all $l\neq k$.) Thus, all entries of $U$ are projections. Furthermore
\begin{align*}
\sum_{s=1}^n u_{is}&=\sum_{s=1}^n\left(\left(\sum_{k=1}^3 \delta_{i\sigma_k(s)} q_k\right)+ \delta_{is} (\vec{1}-q_1-q_2-q_3)\right)\\
&=\left(\sum_{k=1}^3 \sum_{s=1}^n\delta_{i\sigma_k(s)} q_k\right)+ \sum_{s=1}^n\delta_{is} (\vec{1}-q_1-q_2-q_3)\\
&= q_1+q_2+q_3+\vec{1}-q_1-q_2-q_3\\
&=\vec{1}
\end{align*}
and similarly $\sum_{s=1}^n u_{si}=\vec{1}$. Also, $U$ commutes with $(A_G\otimes \vec{1})$, since
\begin{align*}
    U(A_G\otimes \vec{1})&=\left(\sum_{k=1}^3 P^{\sigma_k} A_G \otimes q_k\right)+ A_G\otimes (\vec{1}-q_1-q_2-q_3)\\
    &=\left(\sum_{k=1}^3 A_G P^{\sigma_k}  \otimes q_k\right)+ A_G\otimes (\vec{1}-q_1-q_2-q_3)\\
    &=(A_G\otimes \vec{1})U.
\end{align*}

Now, take the magic unitary $U$ and suppose that the associated correlation (using the tracial state $\mathrm{tr}$) is classical. Then, there exists a commutative magic unitary $V$ and a tracial state $\tau$ with $\tau(v_{ij}v_{kl})=\tr(u_{ij}u_{kl})$ for all $1 \le i,j,k,l \le n$. By \eqref{eq:entries} and since the automorphisms $\sigma_k$, $k=1,2,3$, are non-trivial, there exist $1\le a,b,c \le n$ such that $u_{\sigma_1(a)a}=q_1$, $u_{\sigma_2(b)b}=q_2$ and $u_{\sigma_3(c)c}=q_3$. Letting $p_1 = v_{\sigma_1(a)a}$, $p_2 = v_{\sigma_2(b)b}$, and $p_3 = v_{\sigma_3(c)c}$, this implies that $\tau(p_ip_j) = \tr(q_iq_j)$ which is $1/2$ if $i = j$ and $1/8$ otherwise. 

Since the $p_i$ all commute, their products are all projections. Moreover, it is straightforward to check that
\[r_1 = p_1, \quad r_2 = p_2 - p_1p_2, \quad r_3 = \vec{1} - p_1 - p_2 + p_1p_2\]
are pairwise orthogonal projections that sum to the identity. Using the known values of $\tau(p_ip_j)$, we also have that
\[\tau(r_1) = 1/2, \quad \tau(r_2) = 3/8, \quad \tau(r_3) = 1/8.\]
Now considering $p_3$, we have that
\[1/2 = \tau(p_3) = \tau((r_1 + r_2 + r_3)p_3) = \tau(r_1p_3) + \tau(r_2p_3) + \tau(r_3p_3).\]
Since $\tau(r_3) = 1/8$, we have that $\tau(r_3p_3) \le 1/8$, and so $\tau(r_1p_3) + \tau(r_2p_3) \ge 3/8$. This implies that at least one of
\[\tau(r_1p_3) \ge 3/16, \quad \tau(r_2p_3) \ge 3/16\]
holds. If the former holds, then $\tau(p_1p_3) \ge 3/16 > 1/8 = \tr(q_1q_3)$ which is a contradiction. If the latter holds, then $\tau(p_2p_3) \ge \tau((p_2 - p_1p_2)p_3) = \tau(r_2p_3) \ge 3/16 > 1/8$ which is again a contradiction. Therefore the correlation produced by $U$ must be non-classical.
\qeds



\section{Products}\label{sec:products}

In this section, we look at products of quantum permutation groups and study the correlations associated to their fundamental representations. We divide the section into three subsections dealing with the free product, the free wreath product and the tensor product of quantum permutation groups, respectively. Those products are defined as follows. 

\begin{prop}
Let $\mathbb{G}=(C(\mathbb{G}),U)$ and $\mathbb{H}=(C(\mathbb{H}), V)$ be quantum permutation groups with $U \in M_n(C(\mathbb{G}))$, $V \in M_m(C(\mathbb{H}))$.
\begin{itemize}
\item[(i)] 
 The pair $\mathbb{G}\times \mathbb{H}:= (C(\mathbb{G}) \otimes_{\mathrm{max}} C(\mathbb{H}), U \otimes V)$
is a quantum permutation group. Here $C(\mathbb{G}) \otimes_{\mathrm{max}} C(\mathbb{H})$ is the universal $C^*$-algebra with generators $u_{ij}v_{kl}$ such that $u_{ij}v_{kl}=v_{kl}u_{ij}$ for all $i,j,k,l$ and $u_{ij}$, $v_{kl}$ fulfill the relations of $C(\mathbb{G})$ and $C(\mathbb{H})$, respectively, where additionally $\vec{1}_{C(\mathbb{G})}=\vec{1}_{C(\mathbb{H})}$.
The quantum group $\mathbb{G} \times \mathbb{H}$ is called the \emph{tensor product} of $\mathbb{G}$ and $\mathbb{H}$, see \cite{Watensor}.
\item[(ii)] 
The pair $\mathbb{G}*\mathbb{H}:= (C(\mathbb{G})*C(\mathbb{H}), U \oplus V)$ is a quantum permutation group. Here $C(\mathbb{G})*C(\mathbb{H})$ is the universal $C^*$-algebra with generators $u_{ij}$ and $v_{kl}$ such that $u_{ij}$ and $v_{kl}$ fulfill the relations of $C(\mathbb{G})$ and $C(\mathbb{H})$, respectively, where additionally $\vec{1}_{C(\mathbb{G})}=\vec{1}_{C(\mathbb{H})}$. We call $\mathbb{G}*\mathbb{H}$ the \emph{free product} of $\mathbb{G}$ and $\mathbb{H}$, see \cite{Wafree}.
\item[(iii)]
The pair $\mathbb{G} \wr_* \mathbb{H}:=(C(\mathbb{G})*_w C(\mathbb{H}), W)$ is a quantum permutation group with $W=(w_{ia,jb})=(u_{ij}^{(a)} v_{ab})\in M_{nm}(C(\mathbb{G})*_w C(\mathbb{H}))$, where $U^{(a)} = (u_{ij}^{(a)})$ are copies of $U$. Here $C(\mathbb{G})*_w C(\mathbb{H})$ is the universal $C^*$-algebra generated by $u_{ij}^{(a)}, v_{ab}$ with the relations of the $C^*$-algebra $C(\mathbb{G})^{*m}*C(\mathbb{H})$, where additionally $u_{ij}^{(a)}v_{ab}= v_{ab}u_{ij}^{(a)}$. The quantum group $\mathbb{G} \wr_* \mathbb{H}$ is called the \emph{free wreath product}, see \cite{Bicfreewreath}.
\end{itemize}
\end{prop}

We use our findings for results on the nonlocal symmetry of graph unions and products. For example, we show that the disjoint union of two copies of a graph $G$ does not have nonlocal symmetry, if the graph $G$ has no quantum symmetry. 

\subsection{Free product}\label{subsec:freeproduct}

We start by showing that correlations associated to the direct sum of two commutative magic unitaries are classical. This yields that the free product of two quantum permutation groups with commutative $C^*$-algebras has  no nonlocal symmetry. Applying this result to graphs, we see that the disjoint union of two non-quantum isomorphic graphs has no nonlocal symmetry, if the graphs involved have no quantum symmetry. On the other hand, using the three disjoint automorphisms criterion, we get that the disjoint union of three graphs with non-trivial automorphism group does have nonlocal symmetry. 

\begin{lemma}\label{lem:directsum}
Let $\mathcal{A}$ be a $C^*$-algebra and $W = U \oplus V \in M_{m+n}(\mathcal{A})$ be a magic unitary where $U \in M_m(\mathcal{A})$ and $V \in M_n(\mathcal{A})$ are magic unitaries. Moreover assume that all of the entries of $U$ commute and likewise for $V$. Then for any tracial state $\tau$ on $\mathcal{A}$, the correlation $p(l,k|i,j) = \tau(w_{il}w_{jk})$ is classical.
\end{lemma}
\proof
For notational convenience we will let the entries of $U$ and $V$ be indexed by the sets $X$ and $Y$ respectively, and we will generally use $i,j,l,k$ for elements of $X$ and $a,b,c,d$ for elements of $Y$. Since the entries of $U$ commute, for each element $\pi \in Sym(X)$ we can define
\[u_{\pi} = \prod_{i \in X} u_{i\pi(i)},\]
without ambiguity (i.e., the order of the product does not matter). We can similarly define $v_\sigma$ for all $\sigma \in Sym(Y)$. Note that $u_\pi$ is a (possibly zero) projection for each $\pi \in Sym(X)$ and similarly for $v_\sigma$. To prove our result, we will first need to prove the following properties of the $u_\pi$ (resp.~$v_\sigma$) defined above:
\begin{enumerate}
\item $u_\pi u_{\pi'} = \delta_{\pi \pi'} u_\pi$;
\item $\sum_{\pi \in Sym(X)} u_\pi = \vec{1}$;
\item $u_{ij} = \sum_{\pi \in Sym(X): \pi(i) = j} u_\pi$.
\end{enumerate}
For (1), since $u_\pi$ is a projection we have that $u_\pi u_\pi = u_\pi$. On the other hand, if $\pi \ne \pi'$, then there exists $i \in X$ such that $\pi(i) \ne \pi'(i)$ and thus
\[u_\pi u_{\pi'} = \prod_{j \in X}u_{j\pi(j)} \prod_{l \in X} u_{l\pi'(l)}  = \prod_{j \in X} u_{j\pi(j)}u_{j\pi'(j)} = u_{i\pi(i)}u_{i\pi'(i)}\left(\prod_{j \ne i} u_{j\pi(j)}u_{j\pi'(j)}\right) = 0.\]
For (2), let $End(X)$ denote the set of all functions from $X$ to itself. We have that
\[\vec{1} = \prod_{i \in X} \sum_{j \in X} u_{ij} = \sum_{f \in End(X)} \prod_{i \in X} u_{if(i)} = \sum_{\pi \in Sym(X)} \prod_{i \in X} u_{i\pi(i)} = \sum_{\pi \in Sym(X)} u_\pi,\]
where the next to last equality comes from the fact that for any function $f \in End(X)$ that is not a permutation there are $i,j \in X$ such that $f(i) = f(j)$ and thus $u_{if(i)}u_{jf(j)} = 0$. Finally, for (3), we can use (2) to see that
\[u_{ij} = u_{ij} \sum_{\pi \in Sym(X)} u_{\pi} = \sum_{\pi \in Sym(X)} u_{ij} \prod_{l \in X} u_{l\pi(l)} = \sum_{\pi \in Sym(X) : \pi(i) = j} \prod_{l \in X} u_{l\pi(l)} = \sum_{\pi \in Sym(X) : \pi(i) = j} u_\pi\]
where the next to last equality follows since for any $\pi \in Sym(X)$ with $\pi(i) \ne j$, the product $u_{ij}\prod_l u_{l\pi(l)}$ is zero since $u_{ij}u_{i\pi(i)} = 0$ in this case. So we have proven that all of properties (1), (2), and (3) hold for the $u_\pi$ and analogously they hold for the $v_\sigma$.

Now let $\tau$ be a tracial state on $\mathcal{A}$ and $p(s,t|q,r) = \tau(w_{qs}w_{rt})$ as in the lemma statement. We will show that we can construct a classical correlation equal to $p$. For each $\pi \in Sym(X)$ and $\sigma \in Sym(Y)$ let $\pi \oplus \sigma$ denote the permutation of $Z := X \cup Y$ that maps $i \in X$ to $\pi(i)$ and $a \in Y$ to $\sigma(a)$. Further let $p_{\pi,\sigma}$ be the deterministic classical correlation corresponding to the permutation $\pi \oplus \sigma$. Define the correlation $q$ as the convex combination
\[q = \sum_{\pi,\sigma} \tau(u_\pi v_\sigma)p_{\pi,\sigma}.\]
Actually we must verify that this is in fact a convex combination. First, since $u_\pi$ and $v_\sigma$ are projections as noted above, we have that the coefficients are nonnegative. Also, using property (2) from above, we see that the sum of the coefficients is
\[\sum_{\pi, \sigma} \tau(u_\pi v_\sigma) = \tau\left(\left(\sum_\pi u_\pi\right)\left(\sum_\sigma v_\sigma\right)\right) = \tau(\vec{1}) = 1.\]
Thus $q$ is indeed a classical correlation.

Now we must show that $p = q$. Note that if $i \in X$ and $a \in Y$ then $p(a,t|i,r) = 0 = q(a,t|i,r)$. For $p$ this follows from the fact that $w_{ia} = 0$, and for $q$ it follows because the permutations $\pi \oplus \sigma$ for $\pi \in Sym(X)$, $\sigma \in Sym(Y)$ setwise fix $X$ and $Y$. We similarly have that $p$ and $q$ are both zero on all input/output pairs corresponding to mapping some element of $X$ to $Y$ or vice versa. So we can restrict to those values of $p$ and $q$ that correspond to fixing these sets. In the first case, we consider $p(l,k|i,j)$ and $q(l,k|i,j)$ for $l,k,i,j \in X$. In this case we have that $p(l,k|i,j) = \tau(u_{il}u_{jk})$ and
\[q(l,k|i,j) = \sum_{\pi: \pi(i) = l, \pi(j) = k} \sum_{\sigma \in Sym(Y)} \tau(u_{\pi}v_{\sigma}) = \sum_{\pi: \pi(i) = l, \pi(j) = k}  \tau(u_{\pi}).\]
But using (3) and then (1), we see that
\[u_{il}u_{jk} = \sum_{\pi : \pi(i) = l} u_{\pi} \sum_{\pi' : \pi'(j) = k} u_{\pi'} = \sum_{\pi : \pi(i) = l, \pi(j) = k} u_{\pi}.\]
Thus $q(l,k|i,j) = p(l,k|i,j)$ in this case as desired. It similarly holds that $p(b,d|a,c) = q(b,d|a,c)$ for all $a,b,c,d \in Y$.

Now consider $p(l,b|i,a)$ and $q(l,b|i,a)$ for $i,l \in X$ and $a,b \in Y$. We have that $p(l,b|i,a) = \tau(u_{il}v_{ab})$ and
\[q(l,b|i,a) = \sum_{\pi : \pi(i) = l} \sum_{\sigma : \sigma(a) = b} \tau(u_\pi v_\sigma) = \tau\left(\left(\sum_{\pi:\pi(i)=l} u_\pi\right)\left(\sum_{\sigma:\sigma(a)=b} v_\sigma\right)\right) = \tau(u_{il}v_{ab}).\]
Thus $p(l,b|i,a) = q(l,b|i,a)$ in this case, and analogously $p(b,l|a,i) = q(b,l|a,i)$. So we have shown that $p = q$ and thus $p$ is a classical correlation as desired.\qeds


We first apply our lemma to the free product of quantum permutation groups.

\begin{prop}\label{prop:freeprodqpg}
Let $\mathbb{G}$ and $\mathbb{H}$ be two quantum permutation groups with commutative $C^*$-algebras $C(\mathbb{G})$ and $C(\mathbb{H})$. Then, their free product $\mathbb{G}*\mathbb{H}$ has no nonlocal symmetry.
\end{prop}

\proof
This follows directly from Lemma \ref{lem:directsum}.
\qeds

It has the following consequence for nonlocal symmetries of graphs.

\begin{prop}\label{prop:nonisodisunion}
Let $G_1$ and $G_2$ be connected graphs that are not quantum isomorphic. If $G_1$ and $G_2$ have no quantum symmetry, then their disjoint union $G_1 \cup G_2$ has no nonlocal symmetry.
\end{prop}

\proof
By \cite[Theorem 4.5]{qperms}, we know that the fundamental representation of $\qut(G_1 \cup G_2)$ is of the form $W= \begin{pmatrix} U&0 \\ 0& V\end{pmatrix}$. It is easy to see that $U$ and $V$ are magic unitaries commuting with $A_{G_1}$ and $A_{G_2}$, respectively. Since $G_1$ and $G_2$ have no quantum symmetry, we know that the entries of $U$ commute and likewise the entries of $V$. Any quantum correlation $p$ for the $G_1 \cup G_2$ automorphism game has the form $p(l,k|i,j) = \tau(w_{il}w_{jk})$ for a tracial state $\tau$ on $C(\qut(G_1 \cup G_2))$. By Lemma~\ref{lem:directsum}, any such correlation $p$ is classical and thus we are done.
\qeds

Given Lemma \ref{lem:directsum}, it is natural to ask if it can be extended to the direct sum of 3 or more (individually) commuting magic unitaries. It turns out that the answer is no, this follows from Theorem \ref{thm:threedisjaut}.

\begin{cor}\label{cor:graphsdisjunion}
Let $G_1$, $G_2$ and $G_3$ be finite graphs with non-trivial automorphism group. Then $G_1\cup G_2 \cup G_3$ does have nonlocal symmetry.
\end{cor}

\proof
We get three non-trivial, disjoint automorphisms of $G_1\cup G_2 \cup G_3$ by choosing any non-trivial automorphism of each graph that fix the other ones (those exists since we assumed that the graphs have non-trivial automorphism group). Then Theorem \ref{thm:threedisjaut} yields the assertion. 
\qeds

\subsection{Free wreath product}\label{subsec:wreathproduct}

In this subsection, we first show how correlations associated to a certain magic unitary look. We relate this to the free wreath product of quantum permutation groups. We will then see that the disjoint union of two copies of a graph does not have nonlocal symmetry if the graph has no quantum symmetry. 

\begin{lemma}\label{lem:freewreath}
Let $U = \begin{pmatrix}U^{11} & U^{12} \\ U^{21} & U^{22}\end{pmatrix}$ be a magic unitary with entries from a $C^*$-algebra $\mathcal{A}$. Let $\{X, Y\}$ be the partition of the indices of the rows/columns of $U$ that corresponds to the given partition of $U$. Suppose that the entries of $U^{11}$ and $U^{22}$ are orthogonal to the entries of $U^{12}$ and $U^{21}$ and let $\mathcal{A}'$ and $\mathcal{A}''$ be the $C^*$-subalgebras of $\mathcal{A}$ generated by the entries of these pairs of submatrices respectively. Then $U' = \begin{pmatrix}U^{11} & 0 \\ 0 & U^{22}\end{pmatrix}$ and $U'' = \begin{pmatrix}0 & U^{12} \\ U^{21} & 0\end{pmatrix}$ are magic unitaries (or zero matrices) over $\mathcal{A}'$ and $\mathcal{A''}$ respectively, and for any tracial state $\tau$ on $\mathcal{A}$, there exist tracial states $\tau'$ and $\tau''$ on $\mathcal{A}'$ and $\mathcal{A}''$ respectively and $0 \le \gamma \le 1$ such that
\begin{equation}\label{eq:tau}
\tau(u_{il}u_{jk}) = \begin{cases} \gamma \tau'(u_{il}u_{jk}) & \text{if } (i,j),(l,k) \in \left(X \times X\right) \cup \left(Y \times Y\right) \\ (1-\gamma) \tau''(u_{il}u_{jk}) & \text{if } (i,j),(l,k) \in (X \times Y) \cup (Y \times X) \\ 0 & \text{otherwise.}\end{cases}
\end{equation}
It follows that if $p,p',p''$ are the correlations corresponding the the magic unitaries $U,U',U''$ and tracial states $\tau,\tau',\tau''$, then $p = \gamma p' + (1-\gamma)p''$.
\end{lemma}
\proof
From the proof of Theorem 4.4 from~\cite{qperms}, it follows that there exists a projection $\rho \in \mathcal{A}$ such that the rows and columns of $U^{11}$ and $U^{22}$ each have sum equal to $\rho$. Thus the rows and columns of $U^{12}$ and $U^{21}$ each have sum equal to $\vec{1} - \rho$. Note that $\rho$ and $\vec{1}-\rho$ are the identity in the subalgebras $\mathcal{A}'$ and $\mathcal{A}''$ respectively, and thus we have proven that $U'$ and $U''$ are magic unitaries (or one is a zero matrix if $\rho = 0$ or $\rho = \vec{1}$). Let $\gamma = \tau(\rho)$ and define $\tau' = (1/\gamma) \tau |_{\mathcal{A}'}$ and $\tau'' = (1/(1-\gamma)) \tau |_{\mathcal{A}''}$ (if $\gamma = 0$ then $\tau'$ is defined as the zero function and analogously for $\tau''$ if $\gamma = 1$). It is easy to see that, unless they are the zero function, these are tracial states on $\mathcal{A}'$ and $\mathcal{A}''$ respectively. Also, Equation~\eqref{eq:tau} clearly holds.

Now if $p,p',p''$ are the correlations corresponding to the magic unitaries $U,U',U''$ and tracial states $\tau,\tau',\tau''$, then we have that
\begin{align*}
p'(l,k|i,j) &= \begin{cases} \tau'(u_{il}u_{jk}) & \text{if } (i,j),(l,k) \in (X \times X) \cup (Y \times Y) \\ 0 & \text{otherwise,}\end{cases} \\
p''(l,k|i,j) &= \begin{cases} \tau''(u_{il}u_{jk}) & \text{if } (i,j),(l,k) \in (X \times Y) \cup (Y \times X) \\ 0 & \text{otherwise.}\end{cases}
\end{align*}
It then follows that $p = \gamma p' + (1-\gamma)p''$.\qeds

\begin{prop}
Let $\mathbb{G}$ be a quantum permutation group with commutative $C^*$-algebra. Then, the free wreath product $\mathbb{G}\wr_* \Z_2$ has no nonlocal symmetry.
\end{prop}

\proof
The fundamental representation of $\mathbb{G}\wr_* \Z_2$ is of the form $W=(w_{ia,jb})=(u_{ij}^{(a)}v_{ab})$, where $U$ is the fundamental representation of $\mathbb{G}$ and $V=\begin{pmatrix}q &1-q\\1-q&q\end{pmatrix}$, where $q$ is a projection.We order the tuples $ia$, $1\le i\le n$, $1\le a \le 2$ like this: $11<21<\dots<n1<12<22<\dots<n2$. In this ordering, we partition the indices into $X=\{11,\dots, n1\}$ and $Y=\{12,\dots, n2\}$. Using this partition, we write $W=\begin{pmatrix} W^{11}&W^{12}\\W^{21}&W^{22}\end{pmatrix}$. The entries of $W^{11}$ and $W^{22}$ are orthogonal to the entries of $W^{12}$ and $W^{21}$ (every entry of $W^{11}$ and $W^{22}$ can be written as $qu_{ij}^{(a)}$ and every entry of $W^{12}$ and $W^{21}$ as $(1-q)u_{ij}^{(a)}$). Write
$W' = \begin{pmatrix}W^{11} & 0 \\ 0 & W^{22}\end{pmatrix}$ and $W'' = \begin{pmatrix}0 &W^{12} \\ W^{21} & 0\end{pmatrix}$. Since we assumed that $C(\mathbb{G})$ is commutative, we also have that the entries of $W^{11}$ commute and likewise for $W^{22}$, $W^{12}$ and $W^{21}$.
We deduce from Lemma \ref{lem:directsum} that the correlations associated to the magic unitaries $W'$ and $W''$ are classical. Lemma \ref{lem:freewreath} yields that any correlation associated to $W$ is a convex combination of classical correlations, i.e., is classical itself.
\qeds

In the following, we denote by $nG$ the disjoint union of $n$ copies of the graph $G$. 

\begin{prop}(\cite[Theorem 6.1]{BanBicfr})
Let $G$ be a connected graph. Let $W$, $U$ and $V$ be the fundamental representations of $\qut(2G)$, $\qut(G)$ and $\Z_2$, respectively. Then we have the $*$-isomorphism
\begin{align*}
\varphi:C(\qut(2G)) &\to C(\qut(G))*_{w} C(\Z_2)\\w_{ia,jb}&\mapsto u_{ij}^{(a)}v_{ab}.
\end{align*}
\end{prop}

\begin{prop}\label{prop:2disjunion}
Let $G$ be a connected graph that has no quantum symmetry. Then $2G$ has no nonlocal symmetry.
\end{prop}

\proof
Let $U= \begin{pmatrix} U^{11}& U^{12} \\ U^{21}& U^{22}\end{pmatrix}$ (partitioned into the vertex sets of the copies of $G$) be the fundamental representation of $\qut(2G)$. By the $*$-isomorphism of the previous proposition, we see that the entries of $U^{11}$ and $U^{22}$ are orthogonal to the entries of $U^{12}$ and $U^{21}$. Let $\mathcal{A'}$ and $\mathcal{A''}$ be the $C^*$-subalgebras of $C(\qut(2G))$ as in Lemma \ref{lem:freewreath}, and let $U'$ and $U''$ be as they are defined there. By Lemma \ref{lem:freewreath}, the matrices $U'$ and $U''$ are magic unitaries in $\mathcal{A'}$ and $\mathcal{A''}$, respectively. It follows that $U^{11}$, $U^{22}$ are magic unitaries in $\mathcal{A'}$ and $U^{12}$, $U^{21}$ are magic unitaries in $\mathcal{A''}$. 
Let $A_G$ be the adjacency matrix of $G$. Then $A_{2G}=\begin{pmatrix} A_G&0\\ 0 & A_G\end{pmatrix}$ and thus
\begin{align*}
\begin{pmatrix}U^{11} A_G& U^{12}  A_G\\ U^{21} A_G& U^{22} A_G \end{pmatrix}&=\begin{pmatrix}U^{11}& U^{12} \\ U^{21}& U^{22} \end{pmatrix}\begin{pmatrix} A_G&0\\ 0 & A_G\end{pmatrix}\\
&=\begin{pmatrix} A_G&0\\ 0 & A_G\end{pmatrix}\begin{pmatrix} U^{11}& U^{12} \\ U^{21}& U^{22}\end{pmatrix}\\
&=\begin{pmatrix} A_GU^{11}& A_G U^{12} \\  A_GU^{21}&  A_GU^{22} \end{pmatrix}.
\end{align*}
We deduce that $U^{11}$, $U^{22}$, $U^{12}$ and $U^{21}$ are magic unitaries that commute with $A_G$. Since we assumed that $G$ has no quantum symmetry, we get that the entries of $U^{11}$ commute and likewise for $U^{22}$, $U^{12}$ and $U^{21}$. By Lemma \ref{lem:directsum}, we obtain that the correlations associated to the magic unitaries $U'$ and $U''$ are classical. Lemma \ref{lem:freewreath} then yields that any correlation associated to $U$ is classical. 
\qeds 

The above proposition shows that the disjoint union of isomorphic graphs has no nonlocal symmetry, as long as the individual graphs have no quantum symmetry. Recall that Proposition~\ref{prop:nonisodisunion} shows the same holds for the disjoint union of graphs that are not quantum isomorphic. This leaves the case of the disjoint union of two non-isomorphic but quantum isomorphic graphs. In this case one always has nonlocal symmetry, since the quantum isomorphism between the two graphs gives a nonlocal symmetry swapping the two components of the disjoint union.

\begin{cor}
Let $G$ be a graph with $\aut(G) \neq \{e\}$. Then $nG$, $n \geq 3$, does have nonlocal symmetry.
\end{cor}

\proof
This directly follows from Corollary \ref{cor:graphsdisjunion}.
\qeds

\subsection{Tensor product}\label{sec:tensorproduct}
As in Sections \ref{subsec:freeproduct} and \ref{subsec:wreathproduct}, we first work on special magic unitaries and the correlations they produce. We then obtain that the tensor product of a quantum permutation group $\mathbb{G}$ and $\mathbb{H}$ has no nonlocal symmetry if $\mathbb{G}$ has no nonlocal symmetry and $\mathbb{H}$ has no quantum symmetry. We apply our results to the cartesian product and the tensor product of graphs. 

In the following lemma, we use the fact that for an $n \times n$ commutative magic unitary $V = (v_{ab})$, we can define projections $v_\pi = \prod_{a=1}^n v_{a\pi(a)}$ for all $\pi \in S_n$. Recall that this was also used in the proof of Lemma~\ref{lem:directsum}.

\begin{lemma}\label{lem:tensorproduct}
Let $W=U\otimes V\in M_{mn}(\mathcal{A})$ be a magic unitary with entries from a $C^*$-algebra $\mathcal{A}$, where $U\in M_m(\mathcal{A})$ is a magic unitary and $V \in M_n(\mathcal{A})$ is a commutative magic unitary. For each $\pi \in S_n$ let $v_\pi = \prod_{a=1}^n v_{a\pi(a)}$, and let $\mathcal{A}_\pi$ be the $C^*$-subalgebra of $\mathcal{A}$ generated by $u_{ij}v_\pi$, $1\le i,j\le m$. Then the matrix $(u_{ij}v_\pi)_{ij}$ is a magic unitary over $\mathcal{A}_\pi$ for each $\pi \in S_n$. Furthermore, for every tracial state $\tau$ on $\mathcal{A}$, there exist tracial states (or zero functions) $\tau_\pi$ on $\mathcal{A}_\pi$ and $0\le \gamma_\pi \le 1$ for all $\pi \in S_n$ such that $\sum_{\pi \in S_n} \gamma_\pi = 1$ and
\begin{align*}
\tau(w_{ia,jb}w_{kc,ld})&= \sum_{\pi: \pi(a)=b, \pi(c)=d} \gamma_\pi\tau_\pi(u_{ij}v_\pi u_{kl}v_\pi).
\end{align*}
It follows that if $p$ is the correlation corresponding to $W$ and $p^\pi$ is the correlation corresponding to $(u_{ij}v_\pi)_{ij}$, then
\begin{align*}p(jb,ld|ia,kc)= \sum_{\pi:\pi(a)=b,\pi(c)=d}\gamma_\pi p^\pi(j,l|i,k).\end{align*}
\end{lemma}
\proof
Recall from the proof of Lemma~\ref{lem:directsum} that $v_{ab} = \sum_{\pi: \pi(a) = b} v_{\pi}$ and $v_{\pi}v_{\pi'} = \delta_{\pi\pi'}v_{\pi}$. Moreover, since $u_{ij}v_{ab} = v_{ab}u_{ij}$ for all $i,j,a,b$, we have that
\begin{align*}
    \tau(w_{ia,jb}w_{kc,ld}) &= \tau(u_{ij}v_{ab}u_{kl}v_{cd}) \\
    &= \tau(u_{ij}u_{kl}v_{ab}v_{cd}) \\
    &= \sum_{\pi:\pi(a) = b} \sum_{\pi':\pi'(c)=d} \tau(u_{ij}u_{kl}v_\pi v_{\pi'}) \\
    &= \sum_{\pi: \pi(a)=b,\pi(c)=d}\tau(u_{ij}u_{kl}v_\pi) \\
    &= \sum_{\pi: \pi(a)=b,\pi(c)=d}\tau(u_{ij}v_\pi u_{kl}v_\pi)
\end{align*}
Now let $\gamma_\pi = \tau(v_\pi) \ge 0$ and define $\tau_\pi = (1/\gamma_\pi)\tau|_{\mathcal{A}_\pi}$ if $\gamma_\pi \ne 0$ and define it as the zero function on $\mathcal{A}_\pi$ otherwise. Since each $v_\pi$ is a central projection, these are all tracial states on their respective algebras unless they are the zero function. Recall also from the proof of Lemma~\ref{lem:directsum} that $\sum_\pi v_\pi = \one$ and thus $\sum_\pi \gamma_\pi = 1$ as desired. Continuing the above sequence of equations, we obtain
\[\tau(w_{ia,jb}w_{kc,ld}) = \sum_{\pi:\pi(a)=b,\pi(c)=d}\gamma_\pi \tau_\pi(u_{ij}v_\pi u_{kl}v_\pi).\]
Note that since $v_\pi$ is the identity in $\mathcal{A}_\pi$, we have that $(u_{ij}v_\pi)$ is indeed a magic unitary and thus $p^\pi(j,l|i,k) := \tau_\pi(u_{ij}v_\pi u_{kl}v_\pi)$ is indeed a correlation (unless $\tau(v_\pi) = 0$ in which case it is the zero function). Moreover, from the above we have that
\[p(jb,ld|ia,kc)= \sum_{\pi:\pi(a)=b,\pi(c)=d}\gamma_\pi p^\pi(j,l|i,k),\]
as desired.\qeds

We apply our lemma to the tensor product of two quantum permutation groups without nonlocal symmetry, at least one of which is classical.

\begin{prop}\label{prop:qpgtensor}
Let $\mathbb{G}$ and $\mathbb{H}$ be quantum permutation groups such that $\mathbb{G}$ has no nonlocal symmetry and $\mathbb{H}$ has a commutative $C^*$-algebra. Then $\mathbb{G}\times \mathbb{H}$ has no nonlocal symmetry.
\end{prop}
\proof
Let $p$ be the correlation associated to the fundamental representation $U\otimes V$ of $\mathbb{G} \times \mathbb{H}$ and some tracial state $\tau$. By Lemma \ref{lem:tensorproduct}, we know that 
\begin{equation}\label{eq:tenscorr}
    p(jb,ld|ia,kc) = \sum_{\pi:\pi(a)=b,\pi(c)=d}\gamma_\pi p^\pi(j,l|i,k)
\end{equation}
where $p^\pi$ is the correlation (or zero function) corresponding to the magic unitary $(u_{ij}v_\pi)_{ij}$ and tracial state $\tau_\pi$. Let us define $C \subseteq S_n$ to be the set of permutations $\pi$ such that $p^\pi$ is not the zero function, i.e., such that $\gamma_\pi := \tau(v_\pi) \ne 0$, and let $C_{ab,cd}$ contain the permutations $\pi \in C$ such that $\pi(a) = b$ and $\pi(c) = d$. Since $\mathbb{G}$ has no nonlocal symmetry, if $\pi \in C$ then $p^\pi$ can be written as
\[p^\pi(j,l|i,k) = \sum_{\sigma \in S_m} \alpha^\pi_\sigma p_\sigma(j,l|i,k),\]
for some coefficients $\alpha^\pi_\sigma \ge 0$ such that $\sum_{\sigma \in S_m} \alpha^\pi_\sigma = 1$, where $p_\sigma$ denotes the deterministic classical correlation corresponding to the permutation $\sigma$ as defined in Equation~\eqref{eq:detcorrs}. By Equation~\eqref{eq:tenscorr}, we get
\[p(jb,ld|ia,kc) = \sum_{\pi \in C_{ab,cd}} \gamma_\pi \sum_{\sigma \in S_m} \alpha^\pi_\sigma p_\sigma(j,l|i,k) = \sum_{\pi \in C_{ab,cd}, \sigma \in S_m} \gamma_\pi \alpha^\pi_\sigma p_\sigma(j,l|i,k).\]
Now for each $\sigma \in S_m$ and $\pi \in S_n$, let $\sigma \otimes \pi$ denote the permutation of $[m] \times [n]$ given by $ia \mapsto \sigma(i)\pi(a)$. Note that $p_\sigma(j,l|i,k) = p_{\sigma \otimes \pi}(jb,ld|ia,kc)$ for any $\pi$ such that $\pi(a) = b$ and $\pi(c) = d$. Furthermore, if either $\pi(a) \ne b$ or $\pi(c) \ne d$, then $p_{\sigma \otimes \pi}(jb,ld|ia,kc) = 0$. Therefore we have that
\begin{align*}
    p(jb,ld|ia,kc) &= \sum_{\pi \in C_{ab,cd}, \sigma \in S_m} \gamma_\pi \alpha^\pi_\sigma p_\sigma(j,l|i,k) \\
    &= \sum_{\pi \in C_{ab,cd}, \sigma \in S_m} \gamma_\pi \alpha^\pi_\sigma p_{\sigma \otimes \pi} (jb,ld|ia,kc) \\
    &= \sum_{\pi \in C, \sigma \in S_m} \gamma_\pi \alpha^\pi_\sigma p_{\sigma \otimes \pi} (jb,ld|ia,kc)
\end{align*}
Lastly, since $\sum_{\sigma \in S_m} \alpha^\pi_\sigma = 1$ for all $\pi \in C$ and $\sum_{\pi \in C} \gamma_\pi = \sum_{\pi \in S_n} \gamma_\pi = 1$, we have that
\[\sum_{\pi \in C, \sigma \in S_m} \gamma^\pi \alpha^\pi_\sigma = \sum_{\pi \in C} \gamma^\pi \sum_{\sigma \in S_m} \alpha^\pi_\sigma  = \sum_{\pi \in C} \gamma_\pi = 1.\]
Of course all these coefficients are nonnegative and so we have written $p$ as a convex combination of classical correlations and it is therefore classical. Thus the quantum permutation group $\mathbb{G} \times \mathbb{H}$ has no nonlocal symmetry.\qeds

We will also use the lemma above to obtain results on the nonlocal symmetry of cartesian products and tensor products of graphs. We start with the definition of those products. 

\begin{definition}
Let $G=(V(G), E(G))$, $H=(V(H),E(H))$ be finite graphs. Denote by $A_G\in \mathrm{M}_n(\C)$, $A_H\in \mathrm{M}_m(\C)$ the adjacency matrices of $G$ and $H$, respectively. We have the following graph products. 
\begin{itemize}
\item[(i)] The \emph{cartesian product} $G\square H$ is the graph with vertex set $V(G) \times V(H)$, where $(g_1,h_1)$ and $(g_2, h_2)$ are connected if and only if ($g_1=g_2$ and $h_1 \sim h_2$) or ($g_1 \sim g_2$ and $h_1=h_2$). For the adjacency matrix, we get 
\begin{align*}
A_{G \square H}=A_{G} \otimes \vec{1} +\vec{1}  \otimes A_{H}.
\end{align*}
\item[(ii)] The \emph{tensor product} $G\times H$ is the graph with vertex set $V(G) \times V(H)$, where $(g_1,h_1)$ and $(g_2, h_2)$ are adjacent if and only if $g_1 \sim g_2$ and $h_1 \sim h_2$. Therefore
\begin{align*}
A_{G \times H}=A_{G} \otimes A_{H}.
\end{align*}
\end{itemize}
\end{definition}

We also need the following proposition.

\begin{prop}[Proposition 4.1 in \cite{qauts11}]\label{propositionsurjhoms}
Let $G$ and $H$ be finite graphs. Denote by $U$, $V$, $W$ and $X$ the fundamental corepresentations of $\qut(G)$, $\qut(H)$ , $\qut(G \square H)$ and $\qut(G \times H)$, respectively. Then, we have the surjective $*$-homomorphisms 
\begin{itemize}
 \item[(i)]$\begin{aligned}[t]\varphi_1:C(\qut(G \square H)) &\to C(\qut(G))\otimes_{\mathrm{max}} C(\qut(H)),\\w_{ia,jb}&\mapsto u_{ij}v_{ab}\end{aligned}$
\item[(ii)]$\begin{aligned}[t]\varphi_2:C(\qut(G \times H) &\to C(\qut(G))\otimes_{\mathrm{max}} C(\qut(H)).\\x_{ia,jb}&\mapsto u_{ij}v_{ab}\end{aligned}$
\end{itemize}
\end{prop}

The next theorem shows that for certain graphs $G$ without nonlocal symmetry and graphs $H$ without quantum symmetry, their cartesian and/or tensor product will also have no nonlocal symmetry.

\begin{theorem}\label{thm:prodnononlocal}
Let $G$ be a graph with no nonlocal symmetry and $H$ a graph with no quantum symmetry. 
\begin{itemize}
    \item[(i)] If the map $\varphi_1$ from Proposition \ref{propositionsurjhoms} is an isomorphism, then $G\square H$ has no nonlocal symmetry.
    \item[(ii)]If the map $\varphi_2$ from Proposition \ref{propositionsurjhoms} is an isomorphism, then $G\times H$ has no nonlocal symmetry.
\end{itemize}
\end{theorem}
\proof
We only prove $(i)$, the proof of $(ii)$ is similar. Let $W$ be the fundamental representation of $\qut(G \square H)$, let $t$ be a tracial state on $C(\qut(G \square H))$, and let $p$ be the corresponding correlation, i.e.,
\begin{align*}
p(jb,ld|ia,kc)=t(w_{ia,jb}w_{kc,ld}).
\end{align*}
Using the isomorphism $\varphi_1$, we get 
\begin{align*}
p(jb,ld|ia,kc)=\tau(u_{ij}v_{ab}u_{kl}v_{cd})=\tau(u_{ij}u_{kl}v_{ab}v_{cd}),
\end{align*}
where $\tau=t \circ \varphi_1^{-1}$ is a tracial state on $C(\qut(G))\otimes_{\mathrm{max}} C(\qut(H))$. Since $G$ has no nonlocal symmetry and $H$ has no quantum symmetry, we get that $p$ is a classical correlation by Proposition \ref{prop:qpgtensor}.\qeds

The following result from~\cite{qauts11} gives sufficient conditions for the maps from Proposition~\ref{propositionsurjhoms} to be isomorphisms.

\begin{theorem}[Theorem 4.1 in \cite{qauts11}]\label{thm:spectra}
Let $G$ and $H$ be connected regular graphs. Let $\sigma_{G}=\{\lambda_i \, |\, i=1,\dots, n\}$ be the set of eigenvalues of $A_{G}$ and $\sigma_{H}=\{\mu_j \, |\, j=1,\dots,m\}$ be the set of distinct eigenvalues of $A_{H}$. Then we have
\begin{itemize}
\item[(i)] the map $\varphi_1$ from Proposition \ref{propositionsurjhoms}  is an isomorphism if $\{\lambda_i-\lambda_j \, |\, i,j=1,\dots, n\}\cap \{\mu_k - \mu_l\, |\, k,l=1,\dots, m\} = \{0\}$.
\item[(ii)] the map $\varphi_2$ from Proposition \ref{propositionsurjhoms} is an isomorphism if $\sigma_{G}$ and $\sigma_{H}$ do not contain $0$ and \linebreak$\left\{\frac{\lambda_i}{\lambda_j}\, |\, i,j=1,\dots, n\right\}\cap \left\{\frac{\mu_k}{\mu_l}\, |\, k,l=1,\dots, m\right\} = \{1\}$.
\end{itemize}
\end{theorem}

Combining this with Theorem~\ref{thm:prodnononlocal} we obtain the following corollary:

\begin{cor}\label{cor:nononlocalsymgen}
Let $G$ and $H$ be connected regular graphs such that $G$ has no nonlocal symmetry and $H$ has no quantum symmetry. Let $\sigma_{G}=\{\lambda_i \, |\, i=1,\dots, n\}$ be the set of distinct eigenvalues of $A_{G}$ and $\sigma_{H}=\{\mu_j \, |\, j=1,\dots,m\}$ be the set of distinct eigenvalues of $A_{H}$. Then we have
\begin{itemize}
\item[(i)] $G \square H$ has no nonlocal symmetry if $\{\lambda_i-\lambda_j \, |\, i,j=1,\dots, n\}\cap \{\mu_k - \mu_l\, |\, k,l=1,\dots, m\} = \{0\}$.
\item[(ii)] $G \times H$ has no nonlocal symmetry if $\sigma_{G}$ and $\sigma_{H}$ do not contain $0$ and \linebreak$\left\{\frac{\lambda_i}{\lambda_j}\, |\, i,j=1,\dots, n\right\}\cap \left\{\frac{\mu_k}{\mu_l}\, |\, k,l=1,\dots, m\right\} = \{1\}$.
\end{itemize}
\end{cor}
Applying the above to the special case $H = K_2$, we obtain the following:

\begin{cor}\label{cor:nononlocsym}
Let $G$ be a connected regular graph that has no nonlocal symmetry and let $\sigma_{G}=\{\lambda_i \, |\, i=1,\dots, n\}$ be the set of eigenvalues of $A_{G}$.
\begin{itemize}
    \item[(i)] If $\lambda_i-\lambda_j\notin \{-2,2\}$ for all $1\le i,j \le n$, then $G\square K_2$ has no nonlocal symmetry.
    \item[(ii)] If $\lambda_i\neq -\lambda_j$ for all $1\le i,j \le n$, then $G\times K_2$ has no nonlocal symmetry.
\end{itemize}
\end{cor}

\proof
This follows from Theorem \ref{thm:prodnononlocal}, Theorem \ref{thm:spectra} and the fact that $\sigma_{K_2}=\{-1,1\}$.
\qeds

In contrast to the previously discussed results, the next proposition deals with graphs that have nonlocal symmetry. We show that if one of the graphs involved in the cartesian - or tensor product has nonlocal symmetry, then also the product graph has nonlocal symmetry. 

\begin{prop}\label{prop:prodsym}
Let $G$ and $H$ be finite graphs. If $G$ or $H$ have nonlocal symmetry, then also $G\square H$ and $G \times H$.
\end{prop}

\proof
We prove that $G\square H$ has nonlocal symmetry, the proof for $G \times H$ is identical. Without loss of generality, let $G$ have nonlocal symmetry. Then, there exists a tracial state $\tau$ and a magic unitary $U$ commuting with $A_G$ such that the associated correlation $p(j,l\,|\,i,k)=\tau(u_{ij}u_{kl})$ is nonclassical. Consider the matrix $W$, where $w_{ia,jb}=u_{ij}\delta_{ab}$. It is easy to see that $W$ is a magic unitary that commutes with $A_{G\square H}$. We have the associated correlation
\begin{align}
    p'(jb,ld\,|\, ia,kc)=\tau(w_{ia,jb}w_{kc,ld})= \tau(u_{ij}\delta_{ab}u_{kl}\delta_{cd})=\delta_{ab}\delta_{cd}\tau(u_{ij}u_{kl}).\label{eq:corten}
\end{align}
Assume that $p'$ is classical, i.e., 
\begin{align}
p'=\sum_{\pi \in \aut(G\square H)} \alpha_\pi p_{\pi}\label{eq:corten2}
\end{align}
for some $0\le \alpha_\pi \le 1, \sum_{\pi \in \aut(G\square H)} \alpha_\pi=1$. By equation \eqref{eq:corten}, we know that $p'(jb,ld\,|\, ia,kc)=0$ for $a\neq b$ or $c\neq d$. Thus, using equation \eqref{eq:corten2}, we get $\pi(ia)=ja$ for all $ia \in V(G\square H)$ and $\pi \in \aut(G\square H)$ with $\alpha_\pi \neq 0$. 
Therefore, all automorphisms $\pi$ with $\alpha_\pi \neq 0$ induce an automorphism $\tilde{\pi} \in \aut(G)$, such that 
\begin{align*}
p'(jb,ld\,|\, ia,kc)=\delta_{ab}\delta_{cd}\sum_{\pi \in \aut(G\square H)} \alpha_\pi p_{\tilde{\pi}}(j,l|i,k).
\end{align*}
Then equation \eqref{eq:corten} implies $p(j,l\,|\, i,k)=\sum_{\pi \in \aut(G\square H)} \alpha_\pi p_{\tilde{\pi}}(j,l|i,k)$ contradicting our assumption that $p$ is nonclassical. Thus, $p'$ is nonclassical and $G\square H$ has nonlocal symmetry. 
\qeds

\section{Nonlocal symmetry of small graphs}\label{sec:smallgraphs}

In this section, we use the results from Sections \ref{sec:complete}--\ref{sec:products} to show whether or not specific small graphs have nonlocal symmetry. Note that by Section \ref{sec:complete}, the quantum group $S_n^+$ does not have nonlocal symmetry for $n \le 4$, which implies that all graphs on four or less points do not have nonlocal symmetry. Therefore, we study graphs on five vertices in Section \ref{sec:5vertices} and show that the only graphs on five vertices that have nonlocal symmetry are $K_5$ and $\overline{K_5}$. After that, we deal with all vertex-transitive graphs up to 11 vertices in Section \ref{sec:smalltrans}.

\subsection{Graphs on five vertices}\label{sec:5vertices}

We show that the only graphs on five vertices admitting nonlocal symmetry are $K_5$ and $\overline{K_5}$. We start with the following lemma. 
\begin{lemma}[Lemma 3.2.3 in \cite{fulton_quantum_nodate}]\label{lem:fulton}
Let $G$ be a finite graph, $A_G$ be its adjacency matrix and $u_{ij}$, $1\le i,j \le n$ be the generators of $C(\qut(G))$. If $(A_G^l)_{ii}\neq (A_G^l)_{jj}$ for some $l\in \N$, then $u_{ij}=0$. Especially, if $\deg(i)\neq \deg(j)$, then $u_{ij}=0$.
\end{lemma}

Furthermore, we have to look at quantum subgroups, which we define next.
\begin{definition}
Let $\mathbb{G}=(C(\mathbb{G}), U)$ and $\mathbb{H}=(C(\mathbb{H}), V)$ be compact matrix quantum groups. We say that $\mathbb{G}$ is a \emph{matrix quantum subgroup} of $\mathbb{H}$ if the map $\phi:C(\mathbb{H})\to C(\mathbb{G}), u_{ij} \mapsto v_{ij}$ is a $*$-homomorphism. We then write $\mathbb{G}\subset \mathbb{H}$.
\end{definition}

It is easy to see that if a quantum group has no nonlocal symmetry, then all its matrix quantum subgroups also do not have nonlocal symmetry. Using this, we see that $K_5$ and $\overline{K_5}$ are the only graphs on five vertices that have nonlocal symmetry. 

\begin{prop}
Besides $K_5$ and $\overline{K_5}$, no graph on five vertices has nonlocal symmetry. 
\end{prop}

\proof
Let $G$ be a graph on five vertices. If the graph is not regular, then we either have $\qut(G) \subset S_4^+$ or $\qut(G) \subset S_3 * \Z_2$ by Lemma \ref{lem:fulton}. By Theorem \ref{thm:Q4=C4} and Proposition \ref{prop:freeprodqpg}, the quantum groups $S_4^+$ and $S_3*\Z_2$ do not produce non-classical correlations and since $\qut(G)$ is a matrix quantum subgroup, this quantum group also does not produce any non-classical correlations. 

Now, we assume that $G$ is regular. The only regular graphs on five vertices are $C_5$, $K_5$ and $\overline{K_5}$. Since $C(\qut(C_5))$ is commutative (\cite{qauts11}), we get that $C_5$ has no nonlocal symmetry. Thus $K_5$ and $\overline{K_5}$ are the only graphs on five vertices that have nonlocal symmetry. 
\qeds

\subsection{Small transitive graphs}\label{sec:smalltrans}

Now, we discuss the nonlocal symmetry of all vertex-transitive graphs of order $6\le n \le 11$. Note that the quantum automorphism groups of these graphs were computed by Banica and Bichon in \cite{qauts11}. All graphs without quantum symmetry also have no nonlocal symmetry, so we do not include those graphs in the discussion. We also exclude the complete graphs from our next proposition, since we know from Section \ref{sec:complete} that the complete graphs $K_n$ do have nonlocal symmetry for $n \geq 5$.

\begin{prop}\label{prop:smallvtxtrans}
Let $G \neq K_n$ be a vertex-transitive graph of order $6 \le n\le 11$ that has quantum symmetry. Then we have Table~\ref{table1}. 
\begin{table}[ht]
\begin{center}
\begin{tabular}{llcllc}
\hline
&Name of $G$&Order&$\aut(G)$&$\qut(G)$&nonlocal symmetry\\[.5ex]
\hline
(1)&$2K_3$ &6&$S_3\wr \Z_2$&$S_3\wr_*\Z_2$&no\\
(2)&$3K_2$ &6&$ \Z_2 \wr S_3$& $  \Z_2\wr_*S_3$&yes\\
(3)&Cube $Q_3$ &8&$S_4 \times \Z_2$&$S_4^+ \times \Z_2$&no\\
(4)&$2K_4$&8&$S_4\wr \Z_2$&$S_4^+\wr_*\Z_2$&yes\\
(5)&$2C_4$&8&$(\Z_2 \wr \Z_2) \wr \Z_2$&$(\Z_2 \wr_*\Z_2)\wr_* \Z_2$&yes\\
(6)&$4K_2$&8&$\Z_2\wr S_4$&$\Z_2 \wr_* S_4^+$&yes\\
(7)&$3K_3$&9&$S_3\wr S_3$&$S_3 \wr_* S_3$&yes\\
(8)&$K_5 \square K_2$&10&$S_5\times \Z_2$&$S_5^+ \times \Z_2$&yes\\
(9)&$C_{10}(4)$&10&$\Z_2 \wr D_5$&$\Z_2 \wr_* D_5$&yes\\
(10)&$2C_5$&10&$D_5 \wr Z_2$&$D_5 \wr_* \Z_2$&no\\
(11)&$2K_5$&10&$S_5 \wr \Z_2$&$S_5^+ \wr_* \Z_2$&yes\\
(12)&$5K_2$&10&$\Z_2 \wr S_5$&$\Z_2 \wr_* S_5^+$&yes\\
\hline
\end{tabular}
\end{center}
\vspace{0.1cm}
\caption{Nonlocal symmetry of all vertex-transitive graphs of order $6 \le n\le 11$ that have quantum symmetry, excluding $K_n$.}\label{table1}
\end{table}
\end{prop}

\proof
The quantum automorphism groups of all those graphs are known from \cite{qauts11}. We prove for every row whether or not the given graph has nonlocal symmetry. 
\begin{itemize}
\item[(1)] We know that the graph $K_3$ has no quantum symmetry. By Proposition \ref{prop:2disjunion}, we get that $2K_3$ has no nonlocal symmetry.
\item[(2)] Consider for every $K_2$ the automorphism that exchanges the vertices in $K_2$. Those three automorphisms are non-trivial and disjoint, thus $3K_2$ has nonlocal symmetry by Theorem \ref{thm:threedisjaut}.
\item[(3)] The cube graph $Q_3$ is isomorphic to $K_4 \times K_2$. Since $K_4$ has eigenvalues $-1$ and $3$, we get that $Q_3$ has no nonlocal symmetry by Corollary \ref{cor:nononlocsym} $(ii)$.
\item[(4)] Consider one $K_4$ and choose two non-trivial, disjoint automorphisms, for example $(12)$ and $(34)$. Now, choosing any non-trivial automorphism of the other $K_4$ yields three non-trivial, disjoint automorphisms. We deduce that $2K_4$ has nonlocal symmetry from Theorem \ref{thm:threedisjaut}.
\item[(5)] Consider one $C_4$ and choose two non-trivial, disjoint automorphisms. Now, choosing any non-trivial automorphism of the other $C_4$ yields three non-trivial, disjoint automorphisms. We deduce that $2C_4$ has nonlocal symmetry from Theorem \ref{thm:threedisjaut}.
\item[(6)] Similar to $(2)$, we get that $4K_2$ has three non-trivial, disjoint automorphisms. Theorem \ref{thm:threedisjaut} yields that $4K_2$ has nonlocal symmetry. 
\item[(7)] Consider for every $K_3$ a non-trivial automorphism. Those three automorphisms are non-trivial and disjoint, thus $3K_3$ has nonlocal symmetry by Theorem \ref{thm:threedisjaut}.
\item[(8)] We know that $K_5$ has nonlocal symmetry by Section \ref{sec:K5}. Proposition \ref{prop:prodsym} shows that $K_5 \square K_2$ also has nonlocal symmetry.
\item[(9)] The graph $C_{10}(4)$ is the complement of $K_2 \circ C_5$, where $\circ$ denotes the lexicographic product (the vertex set is given by $V(K_2)\times V(C_5)$). Using this labelling, the transpositions $((1,1),(2,1))$, $((1,2),(2,2))$ and $((1,3),(2,3))$ are non-trivial, disjoint automorphisms of $C_{10}(4)$. Thus Theorem \ref{thm:threedisjaut} yields that $C_{10}(4)$ has nonlocal symmetry. 
\item[(10)] We know that $C_5$ has no quantum symmetry. Therefore $2C_5$ has no nonlocal symmetry by Proposition \ref{prop:2disjunion}.
\item[(11)] Similar to $(4)$, we get that $2K_5$ has three non-trivial, disjoint automorphisms. Using Theorem \ref{thm:threedisjaut}, we obtain that $2K_5$ has nonlocal symmetry. 
\item[(12)] Similar to $(2)$, we get that $5K_2$ has three non-trivial, disjoint automorphisms and thus nonlocal symmetry by Theorem \ref{thm:threedisjaut}.
\end{itemize}
\qeds

\section{Discussion}\label{sec:discuss}

In this work we have established the notion of nonlocal symmetry for graphs and quantum permutation groups which is motivated by physical observations, as opposed to quantum symmetry which is based on purely mathematical considerations. The primary goal of this work was to exhibit differences between these two notions, and we have seen several examples: the complete graph on four vertices does not admit nonlocal symmetry despite having quantum symmetry; the group based construction of Section~\ref{sec:groupconstruct} can produce non-classical quantum Latin squares that nevertheless give rise to local correlations; graphs with two disjoint automorphisms have quantum symmetry but not necessarily nonlocal symmetry; there are 10 \emph{connected} graphs on five vertices with quantum symmetry~\cite{smallgraphs}, but $K_5$ is the only one with nonlocal symmetry. Each of these examples give rise to a scenario in which non-commutativity of measurement operators does not translate into observable nonlocal effects. There is a natural, though somewhat imprecise, question to ask here. Our notion of nonlocality is based on the specific scenario of the isomorphism game. Given a graph $G$ with quantum symmetry, is there \emph{some} physical scenario that can be associated to $G$ (in a natural way) in which non-commutativity of $C(\qut(G))$ always results in non-classical behavior.


Though we have seen many examples where there is quantum but no nonlocal symmetry, we should note that the two notions seem to diverge less often than they coincide, at least based on what we have seen in this work. Indeed, nonlocal and quantum symmetry coincide for all complete graphs $K_n$ (i.e., for all of the quantum symmetric groups $S_n^+$) except when $n=4$. The proportion of non-classical quantum Latin squares produced by the construction of Section~\ref{sec:groupconstruct} that give rise to local correlations seems to diminish as we increase the size of cyclic group used. Though two disjoint automorphisms do not suffice for nonlocal symmetry, three disjoint automorphisms do suffice. Of the vertex transitive graphs on 11 or fewer vertices, quite few have quantum symmetry but not nonlocal symmetry. So though the notions are distinct, it is possible that they in general tend to coincide. We could formally ask the question: given a random graph $G$ on $n$ vertices, does the probability that $G$ has quantum symmetry but not nonlocal symmetry go to zero as $n$ goes to infinity? The answer is yes, but for a boring reason: the probability that $G$ has quantum symmetry at all goes to zero~\cite{qperms}. Thus it may be more interesting to consider restricting to $G$ that do have quantum symmetry. For this question we do not have an answer.

Perhaps the place where we saw the greatest difference between the two notions is in applying the construction of Section~\ref{sec:groupconstruct} to the group $\mathbb{Z}_2^3$. Here, none of the correlations produced by this construction were nonlocal, despite there being hundreds of non-classical quantum Latin squares. This is in stark contrast to the other abelian groups with 8 elements, especially $\mathbb{Z}_8$ which produces only 10 classical correlations out of 380. This of course raises the question of whether this is a more general phenomenon: does the construction of Section~\ref{sec:groupconstruct} produce only local correlations for the groups $\mathbb{Z}_2^d$? Unfortunately checking this for even $d=4$ appears to be currently out of reach. This is not the case for powers of other groups, since the construction does produce nonlocal correlations for $\mathbb{Z}_3^2$.

In addition to graphs, we have also defined nonlocal symmetry for quantum permutation groups. But what about more general compact matrix quantum groups $\mathbb{G}$? Here, non-commutativity of $C(\mathbb{G})$ is still used as the notion of quantumness within the quantum groups literature, but what can we use as our notion of nonlocality? One option is to consider functions $p(l,k|i,j) = \tau(u_{il}u_{jk})$ where $U = (u_{ij})$ is the fundamental representation of $\mathbb{G}$ and $\tau$ is a tracial state on $C(\mathbb{G})$. This is no longer a correlation, as the values need not even be nonnegative, but we can ask if there is such a function $p$ that cannot be obtained in the same way from the abelianization of $C(\mathbb{G})$, i.e., we impose the same relations on the generators $u_{ij}$ as for $C(\mathbb{G})$ with the additional constraint that they commute. If the answer is no, we may say that $\mathbb{G}$ exhibits nonlocality. This gives a well defined notion but perhaps not one deserving of that name, since it is not clear that there is a corresponding physical scenario where nonlocality, in the physics sense, actually occurs. It is an interesting question whether such a scenario can be constructed.

Many of the results of Section~\ref{sec:tensorproduct} are ripe for generalization. For instance, an optimistic goal would be to show that for any connected graphs $G$ and $H$ that do not have nonlocal symmetry and are quantumly coprime\footnote{Any connected graph has a unique prime factorization with respect to the cartesian product~\cite{sabidussiproduct,vizingproduct}, and two graphs are said to be coprime if they share no prime factors. By \emph{quantumly coprime}, we mean that no factor of one graph is quantum isomorphic to a factor of another. It is not {\it{a priori}} obvious that this definition remains unchanged if one only considers prime factors.}, the cartesian product of $G$ and $H$ does not have nonlocal symmetry. One approach to this would be to generalize Theorem~\ref{thm:prodnononlocal} so that the condition of $H$ having no quantum symmetry is replaced with simply having no nonlocal symmetry, and to generalize Theorem~\ref{thm:spectra} to the case where $G$ and $H$ are any two connected graphs that are quantumly coprime. The mentioned generalization of Theorem~\ref{thm:spectra} would be a quantum analog of a classical result of Sabidussi~\cite{sabidussiproduct}.

Another interesting direction to consider is whether one can find a (nice) description of $Q(G)$ for a graph $G$ with $L(G) \ne Q(G)$. An obvious choice would be to look at $Q(5)$. If this proves too difficult then an alternative would be to describe the set of $\mathbb{Z}_5$-invariant correlations in $Q(5)$. One difficulty that may arise here is that we do not expect these sets of quantum correlations to be polytopes when they differ from the classical set. That then raises the question of whether this expectation is always met, or can we find a graph $G$ such that $Q(G) \ne L(G)$ but $Q(G)$ is a polytope.

The notion of nonlocal symmetry that we have defined for quantum permutation groups very much depends of the specific presentation of that quantum group, and thus we do not expect that this property is invariant under isomorphism. In particular, it is possible for such isomorphisms to map generators (i.e., entries of the fundamental representation) to linear combinations of products of generators, and this seems like it could destroy (or create) nonlocal symmetry. However, we do not have an explicit example of two isomorphic quantum permutation groups where one has nonlocal symmetry and the other does not. Thus this remains an open question.

Lastly, we remark that there are some special graphs that we are interested in knowing whether they have nonlocal symmetry. To begin, do the higher order cubes, $Q_d$ for $d \ge 4$, have nonlocal symmetry? All cubes other than $Q_1 = K_2$ have quantum symmetry, but we have seen that both $Q_2$ and $Q_3$ do not have nonlocal symmetry. Does this pattern continue or does nonlocal symmetry eventually emerge? Note that once it occurs for $Q_d$ is occurs for $Q_{d'}$ for all $d' \ge d$ by Proposition~\ref{prop:prodsym} since $Q_{d+1} = Q_d \square K_2$. Additionally, what about the folded cubes $FQ_d$. These are obtained from $Q_d$ by identifying pairs of vertices at maximum Hamming distance, or by adding edges between the vertices at maximum distance in $Q_{d-1}$. They are known to have quantum symmetry for odd $d \ge 3$~\cite{foldedcubes}. For small cases we have $FQ_3 \cong K_4$ which has no nonlocal symmetry by Theorem~\ref{thm:Q4=C4}, and $FQ_4 \cong K_{4,4} \cong \overline{2K_4}$ which has nonlocal symmetry by Theorem~\ref{prop:smallvtxtrans}. Perhaps the folded cube $FQ_d$ continues to have nonlocal symmetry for all $d \ge 4$. However, for folded cubes, there are many aspects in which the even and odd cases differ, and so another possibility is that $FQ_d$ has nonlocal symmetry for even $d \ge 4$ but does not for odd $d \ge 5$. Currently, we do not know whether $FQ_5$, better known as the Clebsch graph, has nonlocal symmetry. We remark that by the techniques of~\cite{qperms} it is possible to show that any quantum automorphism of the cube $Q_{d-1}$ is also a quantum automorphism of $FQ_d$. Therefore, if $FQ_d$ does not have nonlocal symmetry then neither does $Q_{d-1}$. Furthermore, for odd $d$ it is known that $FQ_d \times K_2 \cong Q_d$, and moreover the conditions of Corollary~\ref{cor:nononlocsym} are met in this case. Thus for odd $d$, if $FQ_d$ does not have nonlocal symmetry then neither do $Q_{d-1}$ or $Q_d$. But if $FQ_d$ does have nonlocal symmetry for some odd $d$, then so does $Q_{d'}$ for all $d' \ge d$ by Proposition~\ref{prop:prodsym}. 

\subsection*{Acknowledgments}

\noindent\textbf{DR:} I would like to thank Moritz Weber for inviting me to visit him and his group at Saarland University where this collaboration began. I would also like to acknowledge Robert \v{S}\'{a}mal with whom I discussed and performed some of the computations that appear in Section~\ref{sec:groupconstruct}. The research leading to these results has received funding from the European Union’s Horizon 2020 research and innovation program under the Marie Sklodowska-Curie grant agreement no. 713683 (COFUNDfellowsDTU).\\
\textbf{SS:} Supported by the EPSRC project ``Quantum groups in action'', EP/T03064X/1 and the DFG project ``Quantenautomorphismen von Graphen''.

\bibliographystyle{plainurl}

\bibliography{main.bbl}

\begin{thebibliography}{10}

\bibitem{qiso1}
Albert Atserias, Laura Man\v{c}inska, David~E. Roberson, Robert \v{S}\'{a}mal,
  Simone Severini, and Antonios Varvitsiotis.
\newblock Quantum and non-signalling graph isomorphisms.
\newblock {\em Journal of Combinatorial Theory, Series B}, 136:289 -- 328,
  2019.
\newblock \href {http://dx.doi.org/10.1016/j.jctb.2018.11.002}
  {\path{doi:10.1016/j.jctb.2018.11.002}}.

\bibitem{banicahomogeneous}
Teodor Banica.
\newblock Quantum automorphism groups of homogeneous graphs.
\newblock {\em Journal of Functional Analysis}, 224(2):243 -- 280, 2005.
\newblock \href {http://dx.doi.org/10.1016/j.jfa.2004.11.002}
  {\path{doi:10.1016/j.jfa.2004.11.002}}.

\bibitem{BanBicfr}
Teodor Banica and Julien Bichon.
\newblock Free product formulae for quantum permutation groups.
\newblock {\em J. Inst. Math. Jussieu}, 6(3):381--414, 2007.
\newblock \href {http://dx.doi.org/10.1017/S1474748007000072}
  {\path{doi:10.1017/S1474748007000072}}.

\bibitem{qauts11}
Teodor Banica and Julien Bichon.
\newblock Quantum automorphism groups of vertex-transitive graphs of order
  {$\le11$}.
\newblock {\em Journal of Algebraic Combinatorics}, 26(1):83, 2007.
\newblock \href {http://dx.doi.org/10.1007/s10801-006-0049-9}
  {\path{doi:10.1007/s10801-006-0049-9}}.

\bibitem{Bicfreewreath}
Julien Bichon.
\newblock Free wreath product by the quantum permutation group.
\newblock {\em Algebr. Represent. Theory}, 7(4):343--362, 2004.
\newblock \href {http://dx.doi.org/10.1023/B:ALGE.0000042148.97035.ca}
  {\path{doi:10.1023/B:ALGE.0000042148.97035.ca}}.

\bibitem{Blackadar}
Bruce Blackadar.
\newblock {\em Operator algebras}, volume 122 of {\em Encyclopaedia of
  Mathematical Sciences}.
\newblock Springer-Verlag, Berlin, 2006.
\newblock Theory of $C^*$-algebras and von Neumann algebras, Operator Algebras
  and Non-commutative Geometry, III.
\newblock \href {http://dx.doi.org/10.1007/3-540-28517-2}
  {\path{doi:10.1007/3-540-28517-2}}.

\bibitem{smallgraphs}
Christian Eder, Viktor Levandovskyy, Julien Schanz, Simon Schmidt, Andreas
  Steenpass, and Moritz Weber.
\newblock Existence of quantum symmetries for graphs on up to seven vertices: a
  computer based approach.
\newblock 2019.
\newblock \href {http://arxiv.org/abs/1906.12097} {\path{arXiv:1906.12097}}.

\bibitem{fulton_quantum_nodate}
Melanie Fulton.
\newblock {\em The {quantum} {automorphism} {group} and {undirected} {trees}}.
\newblock PhD thesis, Virginia Polytechnic Institute and State University,
  2006.

\bibitem{qperms}
Martino Lupini, Laura Mančinska, and David~E. Roberson.
\newblock Nonlocal games and quantum permutation groups.
\newblock {\em Journal of Functional Analysis}, 279(5):108592, 2020.
\newblock \href {http://dx.doi.org/10.1016/j.jfa.2020.108592}
  {\path{doi:10.1016/j.jfa.2020.108592}}.

\bibitem{qfuncs}
Benjamin Musto, David Reutter, and Dominic Verdon.
\newblock A compositional approach to quantum functions.
\newblock {\em Journal of Mathematical Physics}, 59(8):081706, 2018.
\newblock \href {http://dx.doi.org/10.1063/1.5020566}
  {\path{doi:10.1063/1.5020566}}.

\bibitem{qlatin}
Benjamin Musto and Jamie Vicary.
\newblock Quantum {L}atin squares and unitary error bases.
\newblock {\em Quantum Information \& Computation}, 16(15-16):1318--1332, 2016.
\newblock \href {http://dx.doi.org/10.26421/QIC16.15-16}
  {\path{doi:10.26421/QIC16.15-16}}.

\bibitem{bisynchronous}
Vern~I. Paulsen and Mizanur Rahaman.
\newblock Bisynchronous games and factorizable maps.
\newblock 2019.
\newblock \href {http://arxiv.org/abs/1908.03842} {\path{arXiv:1908.03842}}.

\bibitem{paulsen2016estimating}
Vern~I. Paulsen, Simone Severini, Daniel Stahlke, Ivan~G Todorov, and Andreas
  Winter.
\newblock Estimating quantum chromatic numbers.
\newblock {\em Journal of Functional Analysis}, 270(6):2188--2222, 2016.
\newblock \href {http://dx.doi.org/10.1016/j.jfa.2016.01.010}
  {\path{doi:10.1016/j.jfa.2016.01.010}}.

\bibitem{sabidussiproduct}
Gert Sabidussi.
\newblock Graph multiplication.
\newblock {\em Mathematische Zeitschrift}, 72(1):446--457, 1959.
\newblock \href {http://dx.doi.org/10.1007/BF01162967}
  {\path{doi:10.1007/BF01162967}}.

\bibitem{foldedcubes}
Simon Schmidt.
\newblock Quantum automorphisms of folded cube graphs.
\newblock {\em Annales de l'Institut Fourier}, 70(3):949--970, 2020.
\newblock \href {http://dx.doi.org/10.5802/aif.3328}
  {\path{doi:10.5802/aif.3328}}.

\bibitem{vizingproduct}
Vladim~G Vizing.
\newblock The cartesian product of graphs.
\newblock {\em Vycisl. Sistemy}, 9(30-43):33, 1963.

\bibitem{Wafree}
Shuzhou Wang.
\newblock Free products of compact quantum groups.
\newblock {\em Communications in Mathematical Physics}, 167(3):671--692, 1995.
\newblock \href {http://dx.doi.org/10.1007/BF02101540}
  {\path{doi:10.1007/BF02101540}}.

\bibitem{Watensor}
Shuzhou Wang.
\newblock Tensor products and crossed products of compact quantum groups.
\newblock {\em Proc. London Math. Soc. (3)}, 71(3):695--720, 1995.
\newblock \href {http://dx.doi.org/10.1112/plms/s3-71.3.695}
  {\path{doi:10.1112/plms/s3-71.3.695}}.

\bibitem{wang}
Shuzhou Wang.
\newblock Quantum symmetry groups of finite spaces.
\newblock {\em Communications in Mathematical Physics}, 195(1):195--211, Jul
  1998.
\newblock \href {http://dx.doi.org/10.1007/s002200050385}
  {\path{doi:10.1007/s002200050385}}.

\bibitem{woronowicz}
Stanis{\l}aw~Lech Woronowicz.
\newblock Compact matrix pseudogroups.
\newblock {\em Communications in Mathematical Physics}, 111(4):613--665, 1987.
\newblock \href {http://dx.doi.org/10.1007/BF01219077}
  {\path{doi:10.1007/BF01219077}}.

\end{thebibliography}

\end{document}